\documentclass[11pt,reqno]{amsart}
\usepackage{amsmath,amssymb,mathtools}
\usepackage{booktabs,array,longtable}
\usepackage{enumitem}
\usepackage{xcolor}
\usepackage[colorlinks=true,linkcolor=blue!60!black,citecolor=blue!60!black,urlcolor=blue!60!black]{hyperref}
\usepackage[capitalise,noabbrev]{cleveref}
\calclayout
\setlist[enumerate]{label=\textup{(\roman*)},leftmargin=*,itemsep=3pt}
\setlist[itemize]{leftmargin=*,itemsep=3pt}
\numberwithin{equation}{section}
\newtheorem{theorem}{Theorem}[section]
\newtheorem{proposition}[theorem]{Proposition}
\newtheorem{lemma}[theorem]{Lemma}
\newtheorem{corollary}[theorem]{Corollary}
\theoremstyle{definition}

\newtheorem{example}[theorem]{Example}

\theoremstyle{remark}
\newtheorem{remark}[theorem]{Remark}
\newcommand{\Qp}{\mathbf Q_p}
\newcommand{\Zp}{\mathbf Z_p}
\newcommand{\Fp}{\mathbf F_p}
\newcommand{\OK}{\mathcal O_K}
\newcommand{\GK}{G_K}
\newcommand{\Cp}{\mathbf C_p}
\newcommand{\Bcris}{\mathbf B_{\mathrm{cris}}}
\newcommand{\Bst}{\mathbf B_{\mathrm{st}}}
\newcommand{\BdR}{\mathbf B_{\mathrm{dR}}}
\newcommand{\Rplus}{\widetilde{\mathbf B}_{\mathrm{rig}}^+}
\newcommand{\Rlog}{\widetilde{\mathbf B}_{\log}^+}
\newcommand{\Dst}{D_{\mathrm{st}}}
\newcommand{\Dcris}{D_{\mathrm{cris}}}
\newcommand{\DdR}{D_{\mathrm{dR}}}
\newcommand{\Fil}{\operatorname{Fil}}
\newcommand{\Hst}{H^1_{\mathrm{st}}}
\newcommand{\Hf}{H^1_f}
\newcommand{\He}{H^1_e}
\newcommand{\Hg}{H^1_g}
\newcommand{\Hom}{\operatorname{Hom}}
\newcommand{\End}{\operatorname{End}}
\newcommand{\im}{\operatorname{im}}
\newcommand{\coker}{\operatorname{coker}}
\newcommand{\Tr}{\operatorname{Tr}}
\newcommand{\tr}{\operatorname{tr}}
\newcommand{\ord}{\operatorname{ord}}
\newcommand{\inv}{\operatorname{inv}}
\newcommand{\ev}{\operatorname{ev}}
\newcommand{\pgl}{\mathfrak{pgl}}
\newcommand{\slie}{\mathfrak{sl}}
\newcommand{\ur}{\mathrm{ur}}
\newcommand{\Tate}{\mathrm{Tate}}
\newcommand{\dd}{\mathrm d}
\newcommand{\pair}[2]{\left\langle #1,#2\right\rangle_K}
\title[Monodromy and local duality for $1$-motives]{Semistable $1$-motives: monodromy, explicit periods, and local duality}
\author{Mohammadreza Mohajer}
\address{Department of Mathematics, Physics and Geology, Cape Breton University,
1250 Grand Lake Road, Sydney, NS B1P 6L2, Canada}
\email{mohammadreza.mohajer96@gmail.com}
\email{Mohammadreza\_Mohajer@cbu.ca}
\begin{document}
\begin{abstract}
Let $K/\Qp$ be unramified, with $p>2$. For a $1$-motive with constant lattice, split torus, and good abelian quotient, we compute the semistable extension complex from toric valuations and formal-group logarithms. A canonical monodromy boundary exists for every semistable representation; it detects splitting over crystalline periods and the increase in monodromy rank. For the stated $1$-motives, this boundary and the unramified quotient recover the restricted local Tate pairing, with exponential radicals. Nonsplit elliptic examples, Raynaud uniformization, and a height-one positive-period bound explain the geometric range. For toric motives we construct actual periods, convergent Robba bases, and integral action matrices. A toric basis family has perfect projective/trace-zero matrix duality, including its exact integral defect. An explicit residue-zero solution operator constructs every dual Herr class. Applied componentwise, the same operator gives relative full dual lifts for general constant-lattice $1$-motives once their coefficient sequence is specified. Finally, finite-base coefficient traces realize division-tower Iwasawa cohomology for any Tate lattice, and the period cup pairing equals the arithmetic group-ring pairing.
\end{abstract}
\maketitle
{\tableofcontents}
\clearpage

\section{Introduction}\label{sec:introduction}

A $1$-motive combines a lattice, a torus, and an abelian variety. Its semistable realization records how these three pieces interact through Frobenius, monodromy, and the Hodge filtration. Local Tate duality pairs the Galois extension classes of that realization with those of its Cartier dual. We study this pairing by first separating the monodromy obstruction from the exponential classes, and then constructing period and cochain representatives wherever the coefficients can be written explicitly.

The aim of this paper is to develop explicit counterparts for
semistable $1$-motives of three classical constructions in local
$p$-adic arithmetic: 
\begin{enumerate}
    \item the valuation and logarithm description of
Kummer classes,
\item local Tate duality with Bloch--Kato local
conditions,
\item and the cyclotomic description of Iwasawa cohomology
by period modules.
\end{enumerate}
We extend the valuation--logarithm calculation to motives with
a good-reduction abelian part and describe how monodromy governs
semistable extension classes and their local duality.
We complement these formulas with explicit period and cochain
representatives in the toric families treated below, and carry
the cyclotomic cohomological and pairing constructions over
finite base fields to the torsion-and-division towers of
$1$-motives.

The geometric setting is
\begin{equation}\label{eq:intro-mixed-motive}
M=[u:\mathbf Z^s\to G],\qquad
0\longrightarrow\mathbf G_m^a\longrightarrow G\longrightarrow A\longrightarrow0,
\end{equation}
over an unramified extension $K/\Qp$, with $p>2$ and $n=[K:\Qp]$. Let $\OK$ be the ring of integers of $K$. We assume that $G$ has a semi-abelian model $\mathcal G$ over $\OK$, given by an extension of an abelian scheme of dimension $g$ by a split torus. The defining points $u(y)$ may be nonintegral. Write $V_M=T_p(M)\otimes\Qp$ for the covariant Tate realization. The hypotheses allow a nonsplit semi-abelian group and nonzero monodromy. They are stated precisely in \eqref{eq:mixed-motive}--\eqref{eq:good-model}.

There are two computable maps attached to these points. The quotient $G(K)/\mathcal G(\OK)$ is $\mathbf Z^a$, and its valuation map gives $R:\mathbf Z^s\to\mathbf Z^a$. Removing the toric point $p^{Ry}$ leaves an integral point $u_0(y)$. Put
\begin{equation}\label{eq:intro-geometric-coordinates}
\mathfrak g=\operatorname{Lie}(G),\qquad
\lambda(y)=-\log_{\mathcal G}(u_0(y)),
\end{equation}
with $\lambda$ extended $K$-linearly. \Cref{thm:mixedrealization} identifies $R$ with Fontaine monodromy and $\lambda$ with the Hodge quotient of the canonical Frobenius lattice lift. The proof uses Raynaud's good-point decomposition, the crystalline comparison for $p$-divisible groups, and the formal-group interpretation of Bloch--Kato's exponential. The filtered abelian module remains part of the realization. The two point maps determine the reduced extension calculation.

In the following rational cohomology formulas, $R$ denotes
its $\Qp$-linear extension $\Qp^s\to\Qp^a$; kernels and
cokernels of $R$ are taken over $\Qp$. Write $\Hst(K,V)$ for the classes in $H^1(K,V)$ with semistable middle representation. A single reduction lemma removes the crystalline abelian block and gives
\begin{equation}\label{eq:intro-mixed-normal}
0\longrightarrow\frac{\Qp^a\oplus\mathfrak g}{(R,\lambda)\Qp^s}
\longrightarrow\Hst(K,V_M)\longrightarrow\ker R\longrightarrow0.
\end{equation}
The final map takes the trace of the lattice Frobenius entry.

Let $d$ denote the rank of $R$, and put $h=\dim H^0(K,V_M)$. Then
\[
\dim\Hst(K,V_M)=n(a+g)+a-d+h.
\]
\Cref{lem:threeblock,thm:threeblockcohomology,thm:mixednormalform} give the maps, normal forms, and proof. Their specialization to the split toric motive
\begin{equation}\label{eq:motive-intro}
M_Q=[\mathbf Z^s\longrightarrow\mathbf G_m^a],\qquad
Q=(q_{ij})_{1\leq i\leq a,\,1\leq j\leq s},
\end{equation}
has $R=(\ord_p q_{ij})$, $L=(\log_p q_{ij})$, and $\lambda=-L$, where $\log_p(p)=0$. Thus \Cref{thm:normalform} gives
\begin{equation}\label{eq:intro-normal}
0\longrightarrow\frac{\Qp^a\oplus K^a}{(R,-L)\Qp^s}
\longrightarrow\Hst(K,V_Q)\longrightarrow\ker R\longrightarrow0.
\end{equation}

The monodromy boundary has a representation-theoretic formulation independent of this geometric setting. For any semistable $V$, let $D=\Dst(V)$ be its semistable period module, with Frobenius $\varphi$ and monodromy $N$, and put
$\mathcal M_N(V)=(D/ND)^{p\bar\varphi=1}$, where $\bar\varphi$ is the induced Frobenius. \Cref{thm:universalboundary} proves
\begin{equation}\label{eq:intro-universal-boundary}
0\longrightarrow\Hf(K,V)\longrightarrow\Hst(K,V)
\xrightarrow{\partial_N}\mathcal M_N(V)\longrightarrow0.
\end{equation}
Here $\Hf$ is the Bloch--Kato finite subgroup. The boundary of an extension is the new monodromy vector modulo $ND$. It is nonzero exactly when the rank of monodromy increases by one. \Cref{prop:obstruction-duality} identifies the dual of this obstruction space with the crystalline Frobenius quotient of the Tate dual.

Let $\He$ denote the Bloch--Kato exponential subgroup. For every semistable representation, \Cref{cor:universalquotientduality} shows that the restricted Tate pairing becomes perfect after dividing both sides by $\He$. This consequence of Berger's extension theorem and Bloch--Kato orthogonality connects the universal monodromy sequence with arithmetic duality.

For \eqref{eq:intro-mixed-motive}, the target of \eqref{eq:intro-universal-boundary} is $\coker R$. The unramified quotient character defines $\alpha_{V_M}:\Hst\to\ker R$. \Cref{thm:boundary} proves
\begin{equation}\label{eq:intro-boundary}
\He(K,V_M)\simeq\mathfrak g/\lambda(\ker R),\qquad
\Hst(K,V_M)/\He(K,V_M)\simeq\ker R\oplus\coker R.
\end{equation}
For $V_M^\dagger=V_M^\vee(1)$, use dual evaluation frames, so $R^\dagger=-R^{\mathsf t}$. \Cref{thm:semistablereciprocity} recovers the restriction of the local Tate pairing by
\begin{equation}\label{eq:intro-mixed-reciprocity}
\langle x,y\rangle_K
=\alpha_{V_M^\dagger}(y)(\partial_Nx)
 -\alpha_{V_M}(x)(\partial_{N^\dagger}y).
\end{equation}
Its radicals are exactly the exponential subgroups. The good abelian part contributes through its tangent space and logarithms, while the displayed quotient pairing uses the outer weight blocks. \Cref{ex:goodabelian,ex:mixedelliptic} compute the pure abelian case and a nonsplit extension of an elliptic curve by $\mathbf G_m$. \Cref{prop:uniformizationtransfer} explains the extension through Raynaud uniformization and unramified descent.

For the toric motive, \Cref{thm:actualperiods} constructs actual Galois matrices and invariant period lifts for every class in \eqref{eq:intro-normal}. The construction uses an unramified Witt-vector period and crystalline logarithms with prescribed de Rham values. It verifies a sharp one-pole bound in the positive logarithmic ring. The same ambient bound holds for the general mixed motive and all its semistable extensions by the positive comparison, as shown in \Cref{prop:generalpositive}. Writing $t$ for the Fontaine's cyclotomic period, this gives $t\Dcris(V_M)\subset\widetilde{\mathbf B}_{\mathrm{rig}}^+\otimes_{\Qp}V_M$. For a class $x\in\Hst(K,V_M)$, represented by a semistable extension
$0\to V_M\to E_x\to\Qp\to0$, choose
$h\in\Dst(E_x)$ lifting $1$.
\Cref{prop:boundarypositive} gives
\[
th=A_h-\ell_p\,tNh,
\qquad
A_h\in\widetilde{\mathbf B}_{\mathrm{rig}}^+
\otimes_{\Qp}E_x,
\]
where $\ell_p$ is the logarithmic period normalized by
$N(\ell_p)=-1$.
The boundary $\partial_Nx$ is the class of $Nh$ in
$\Dst(V_M)/N\Dst(V_M)$.
Consequently, changing the lift $h$ removes the logarithmic term
precisely when $x$ belongs to the Bloch--Kato finite subgroup
$\Hf(K,V_M)$.
The explicit toric bases over the cyclotomic Robba ring are
constructed in \Cref{thm:Robbabasis}.

For a finite free module $W$ over a commutative ring $O$, write $\pgl_O(W)=\End_O(W)/O I$ and $\slie_O(W)=\ker(\tr:\End_O(W)\to O)$. The full integral calculation concerns a toric basis family. Let $r=n+1$ and choose $q_0=p,q_1,\ldots,q_n$ forming a $\Zp$-basis of the pro-$p$ completion $W$ of $K^\times$. For $T=T_p([\mathbf Z\to\mathbf G_m^r])$ with this column of parameters, \Cref{thm:integraltrace} gives
\begin{equation}\label{eq:introtrace}
H^1(K,T)=\pgl_{\Zp}(W),\qquad
H^1(K,T^\dagger)=\slie_{\Zp}(W),\qquad
\langle[A],B\rangle_K=\tr(BA).
\end{equation}
Both lattices have rank $r^2-1$. Their duality remains perfect when $p\mid r$. The natural map $\slie_{\Zp}(W)\to\pgl_{\Zp}(W)$ has cokernel $\Zp/r\Zp$, which measures the obstruction to choosing integral trace-zero representatives; see
\Cref{prop:tracedefect}. The abelian product extension in \Cref{cor:integralproduct} retains this exact defect and adds the usual abelian Tate pairing under explicit torsion hypotheses.

Full dual cohomology requires more classes than the semistable complex supplies. \Cref{cor:dimensionobstruction} computes the missing dimensions in the basis family. \Cref{lem:unitinterpolation,lem:valuationinterpolation,thm:Robbabasis} construct the toric Robba module from convergent cyclotomic interpolation. \Cref{prop:integralRobba} gives its integral action matrices and trace operator. \Cref{thm:explicitprimitive} solves every Tate-valued degree-two equation with zero trace-residue by Laurent integration and convergent resolvents. It gives all dual Herr classes and the full residue identity
\begin{equation}\label{eq:intro-residue}
\operatorname{Inv}_{K,\gamma}(X_A\cup Y_B)=\tr(BA)
\end{equation}
in \Cref{thm:fulldualrepresentatives,thm:fullresiduereciprocity}, including the nonsemistable dual classes.

There is a relative extension of this construction to general $1$-motives. With $U=V_p(G)$, the obstruction for $z\in H^1(K,U^\dagger)$ to lift to $H^1(K,V_M^\dagger)$ is the vector
\[
y\longmapsto\inv_K(\kappa_G(u(y))\cup z).
\]
\Cref{prop:fullmixedobstruction} proves this criterion and computes its rank. Given the actual coefficient sequence for $U$ and $V_M$, \Cref{thm:relativeHerrlift} constructs the lift by applying the same Tate-valued solution operator to each lattice coordinate. The trace-zero equation of \eqref{eq:introtrace} is its toric basis-family specialization. This relative theorem specifies its coefficient inputs; the maps $(R,\lambda)$ alone do not determine full abelian period matrices.

Finally, let $\mathcal L/K$ be a Galois extension containing the cyclotomic tower, and let $F_m/K$ run through cofinal finite Galois subextensions. For any finite free Galois lattice $T$, \Cref{thm:towerperiodrealization} proves
\begin{equation}\label{eq:intro-tower}
H^1_{\mathrm{Iw}}(\mathcal L/K,T)
\simeq\varprojlim_m D_{F_m}(T)^{\psi=1},
\end{equation}
with coefficient-trace and cyclotomic-shift transitions. \Cref{thm:periodtowerreciprocity} identifies its explicit period cup pairing with the group-ring arithmetic pairing. \Cref{cor:motivicdivisiontower} applies both statements to the torsion-and-division tower of every $1$-motive. The toric Lie-dimension and orientation formulas require their own tower hypotheses; they are not reused for an abelian torsion tower. The coefficient inputs, integral qualifications, and further analytic comparisons are collected in \Cref{subsec:mixedscope,sec:remaining}. Appendix~\ref{app:checks} collects algebraic checks of the explicit formulas
and records the sign and normalization conventions used in the
cohomology, period, and duality calculations.

\subsection{Relation with existing theory}\label{subsec:related}

Deligne's $1$-motives and Cartier duality are used through the precise constructions in Harari--Szamuely \cite[Section 1]{HS05}. Their local duality theorem concerns the hypercohomology of the $1$-motive; the degree-one pairings here are the local Tate pairings of its rational and integral realizations. Raynaud gives the geometric monodromy, uniformizer decomposition, and strictification preserving Tate realizations \cite[Sections 4.2--4.5]{Ray94}. W\"urthen--Zhao construct the associated logarithmic finite groups and $p$-divisible groups \cite{WZ24}. The equivalence with semistable representations of Hodge height one and the compatibility of their monodromies are established by Bertapelle--Wang--Zhao \cite[Theorems B and F]{BWZ26}. Our coordinate calculations use these established geometric foundations and record their role at the proof steps.

Fontaine's rings and Kummer logarithms underlie the period representatives \cite[Chapter 6]{FO}. The semistable extension complex is classical; \cite[Theorem 1.1.4]{BLZ16} records the Bloch--Kato/Nekov\'a\v r Ext comparison. Bloch--Kato local conditions and orthogonality retain their original attribution \cite[Section 3]{BK90}. The three-block reduction makes these inputs explicit for mixed motives and specializes to the toric normal form without a second proof.

For semistable $V$, the equality
$\Hst(K,V)=\Hg(K,V)$ is due to Hyodo in the local-field
setting; we use Berger's extension theorem
\cite[Theorem 6.2]{Ber02}. Berger's semistable comparison and filtered reconstruction supply the coefficient-module comparison \cite{Ber02,Ber08}. Colmez's course \cite[Sections 8.3.3--8.4.3]{Col04} places the positive ring $\Rplus$ beside overconvergent modules and reciprocity; \Cref{prop:generalpositive} explains its role throughout the generalization. The interpolation formulas have rank-one antecedents in \cite[Sections 2.6--2.7]{Col08}. The cyclotomic selectors and logarithmic averaging inputs are explicitly attributed to Chitrao--Ghate--Yasuda and Benois at their uses \cite[Section 3]{CGY25}. Cherbonnier--Colmez's inverse theorem supplies the integral $\psi=0$ step \cite{CC99}. The residue coordinates agree with \cite[Corollary 1.5.5]{Ben11}.

Herr and Liu establish the full coefficient cohomology, cup comparison, and local duality \cite{Her01,Liu08}. We use those comparisons to identify the lifting obstruction, and then solve its scalar equations explicitly. The relative construction of \Cref{thm:relativeHerrlift} uses the same solution operator with a general quotient coefficient module. This is why the scalar residue correction extends while the projective/trace-zero description remains a property of the stated basis family.

Perrin-Riou's semistable cyclotomic regulator and its reciprocity theorem already provide an analytic theory \cite[Section 5.3]{PR01}; the crystalline and de Rham settings are treated in \cite{PR94,Nak14}. Riedel's thesis is a further source for the exponential and reciprocity formalism \cite{Rie13}. Formal-group, Lubin--Tate, and locally analytic extensions have substantial foundations in \cite{TR11,SV26,Poy25,Por24}. Gao--Zhao identify the appropriate cohomological comparisons for $(\varphi,\tau)$-modules and limitations of several replacement complexes \cite{GZ25}. Our division-tower limit retains Fontaine's cyclotomic $\psi=1$ theorem, as proved in \cite[Theorem II.1.3]{CC99}, at each finite base. It identifies transition traces before passing to the limit. The analytic-regulator and equal-characteristic comparisons in \Cref{sec:remaining} are stated with their additional requirements.

\section{Arithmetic conventions and motivic realizations}\label{sec:conventions}

\subsection{The multiplicative group and reciprocity}

Throughout, $K/\Qp$ is unramified of degree $n$, $p>2$, $\sigma$ is the arithmetic Frobenius of $K/\Qp$, and $\ord_p(p)=1$. Fix an algebraic closure $\overline K$ and put $G_K=\operatorname{Gal}(\overline K/K)$. We write $H^i(K,-)$ for continuous Galois cohomology. Define
\[
W=\varprojlim_m K^\times/(K^\times)^{p^m},\qquad W_{\Qp}=W\otimes_{\Zp}\Qp.
\]
We write $W$ additively. The multiplicative decomposition of $K^\times$ gives isomorphisms
\begin{equation}\label{eq:val-log}
W\xrightarrow{\ (\ord_p,\log_p)\ }\Zp\oplus p\OK,
\qquad W_{\Qp}\simeq\Qp\oplus K.
\end{equation}
Indeed, prime-to-$p$ torsion disappears in the completion and logarithm is an isomorphism $1+p\OK\to p\OK$. An unramified extension of $\Qp$ contains no nontrivial $p$-power root of unity.

For $O=\Zp$ the local cohomology groups needed later are
\begin{equation}\label{eq:standardcohom}
\begin{array}{c|ccc}
 & H^0(K,-)&H^1(K,-)&H^2(K,-)\\ \hline
 O&O&W^*&0\\
 O(1)&0&W&O,
\end{array}
\qquad W^*=\Hom_O(W,O).
\end{equation}
We use continuous cohomology for compact modules. Kummer theory identifies $H^1(K,O(1))$ with $W$, and local class field theory identifies $H^1(K,O)$ with $W^*$. We fix the sign of the latter identification by
\begin{equation}\label{eq:reciprocity-convention}
\inv_K\bigl(\kappa(a)\cup\lambda\bigr)=\lambda(a).
\end{equation}
Here $\inv_K:H^2(K,O(1))\to O$ is the local invariant. For a representation and its Tate dual we write $\pair{x}{y}=\inv_K(x\cup y)$, with evaluation on the coefficient modules understood. Local Tate duality, including compatibility with restriction and corestriction, is used in this convention; see \cite[Chapter I, Section 2]{Milne06}.

Let $\lambda_0\in H^1(K,O)$ be the character corresponding to $\ord_p$ in \eqref{eq:reciprocity-convention}. It is unramified. Let $\operatorname{Frob}_K$ denote arithmetic residue Frobenius, and put
\begin{equation}\label{eq:epsilon}
\epsilon=\lambda_0(\operatorname{Frob}_K)\in\{1,-1\}.
\end{equation}
Recording $\epsilon$ keeps the period construction independent of the convention choosing arithmetic or geometric Frobenius in the reciprocity map. All final arithmetic formulas use \eqref{eq:reciprocity-convention}.

\subsection{One-motives and their realizations}\label{subsec:motivedefinition}

A $1$-motive over $K$ is a two-term complex $M=[u:Y\to G]$, in degrees $-1$ and $0$. Here $Y$ is an \'etale lattice, namely a finite free abelian group with continuous finite Galois action, and $G$ is a semi-abelian variety, an extension of an abelian variety $A$ by a torus $T$. Write $s=\operatorname{rank}Y$, $a=\dim T$, and $g=\dim A$. This is Deligne's definition, recalled with Cartier duality in \cite[Section 1]{HS05}. The finite $p^m$-realization can be described by pairs $(y,P)$ with $[p^m]P=u(y)$, modulo pairs $(p^my',u(y'))$. It fits into
\begin{equation}\label{eq:finite-motive-realization}
0\longrightarrow G[p^m]\longrightarrow M[p^m]
\longrightarrow Y/p^mY\longrightarrow0.
\end{equation}
Set $T_p(M)=\varprojlim_m M[p^m]$ and $V_p(M)=T_p(M)\otimes\Qp$. For a semi-abelian variety $H$, we use $T_p(H)=\varprojlim_m H[p^m]$ and $V_p(H)=T_p(H)\otimes\Qp$. The Tate lattice of $M$ has rank $a+2g+s$. Its weight filtration is
\[
W_{-2}T_p(M)=T_p(T),\qquad W_{-1}T_p(M)=T_p(G),\qquad W_0T_p(M)=T_p(M).
\]
The graded pieces are $T_p(T)$, $T_p(A)$, and $Y\otimes\Zp$.

The Cartier dual is denoted $M^\vee$. Its lattice is $X^*(T)=\Hom(T,\mathbf G_m)$, its torus is $\Hom(Y,\mathbf G_m)$, and its abelian quotient is $A^\vee$. Its construction uses the Poincar\'e biextension, the universal line bundle on $A\times A^\vee$ with compatible group laws in both variables. The finite realizations have a perfect pairing into $\mu_{p^m}$, hence
\[
T_p(M^\vee)\simeq T_p(M)^\dagger,
\qquad T^\dagger:=\Hom_{\Zp}(T,\Zp)(1).
\]
We use $V^\dagger=V^\vee(1)$ for a rational representation. These realization pairings are distinguished from the $1$-motive hypercohomology pairings of \cite{HS05}.

\subsection{The toric representation}

For $Q$ in \eqref{eq:motive-intro}, write
\[
0\longrightarrow\Zp(1)^a\longrightarrow T_Q
\longrightarrow\Zp^s\longrightarrow0,
\qquad V_Q=T_Q\otimes\Qp.
\]
The covariant Tate module convention makes its extension matrix the Kummer matrix $(\kappa(q_{ij}))$. Write $\chi:\GK\to\Zp^\times$ for the cyclotomic character, defined by the action on $p$-power roots of unity, and set $\Gamma_K=\operatorname{Gal}(K(\mu_{p^\infty})/K)$. Choose vectors $v_i,w_j$ with
\begin{equation}\label{eq:galois-toric}
g v_i=\chi(g)v_i,\qquad
g w_j=w_j+\sum_i c_{ij}(g)v_i.
\end{equation}
The scalar Kummer cocycles satisfy
\[
c_{ij}(gh)=c_{ij}(g)+\chi(g)c_{ij}(h).
\]
For computations over $\Qp$, one may choose Kummer cocycles linearly on a basis of $W_{\Qp}$ and use these cocycles for all entries. Changing these choices gives isomorphic extensions.

We use Fontaine's rings $\Bcris^+\subset\Bst^+\subset\BdR^+$, their localizations by the cyclotomic period $t$, and the map $\theta:\BdR^+\to\Cp$. The rings $\Bcris$ and $\Bst$ carry Frobenius $\varphi$, and $\Bst$ carries monodromy $N$. One can take $t=\log[\varepsilon]$ for a compatible system $\varepsilon=(1,\zeta_p,\zeta_{p^2},\ldots)$ of roots of unity. Our convention is
\[
g t=\chi(g)t,\qquad\varphi(t)=pt,
\qquad N(\ell_p)=-1,\quad\varphi(\ell_p)=p\ell_p,
\quad\theta(\ell_p)=0.
\]
The element $\ell_p$ is the logarithmic period associated with compatible roots of $p$, with the chosen branch of logarithm. Standard Kummer periods $\ell_{ij}$ satisfy
\begin{equation}\label{eq:kummerperiod}
g\ell_{ij}=\ell_{ij}+c_{ij}(g)t,\quad
\varphi(\ell_{ij})=p\ell_{ij},\quad
N\ell_{ij}=-R_{ij},\quad\theta(\ell_{ij})=L_{ij}.
\end{equation}
For a unit the logarithmic period belongs to $\Bcris^+$. A general period is its valuation times $\ell_p$ plus a unit period. This is the crystalline logarithm construction of \cite[Section 6.1.2]{FO}, extended by the logarithmic generator of $\Bst$.

For $\bullet\in\{\mathrm{cris},\mathrm{st},\mathrm{dR}\}$, the period functor is
\[
D_\bullet(V)=(\mathbf B_\bullet\otimes_{\Qp}V)^{\GK}.
\]
The first two modules have coefficient field $K_0=K$ here, and the last has coefficient field $K$. The filtration is induced by the powers of $t$ in $\BdR$. Semistability means $\dim_K\Dst(V)=\dim_{\Qp}V$; crystallinity is the corresponding equality for $\Dcris$. A representation is de Rham if $\dim_K\DdR(V)=\dim_{\Qp}V$.

\begin{proposition}\label{prop:toricmodule}
Put $D=\Dst(V_Q)$. It has a $K$-basis $e_1,\ldots,e_a,f_1,\ldots,f_s$ with
\begin{align}\label{eq:phiN}
\varphi(e_i)&=p^{-1}e_i,&\varphi(f_j)&=f_j,&
Ne_i&=0,&Nf_j&=\sum_i R_{ij}e_i.
\end{align}
Frobenius is $\sigma$-semilinear. Under $D=K^a\oplus K^s$, the filtration is
\begin{equation}\label{eq:filtration}
\Fil^{-1}D=D,\qquad
\Fil^0D=\{(Ly,y):y\in K^s\},\qquad\Fil^1D=0.
\end{equation}
In particular $t_D=D/\Fil^0D$ is identified with $K^a$ by $(x,y)\mapsto x-Ly$.
\end{proposition}
\begin{proof}
In $\Bst\otimes V_Q$ set
\[
e_i=t^{-1}v_i,\qquad f_j=w_j-\sum_i\ell_{ij}e_i.
\]
Equations \eqref{eq:galois-toric} and \eqref{eq:kummerperiod} show that these vectors are Galois invariant. Their change-of-basis matrix is invertible over $\Bst$, so they give the full semistable comparison and prove semistability. Applying $\varphi$ and $N$ gives \eqref{eq:phiN}. Moreover
\[
f_j+\sum_iL_{ij}e_i
=w_j-\sum_i\frac{\ell_{ij}-L_{ij}}{t}v_i
\in\BdR^+\otimes V_Q.
\]
Conversely, an element $\sum_i x_i e_i+\sum_j y_jf_j$ has no negative $t$-coefficient precisely when $x=Ly$. This proves \eqref{eq:filtration} and the description of $t_D$.
\end{proof}

It follows that $\DdR(V_Q)=K\otimes_{K_0}\Dst(V_Q)$ with the filtration \eqref{eq:filtration}, and
\begin{equation}\label{eq:crystallinetoric}
\Dcris(V_Q)=\ker N=K^a\oplus(K\otimes_{\Qp}\ker R),\qquad
\dim_K\Dcris(V_Q)=a+s-\operatorname{rank}R.
\end{equation}
The displayed identity for $\DdR$ uses $K_0=K$ and the comparison just constructed.

In the dual bases giving positive coefficient evaluation, the Tate dual $V_Q^\dagger=\Hom(V_Q,\Qp)(1)$ is toric of the same type, with torus rank $s$, lattice rank $a$, and framed parameter matrix $Q^\dagger=(q_{ij}^{-1})^{\mathsf t}$. This is a framed model of the Cartier-dual realization; a geometric dual convention can use the opposite lattice frame. Thus
\begin{equation}\label{eq:dualmatrices}
R^\dagger=-R^{\mathsf t},\qquad L^\dagger=-L^{\mathsf t}.
\end{equation}
The period pairing takes values in $\Dst(\Qp(1))$. Its monodromy operators satisfy $b(Nx,y)+b(x,N^\dagger y)=0$.

\section{A normal form for semistable extensions}\label{sec:normalform}

\subsection{The extension complex}

Let $V$ be a semistable representation, put $D=\Dst(V)$, and write $t_D=D_K/\Fil^0D_K$, where $D_K=K\otimes_{K_0}D$ and $K_0$ is the maximal unramified subfield of $K$. The group $\Hst(K,V)$ consists of the classes in $H^1(K,V)$ whose middle representation in $0\to V\to E\to\Qp\to0$ is semistable. The classical reduced semistable complex is
\begin{equation}\label{eq:cst}
C_{\mathrm{st}}^\bullet(D):\quad
D\xrightarrow{d^0}D\oplus D\oplus t_D\xrightarrow{d^1}D,
\end{equation}
with
\[
d^0x=((1-\varphi)x,Nx,\bar x),\qquad
d^1(a,b,c)=Na-(1-p\varphi)b.
\]
The quotient $\bar x$ uses $D\to D_K\to t_D$. All spaces in the complex are $\Qp$-vector spaces. The relation $N\varphi=p\varphi N$ gives $d^1d^0=0$. Eliminating the filtration summand from the usual complex with $D\oplus\Fil^0D_K$ in degree zero gives \eqref{eq:cst}. Its $H^0$ is $H^0(K,V)$ and its $H^1$ is $\Hst(K,V)$. These are the classical Ext identifications of Bloch--Kato and Nekov\'a\v r, recalled in \cite[Theorem 1.1.4]{BLZ16}. In particular, this comparison concerns extensions with semistable middle representation; it need not exhaust $H^1(K,V)$.

For the admissibility input, a filtered $(\varphi,N)$-module means a finite-dimensional $K_0$-space with bijective semilinear Frobenius, nilpotent monodromy satisfying $N\varphi=p\varphi N$, and an exhaustive separated decreasing filtration on $D_K$. It is weakly admissible if $t_H(D)=t_N(D)$ and $t_H(D')\leq t_N(D')$ for every Frobenius- and monodromy-stable subspace with induced filtration. Here $t_H(D)=\sum_i i\dim_K\operatorname{gr}^iD_K$, and $t_N(D)=f^{-1}v_p(\det_{K_0}\varphi^f)$ for $f=[K_0:\Qp]$. The Colmez--Fontaine theorem identifies weakly admissible objects with semistable representations; we use the formulation in \cite[Theorem 1.1.2]{BLZ16}.

A lift $h\in\Dst(E)$ of $1$ and a filtered lift $h_0\in\Fil^0\DdR(E)$ give the cocycle
\begin{equation}\label{eq:extension-triple}
\bigl((1-\varphi)h,Nh,\overline{h-h_0}\bigr).
\end{equation}
Changing $h$ by $x\in D$ adds $d^0x$. We use this convention throughout, including the sign of the Hodge coordinate.

For a de Rham representation $V$, the Bloch--Kato local conditions are as follows:
\begin{equation}\label{eq:BK-localconditions}
\He(K,V)=\ker\bigl(H^1(K,V)\to H^1(K,\Bcris^{\varphi=1}\otimes V)\bigr),
\end{equation}
\[
\Hf(K,V)=\ker\bigl(H^1(K,V)\to H^1(K,\Bcris\otimes V)\bigr),
\]
\[
\Hg(K,V)=\ker\bigl(H^1(K,V)\to H^1(K,\BdR\otimes V)\bigr).
\]
The first group is the image of the Bloch--Kato exponential. We use their original definitions and fundamental exact sequences \cite[Section 3]{BK90}. For semistable $V$, the inclusions $\He\subset\Hf\subset\Hst$ follow either from these splitting conditions or from \eqref{eq:cst}. A Hodge-only cocycle $(0,0,z)$ represents an exponential class. Conversely, an invariant lift with coefficients in $\Bcris^{\varphi=1}$ has zero Frobenius and monodromy entries in \eqref{eq:extension-triple}, so every exponential class has such a representative.

\begin{lemma}\label{lem:frobenius}
For unramified $K/\Qp$, the maps $1-p^{-1}\sigma$ and $1-p\sigma$ on $K$ are invertible over $\Qp$. The sequence
\[
0\longrightarrow\Zp\longrightarrow\OK
\xrightarrow{1-\sigma}\OK\xrightarrow{\Tr_{K/\Qp}}\Zp
\longrightarrow0
\]
is exact. There is therefore $\vartheta\in\OK$ with $\Tr(\vartheta)=1$, even when $p\mid n$.
\end{lemma}
\begin{proof}
Since $\sigma^n=1$, geometric sums give inverses, with nonzero scalar factors $1-p^{-n}$ and $1-p^n$. The fixed ring of $\sigma$ on $\OK$ is $\Zp$. Reduction of trace modulo $p$ is the surjective finite-field trace, so integral trace is surjective. If $\Tr(a)=0$, solve $(1-\sigma)x=a$ modulo $p$ by the finite-field trace criterion. Subtract this solution and divide the error by $p$. The new error still has trace zero. Repetition and completeness give an integral solution. Finally, lift a residue element of nonzero trace and divide by its integral trace, a unit, to obtain $\vartheta$.
\end{proof}

\subsection{Removing a crystalline middle block}

The following linear-algebra calculation will be used twice: first for the toric representation, and then for a $1$-motive with an abelian quotient. Suppose that $K$ is unramified and that a filtered $(\varphi,N)$-module has a specified decomposition
\begin{equation}\label{eq:three-block-data}
D=D_T\oplus D_A\oplus D_Y,
\quad D_T=K^a,\quad D_Y=K^s,
\quad \varphi=(p^{-1}\sigma,\Phi_A,\sigma).
\end{equation}
Here $D_A$ is a finite-dimensional $K$-space and $\Phi_A$ is invertible and $\sigma$-semilinear. Assume
\begin{equation}\label{eq:middle-hypotheses}
N(x_T,x_A,y)=(Ry,0,0),\quad R\in M_{a,s}(\Qp),
\qquad 1-\Phi_A\ \hbox{and}\ 1-p\Phi_A\ \hbox{are bijective}.
\end{equation}
The last two assertions refer to $\Qp$-linear maps. Write $q:D\to t_D$ for the quotient map and put $\lambda=q|_{D_Y}:K^s\to t_D$.

\begin{lemma}[Three-block reduction]\label{lem:threeblock}
Assume \eqref{eq:three-block-data}--\eqref{eq:middle-hypotheses}. The complex \eqref{eq:cst} is quasi-isomorphic to
\begin{equation}\label{eq:mixed-reduced}
K^s\xrightarrow{\delta^0}K^s\oplus K^a\oplus t_D
\xrightarrow{\delta^1}K^a,
\qquad
\begin{aligned}
\delta^0 y&=((1-\sigma)y,Ry,\lambda y),\\
\delta^1(u,b,z)&=Ru-(1-\sigma)b.
\end{aligned}
\end{equation}
More precisely, the reduction removes acyclic direct summands after the Hodge-coordinate change
\begin{equation}\label{eq:hodge-elimination}
c\longmapsto c-q_T(1-p^{-1}\sigma)^{-1}a_T
                    -q_A(1-\Phi_A)^{-1}a_A.
\end{equation}
Here $q_T,q_A$ are the restrictions of $q$, and $a_T,a_A$ are the indicated coordinates of the first degree-one component.
\end{lemma}
\begin{proof}
Write the first two degree-one entries as $a=(a_T,a_A,u)$ and $b=(b_T,b_A,b_Y)$. Their differential is
\[
\bigl(Ru-(1-\sigma)b_T,-(1-p\Phi_A)b_A,-(1-p\sigma)b_Y\bigr).
\]
For a degree-zero vector $(x_T,x_A,y)$, the new Hodge coordinate in \eqref{eq:hodge-elimination} is $\lambda y$. Thus $(x_T,a_T)$ and $(x_A,a_A)$ form two-term acyclic summands, since $1-p^{-1}\sigma$ and $1-\Phi_A$ are bijective. The pairs formed by $b_A,b_Y$ and their degree-two targets are acyclic by the other two bijectivity assertions. The remaining coordinates and differentials are exactly \eqref{eq:mixed-reduced}. This proves the assertion as an explicit decomposition of complexes.
\end{proof}

\begin{theorem}[Cohomology of the reduced complex]\label{thm:threeblockcohomology}
Let $m=\dim_K t_D$, $d=\operatorname{rank}R$, and
\[
h=\dim_{\Qp}\{y\in\Qp^s:Ry=0,\ \lambda y=0\}.
\]
The complex \eqref{eq:mixed-reduced} has
\begin{equation}\label{eq:mixed-hzero-two}
H^0=\ker(R,\lambda)|_{\Qp^s},\qquad
H^2=\coker(R:\Qp^s\to\Qp^a).
\end{equation}
There is an exact sequence
\begin{equation}\label{eq:mixed-normal}
0\longrightarrow
\frac{\Qp^a\oplus t_D}{(R,\lambda)\Qp^s}
\longrightarrow H^1(C_{\mathrm{st}}^\bullet(D))
\xrightarrow{\Tr(u)}\ker R\longrightarrow0,
\end{equation}
split by $\alpha\mapsto[(\vartheta\alpha,0,0)]$ once $\Tr(\vartheta)=1$ is chosen. In particular,
\begin{equation}\label{eq:mixed-dim-general}
\dim_{\Qp}H^1=nm+a-d+h,\qquad \dim_{\Qp}H^2=a-d.
\end{equation}
\end{theorem}
\begin{proof}
Degree zero requires $y\in\Qp^s$ with $Ry=\lambda y=0$. In degree two, quotient first by $(1-\sigma)K^a$ and then use trace; the result is $\Qp^a/R\Qp^s$. A degree-one cocycle satisfies $Ru=(1-\sigma)b$, so $R\Tr(u)=0$. Set $\alpha=\Tr(u)$ and choose $y\in K^s$ with $(1-\sigma)y=u-\vartheta\alpha$, using \Cref{lem:frobenius}. Subtracting $\delta^0y$ leaves $(\vartheta\alpha,b',z')$ with $b'\in\Qp^a$. Further changes have $y\in\Qp^s$ and identify $(b',z')$ modulo $(R,\lambda)y$. This proves \eqref{eq:mixed-normal} and its splitting. The map $(R,\lambda)$ on $\Qp^s$ has rank $s-h$. The left term in \eqref{eq:mixed-normal} has dimension $a+nm-s+h$, while $\dim\ker R=s-d$. Adding gives \eqref{eq:mixed-dim-general}.
\end{proof}

\begin{corollary}[Toric normal form]\label{thm:normalform}
For the toric module of \Cref{prop:toricmodule}, the reduced complex is
\begin{equation}\label{eq:reducedcomplex}
K^s\xrightarrow{\delta^0}K^s\oplus K^a\oplus K^a
\xrightarrow{\delta^1}K^a,
\quad
\delta^0y=((1-\sigma)y,Ry,-Ly),\quad
\delta^1(u,b,z)=Ru-(1-\sigma)b.
\end{equation}
Its cohomology is
\begin{align}\label{eq:hzero-two}
H^0&=\ker((R,L):\Qp^s\to\Qp^a\oplus K^a),&
H^2&=\coker R,
\end{align}
and
\begin{equation}\label{eq:h1-exact}
0\longrightarrow\frac{\Qp^a\oplus K^a}{(R,-L)\Qp^s}
\longrightarrow\Hst(K,V_Q)\xrightarrow{\Tr(u)}\ker R\longrightarrow0.
\end{equation}
Writing $d=\operatorname{rank}R$ and $h=\dim\ker(R,L)$ gives
\begin{equation}\label{eq:dimension-general}
\dim\Hst(K,V_Q)=an+a-d+h,\qquad\dim H^2=a-d.
\end{equation}
\end{corollary}
\begin{proof}
Apply \Cref{lem:threeblock,thm:threeblockcohomology} with $D_A=0$, $t_D=K^a$, and $\lambda=-L$, as in \Cref{prop:toricmodule}. The Ext comparison following \eqref{eq:cst} identifies degree one with $\Hst$.
\end{proof}

The cokernel of the Kummer map
\begin{equation}\label{eq:kummer-quotient}
Q:\Qp^s\longrightarrow W_{\Qp}^a
\end{equation}
is the left term of \eqref{eq:h1-exact}: use $(\ord_p,-\log_p)$ on each component. The long exact Galois cohomology sequence injects this cokernel into $H^1(K,V_Q)$. This determines the arithmetic meaning and sign of the two normal-form coordinates. For a general $1$-motive the corresponding Hodge space will be $\operatorname{Lie}(G)$ rather than $K^a$; \Cref{sec:mixedmotives} makes that replacement geometrically.

\section{Actual toric period representatives}\label{sec:actualperiods}

\subsection{Two elementary period constructions}

We use the following standard exact sequence of crystalline logarithms:
\begin{equation}\label{eq:fundamental-weightone}
0\longrightarrow\Qp t\longrightarrow
(\Bcris^+)^{\varphi=p}\xrightarrow{\theta}\Cp
\longrightarrow0.
\end{equation}
For a unit $u\in1+p\mathcal O_{\Cp}$, choose compatible $p$-power roots, let $\widetilde u$ be the resulting tilt element, and form $\log[\widetilde u]$. It belongs to $(\Bcris^+)^{\varphi=p}$ and has $\theta$-value $\log_p u$. The construction and exact sequence are given in \cite[Section 6.1.2, especially Remark 6.9]{FO}. We need the exact kernel in \eqref{eq:fundamental-weightone}, not just existence of a de Rham logarithm.

\begin{lemma}[An unramified period]\label{lem:tau}
Fix $\vartheta\in\OK$ with trace one. There exists
$\tau\in W(\overline{\Fp})\subset\Bcris^+$ such that
\begin{equation}\label{eq:tau}
\varphi\tau-\tau=\epsilon\vartheta,\qquad
g\tau-\tau=\lambda_0(g)\quad(g\in\GK).
\end{equation}
Moreover $\tau\in\Rplus$ and $\bar\tau=\theta(\tau)$ is integral in $\Cp$.
\end{lemma}
\begin{proof}
The map $\varphi-1$ on $W(\overline{\Fp})$ is surjective. To see this, solve $x^p-x=a$ in the algebraically closed residue field, lift a solution to Witt vectors, subtract its error, and repeat modulo successive powers of $p$. Completeness gives a solution of $(\varphi-1)\tau=\epsilon\vartheta$. Iteration yields $\varphi^n\tau-\tau=\epsilon\Tr(\vartheta)=\epsilon$. On the constant Witt-vector ring, arithmetic residue Frobenius of $K$ acts as $\varphi^n$. Hence its difference on $\tau$ is $\epsilon$. The asserted Galois identity follows first for integral powers of Frobenius and then for their profinite closure by continuity. Integrality follows from the construction. Finally
\[
\varphi^{-m}\tau=\tau-\epsilon\sum_{j=1}^m\sigma^{-j}(\vartheta)
\in\Bcris^+
\]
for every $m\geq0$. The characterization of $\Rplus$ in \eqref{eq:positive-rings} below gives the claim.
\end{proof}

Fix, linearly on $W_{\Qp}$, a crystalline or semistable Kummer logarithm $\ell(a)$ for each $a\in W_{\Qp}$. Thus
\[
\theta\ell(a)=\log_p(a),\quad N\ell(a)=-\ord_p(a),\quad
\varphi\ell(a)=p\ell(a),\quad
g\ell(a)-\ell(a)=c_a(g)t.
\]
Linearity can be obtained by fixing these logarithms on a basis of $W$.

\begin{theorem}[Construction of all semistable extension classes]\label{thm:actualperiods}
Let $\alpha\in\ker(R:\Qp^s\to\Qp^a)$ and $a'=(a'_i)\in W_{\Qp}^a$. Put
\[
\ell_i(\alpha)=\sum_j\alpha_j\ell_{ij}\in(\Bcris^+)^{\varphi=p},
\qquad \omega_i=\ell(a'_i),\quad m_i=\ord_p(a'_i),\quad z_i=\log_p(a'_i).
\]
Choose $\rho_i\in(\Bcris^+)^{\varphi=p}$ with
\begin{equation}\label{eq:rho}
\theta(\rho_i)=\bar\tau(L\alpha)_i.
\end{equation}
Then
\begin{equation}\label{eq:di}
d_i(g)=\frac{g\rho_i-\rho_i-\lambda_0(g)\ell_i(\alpha)}{t}
\end{equation}
belongs to $\Qp$. The following formulas define a continuous Galois representation $E_{a',\alpha}$ extending $\Qp$ by $V_Q$:
\begin{equation}\label{eq:actualextension}
g z=z+\lambda_0(g)\sum_j\alpha_jw_j+
\sum_i\bigl(c_{a'_i}(g)+d_i(g)\bigr)v_i.
\end{equation}
Here the letter $z$ without an index denotes the new representation basis vector, while $(z_i)$ is the logarithmic vector above. An invariant period lift of $1$ is
\begin{equation}\label{eq:actualh}
h=z-\tau\sum_j\alpha_j f_j-\sum_i(\omega_i+\rho_i)e_i.
\end{equation}
In $\Dst(E_{a',\alpha})$, this lift satisfies
\begin{equation}\label{eq:actualhdata}
(1-\varphi)h=\epsilon\vartheta\sum_j\alpha_jf_j,
\qquad Nh=\sum_i m_i e_i,
\qquad h+\sum_i z_i e_i\in\Fil^0\DdR(E_{a',\alpha}).
\end{equation}
Every element of $\Hst(K,V_Q)$ has such a representative.
\end{theorem}
\begin{proof}
Since $R\alpha=0$, the logarithmic $\ell_p$ coefficients in $\ell_i(\alpha)$ cancel. Thus $\ell_i(\alpha)$ is crystalline and has Frobenius eigenvalue $p$. The existence of $\rho_i$ follows from \eqref{eq:fundamental-weightone}. It has a concrete logarithmic construction: for integral $\alpha$, take compatible roots of
\[
u_i=\exp\bigl(\bar\tau(L\alpha)_i\bigr)\in1+p\mathcal O_{\Cp}
\]
and put $\rho_i=\log[\widetilde u_i]$. For arbitrary $\alpha$, choose these lifts on a $\Zp$-basis of $\ker R\cap\Zp^s$ and extend $\Qp$-linearly. This choice makes the construction linear in $\alpha$.

The numerator of \eqref{eq:di} lies in $(\Bcris^+)^{\varphi=p}$. Its $\theta$-value is
\[
g\bigl(\bar\tau(L\alpha)_i\bigr)-\bar\tau(L\alpha)_i
-\lambda_0(g)(L\alpha)_i=0.
\]
Exactness of \eqref{eq:fundamental-weightone} therefore proves $d_i(g)\in\Qp$. Write $c_i(\alpha,g)=\sum_j\alpha_jc_{ij}(g)$. Directly from \eqref{eq:di},
\begin{equation}\label{eq:d-cocycle}
d_i(gh)=d_i(g)+\chi(g)d_i(h)+\lambda_0(h)c_i(\alpha,g).
\end{equation}
The extra last term is exactly the term obtained by applying $g$ to $\sum_j\alpha_jw_j$ in \eqref{eq:actualextension}. This proves the representation identity. The maps are continuous because the period actions, the Kummer cocycles, and the kernel identification in \eqref{eq:fundamental-weightone} are continuous.

Set $f_\alpha=\sum_j\alpha_jf_j$. Since $f_\alpha$ and each $e_i$ are Galois invariant, the Galois difference of \eqref{eq:actualh} is
\[
\lambda_0(g)\sum_j\alpha_j(w_j-f_j)
+\sum_i\bigl(t d_i(g)-(g\rho_i-\rho_i)\bigr)e_i=0.
\]
Here $w_j-f_j=\sum_i\ell_{ij}e_i$ and \eqref{eq:di} cancel the two sums. The vectors $e_i,f_j,h$ have invertible comparison matrix over $\Bst$, so they prove semistability and form a basis of $\Dst(E_{a',\alpha})$. Frobenius gives the first equality in \eqref{eq:actualhdata}. Monodromy gives the second because $Nf_\alpha=0$, $N\rho_i=N\tau=0$, and $N\omega_i=-m_i$.

To verify the filtration, expand $h+\sum_i z_i e_i$ in the original representation basis. The coefficient numerator multiplying $t^{-1}v_i$ is
\[
\tau\ell_i(\alpha)-\rho_i-\omega_i+z_i.
\]
It has $\theta$-value zero by \eqref{eq:rho}, so its quotient by $t$ belongs to $\BdR^+$. This proves the last assertion of \eqref{eq:actualhdata}. The class of the extension in the reduced complex is consequently
\begin{equation}\label{eq:period-normalform}
\bigl(\epsilon\vartheta\alpha,m,-z\bigr).
\end{equation}
Every normal form of \Cref{thm:normalform} has this shape: choose $\alpha$ from its trace coordinate, with the factor $\epsilon^{-1}$, and use \eqref{eq:val-log} to choose $a'_i$ with any specified valuation and logarithm. Hence the construction exhausts $\Hst(K,V_Q)$.
\end{proof}

\begin{remark}\label{rem:choice-rho}
Two lifts in \eqref{eq:rho} differ by $t\Qp$. The corresponding changes in $d_i$ are scalar coboundaries. Their extension classes agree. The splitting by $\alpha$ depends on the trace-one element $\vartheta$ and the unramified period, just as the splitting in \Cref{thm:normalform} does. The boundary and quotient maps used in the duality theorem below are independent of that splitting.
\end{remark}

\subsection{Positive periods and the one-pole bound}

We use the rings in Colmez's diagram \cite[Section 8.3.4]{Col04}:
\begin{equation}\label{eq:positive-rings}
\Rplus=\bigcap_{m\geq0}\varphi^m(\Bcris^+),\qquad
\Rlog=\bigcap_{m\geq0}\varphi^m(\Bst^+)=\Rplus[\ell_p].
\end{equation}
These intersections use the standard common period overring and the chosen logarithmic coordinate. Frobenius is an automorphism of $\Rplus$. Thus its only left inverse on this perfect ambient ring is $\varphi^{-1}$. The trace operator customarily called $\psi$ on a decompleted Robba ring belongs to a different coefficient construction.

\begin{proposition}\label{prop:positivebound}
For $V_Q$ and every extension constructed in \Cref{thm:actualperiods}, the displayed period bases satisfy
\[
t e_i,\quad t f_j,\quad t h\ \in\ \Rlog\otimes_{\Qp}E_{a',\alpha}.
\]
Consequently $t\Dst(E_{a',\alpha})\subset\Rlog\otimes E_{a',\alpha}$ under the semistable comparison. When $a>0$, exponent one is the least nonnegative uniform exponent for the full period space.
\end{proposition}
\begin{proof}
A crystalline logarithm $u$ with $\varphi u=pu$ lies in every $\varphi^m(\Bcris^+)$ because $\varphi^{-m}u=p^{-m}u\in\Bcris^+$. The same argument applies to $t$. General Kummer logarithms lie in $\Rlog$, and \Cref{lem:tau} gives $\tau\in\Rplus$. Multiplying the definitions of $e_i,f_j,h$ by $t$ now gives the claim term by term. In particular no inverse Frobenius iteration introduces an additional power of $t^{-1}$.

The chosen embedding sends $\Rlog$ into $\BdR^+$. Since $e_i=t^{-1}v_i$ and $t^{-1}\notin\BdR^+$, the full period space is not contained in $\Rlog\otimes E_{a',\alpha}$ when a toric vector is present. This proves minimality.
\end{proof}

This bound is an ambient rational period bound. It does not identify an integral overconvergent lattice or the denominators of a residue regulator.

\begin{proposition}[The positive bound at Hodge height one]\label{prop:generalpositive}
Let $W$ be a semistable representation of a finite extension of $\Qp$ such that $\Fil^{-1}\DdR(W)=\DdR(W)$ and $\Fil^1\DdR(W)=0$. In the positive rings \eqref{eq:positive-rings},
\begin{equation}\label{eq:general-positive-bound}
t\Dst(W)\subset\Rlog\otimes_{\Qp}W,
\qquad
t\Dcris(W)\subset\Rplus\otimes_{\Qp}W.
\end{equation}
If $\DdR(W)/\Fil^0\DdR(W)\ne0$, exponent one is the least nonnegative exponent for the full semistable period space. In particular, this bound applies to $V_M$ under \Cref{thm:mixedrealization} and to every semistable extension of $\Qp$ by $V_M$; it is optimal when $a+g>0$.
\end{proposition}
\begin{proof}
We use the standard positive-period comparison: if $\Fil^0\DdR(W')=\DdR(W')$, its semistable periods lie in $\Rlog\otimes W'$, and its crystalline periods lie in $\Rplus\otimes W'$. This is the positive case and Frobenius regularization recalled in \cite[Section 3.2, preceding Proposition 3.4]{Ber02}; see \cite[Sections IV.3.1--IV.3.3]{BerIntro02}. Apply it to $W'=W(-1)$, whose zeroth filtration is the whole module. The period of $\Qp(-1)$ is $t$, so untwisting gives \eqref{eq:general-positive-bound}. If the full space $\Dst(W)$ lay in $\Rlog\otimes W$, its image under the chosen embedding into $\BdR$ would lie in $\BdR^+\otimes W$. Semistable comparison would then give $\Fil^0\DdR(W)=\DdR(W)$, contrary to the hypothesis. Finally, the filtrations for $V_M$ and its extensions have only the two indicated jumps, by the $p$-divisible comparison and strictness. Their de Rham quotients have dimension $a+g$.
\end{proof}

The toric construction in \Cref{prop:positivebound} verifies this bound directly on each displayed coefficient. \Cref{prop:generalpositive} extends the bound using the established positive comparison. It supplies a uniform ambient bound in $\Rplus$ and its logarithmic extension; a formula for a decompleted basis requires further data, specified in \Cref{subsec:relativeHerr,subsec:mixedscope}.

\section{One-motives with good abelian quotient}\label{sec:mixedmotives}

\subsection{The integral model and its point coordinates}

Continue to assume that $K/\Qp$ is unramified and $p>2$. Let
\begin{equation}\label{eq:mixed-motive}
M=[u:Y\to G],\qquad Y=\mathbf Z^s,\qquad
0\longrightarrow T\xrightarrow{i}G\longrightarrow A\longrightarrow0,
\quad T=\mathbf G_m^a.
\end{equation}
Assume that this extension has a semi-abelian model
\begin{equation}\label{eq:good-model}
0\longrightarrow\mathcal T\longrightarrow\mathcal G
\longrightarrow\mathcal A\longrightarrow0
\end{equation}
over $\OK$, with $\mathcal T=\mathbf G_{m,\OK}^a$ and $\mathcal A$ an abelian scheme of relative dimension $g$. The point $u(y)$ need not belong to $\mathcal G(\OK)$. Thus the hypothesis on $G$ allows nonzero monodromy for $M$. The semi-abelian extension itself may be nonsplit.

There is a canonical exact sequence
\begin{equation}\label{eq:semiabelian-valuation}
0\longrightarrow\mathcal G(\OK)\longrightarrow G(K)
\xrightarrow{\nu_G}\mathbf Z^a\longrightarrow0,
\end{equation}
whose restriction to $T(K)$ is componentwise valuation. To prove this, properness gives $A(K)=\mathcal A(\OK)$, and $H^1(\OK,\mathcal T)=\operatorname{Pic}(\OK)^a=0$ makes $\mathcal G(\OK)\to A(K)$ surjective. Subtract an integral lift of the image of a given point of $G(K)$. Its remaining toric point determines a valuation vector, independent of the lift. The kernel and surjectivity follow from $K^\times/\OK^\times=\mathbf Z$ in each toric coordinate.

For the fixed basis of $Y$, define
\begin{equation}\label{eq:geometric-R-u0}
R=\nu_G\circ u:\mathbf Z^s\longrightarrow\mathbf Z^a,
\qquad
u_0(y)=u(y)-i(p^{Ry})\in\mathcal G(\OK).
\end{equation}
Here $p^{Ry}$ means the tuple with entries $p^{(Ry)_i}$, and addition denotes the commutative group law on $G$. This is the good-point and uniformizer decomposition of Raynaud \cite[Proposition 4.5.1 and Remark 4.5.2(i)]{Ray94} in the present split integral setting. The preceding argument proves it directly in that setting. The logarithmic refinement is treated in \cite[Propositions 4.2--4.7]{WZ24}.

Let $\mathfrak g=\operatorname{Lie}(G)$, the tangent space at the identity, of dimension $a+g$ over $K$. The formal logarithm on a sufficiently small open subgroup of $\mathcal G(\OK)$ extends uniquely to a homomorphism
\[
\log_{\mathcal G}:\mathcal G(\OK)\longrightarrow\mathfrak g.
\]
Indeed, this group is compact and a sufficiently small formal subgroup has finite index. For $P\in\mathcal G(\OK)$ choose a positive integer $b$ with $[b]P$ in that subgroup and put $\log_{\mathcal G}(P)=b^{-1}\log([b]P)$. The homomorphism property of the formal logarithm proves independence of $b$. On a sufficiently small neighborhood its inverse is the formal exponential. Its restriction to $\mathcal T(\OK)$ is the usual logarithm on units. Define the $K$-linear map
\begin{equation}\label{eq:geometric-lambda}
\lambda:K^s\longrightarrow\mathfrak g,\qquad
\lambda(y)=-\log_{\mathcal G}(u_0(y))\quad(y\in\mathbf Z^s),
\end{equation}
by linear extension. No logarithm of a general nonintegral point is needed in this definition.

\begin{theorem}[Geometric coordinates of the semistable realization]\label{thm:mixedrealization}
For $M$ satisfying \eqref{eq:mixed-motive}--\eqref{eq:good-model}, the representation $V_M=V_p(M)$ is semistable. Put $D_G=\Dcris(V_p(G))$ and $D_A=\Dcris(V_p(A))$. There is a unique Frobenius-compatible splitting of its weight filtration,
\begin{equation}\label{eq:mixed-D}
\Dst(V_M)=K^a\oplus D_A\oplus K^s,
\qquad \varphi=(p^{-1}\sigma,\Phi_A,\sigma).
\end{equation}
In this decomposition, $D_G=K^a\oplus D_A$ and
\begin{equation}\label{eq:mixed-N-Hodge}
N(x_T,x_A,y)=(Ry,0,0),\qquad
t_M:=\DdR(V_M)/\Fil^0\DdR(V_M)\simeq\mathfrak g,
\qquad q|_{K^s}=\lambda.
\end{equation}
The identification of $t_M$ with $\mathfrak g$ is induced by the inclusion of $V_p(G)$ and the de Rham comparison for its $p$-divisible group. Both $1-\Phi_A$ and $1-p\Phi_A$ are bijective as $\Qp$-linear maps. Consequently,
\begin{equation}\label{eq:mixed-cris}
\Dcris(V_M)=D_G\oplus(K\otimes_{\Qp}\ker R),\qquad
\operatorname{rank}_K N=\operatorname{rank}_{\Qp}R.
\end{equation}
\end{theorem}
\begin{proof}
The groups $\mathcal G[p^m]$ form a $p$-divisible group of height $a+2g$ and dimension $a+g$. For the integral $1$-motive $M_0=[Y\xrightarrow{u_0}\mathcal G]$, its finite realizations fit into
\[
0\longrightarrow\mathcal G[p^m]\longrightarrow M_0[p^m]
\longrightarrow Y/p^mY\longrightarrow0.
\]
They form a $p$-divisible group over $\OK$ of the same dimension and height $a+2g+s$. The construction is classical; see \cite[Section 2]{Ray94} and \cite[Proposition 4.5 and Definition 4.6]{WZ24}, whose trivial-log case is the construction just used. The crystalline comparison for $p$-divisible groups makes $V_p(G)$ and $V_p(M_0)$ crystalline, with their covariant de Rham quotients identified with $\mathfrak g$; see \cite[Sections 7.2 and 12.3]{BC09}. These comparison statements are the geometric input to this proof.

The remaining homomorphism $y\mapsto i(p^{Ry})$ is toric. Its Kummer extension is the pushout along $V_p(T)\to V_p(G)$ of the extension attached to $[Y\to T]$ with parameters $p^{R_{ij}}$. Its semistable realization is calculated in \Cref{prop:toricmodule}, with monodromy matrix $R$ and logarithmic matrix zero. Kummer classes are additive in their defining points. Hence the extension of $Y\otimes\Qp$ by $V_p(G)$ attached to $M$ is the Baer sum of this extension and the crystalline extension attached to $M_0$. The semistable Ext comparison following \eqref{eq:cst} is compatible with pushout and addition. It therefore proves semistability of $V_M$.

We next justify the Frobenius splitting, including the spectral assertion needed below. For good-reduction $A$, the eigenvalues of $\Phi_A^n$ in the covariant Tate-module convention are the inverses of the Frobenius eigenvalues of its special fiber. Weil's theorem gives complex absolute value $p^{-n/2}$ for every eigenvalue. Thus the characteristic polynomials of $\Phi_A^n$, $p^{-n}I$, and $I$ are pairwise coprime. This is the usual Weil property under the crystalline--Dieudonn\'e comparison for an abelian scheme; see \cite[Sections 7.2--7.3]{BC09} and \cite[Chapter II, Sections 1--2]{MilneAV}. Polynomial projectors for $\varphi^n$ split the weight filtration uniquely. Their coefficients are $\sigma$-fixed, so the projectors commute with $\varphi$ and the splitting is Frobenius compatible. The same spectral calculation shows that $1-\Phi_A$ and $1-p\Phi_A$ have zero kernel and hence are bijective over $\Qp$.

Monodromy vanishes on $D_G$, since $V_p(G)$ is crystalline. In the toric summand constructed above it sends the lattice basis to $R$ in $K^a$. The crystalline summand contributes zero. Additivity of extension data in \eqref{eq:extension-triple} therefore gives the first formula of \eqref{eq:mixed-N-Hodge}.

Strictness of the de Rham comparison identifies $t_M$ with $t_{D_G}$: the quotient representation $Y\otimes\Qp$ has zero de Rham quotient. The latter space is $\mathfrak g$ by the $p$-divisible comparison. Since $1-\varphi$ on $D_G$ is bijective, the Bloch--Kato exponential identifies $\mathfrak g$ with $\Hf(K,V_p(G))=\He(K,V_p(G))$. On the small formal subgroup, Bloch--Kato's exponential is the formal exponential followed by the Kummer map \cite[Example 3.10.1]{BK90}; see \cite[Section V.2]{BerIntro02} for this identification. Thus the Hodge entry of its Kummer extension, with the convention \eqref{eq:extension-triple}, is minus its formal logarithm. Multiplying any integral point into that subgroup and dividing the resulting rational identity extends the calculation to all integral points. It follows that the $\varphi$-fixed lattice lift for $M_0$ has quotient $-\log_{\mathcal G}(u_0(y))$. The toric $p^{Ry}$ extension has zero Hodge entry by \Cref{prop:toricmodule}. Addition gives $q|_{K^s}=\lambda$. Taking the kernel of the displayed monodromy proves \eqref{eq:mixed-cris}.
\end{proof}

The splitting in \eqref{eq:mixed-D} respects Frobenius, but it need not respect the Hodge filtration. In particular, the filtered module $D_G$ still records the extension of $A$ by $T$. The pair $(R,\lambda)$ determines the reduced extension complex once the geometric identifications above are made; it does not classify all mixed $1$-motives or reconstruct their abelian period matrices.

\begin{theorem}[Semistable extensions of a mixed $1$-motive]\label{thm:mixednormalform}
Under the hypotheses of \Cref{thm:mixedrealization}, the complex computing $\Hst(K,V_M)$ is
\begin{equation}\label{eq:mixed-motivic-complex}
K^s\longrightarrow K^s\oplus K^a\oplus\mathfrak g
\longrightarrow K^a,
\quad
y\longmapsto((1-\sigma)y,Ry,\lambda y),
\quad (u,b,z)\longmapsto Ru-(1-\sigma)b.
\end{equation}
There is an exact sequence
\begin{equation}\label{eq:mixed-motivic-normal}
0\longrightarrow\frac{\Qp^a\oplus\mathfrak g}{(R,\lambda)\Qp^s}
\longrightarrow\Hst(K,V_M)\xrightarrow{\Tr(u)}\ker R\longrightarrow0.
\end{equation}
With $d=\operatorname{rank}R$ and $h=\dim\ker(R,\lambda)|_{\Qp^s}=\dim H^0(K,V_M)$,
\begin{equation}\label{eq:mixed-motivic-dim}
\dim\Hst(K,V_M)=n(a+g)+a-d+h,\qquad
H^2(C_{\mathrm{st}}^\bullet(\Dst(V_M)))=\coker R.
\end{equation}
The trace-one normal forms and the splitting of \Cref{thm:threeblockcohomology} apply without change.
\end{theorem}
\begin{proof}
Equations \eqref{eq:mixed-D}--\eqref{eq:mixed-N-Hodge} verify every hypothesis of \Cref{lem:threeblock}. Now apply \Cref{thm:threeblockcohomology} with $t_D=\mathfrak g$ and $m=a+g$. The Ext comparison gives the asserted interpretation of degree one.
\end{proof}

\begin{corollary}[Changing the logarithm branch]\label{cor:mixedbranch}
Replace $p$ by the uniformizer $p'=p\eta$, where $\eta\in\OK^\times$, and normalize the new toric logarithm by $\log_{p'}(p')=0$. Let $j:\operatorname{Lie}(T)=K^a\to\mathfrak g$ be the tangent map of $i$. The point coordinates obtained from this normalization satisfy
\begin{equation}\label{eq:mixed-branch-change}
R'=R,\qquad \lambda'=\lambda+j\bigl(\log_p(\eta)R\bigr).
\end{equation}
In particular $\lambda|_{\ker R}$ and the quotient $\mathfrak g/\lambda(\ker R)$ are unchanged.
\end{corollary}
\begin{proof}
The new integral point is $u'_0(y)=u_0(y)-i(\eta^{Ry})$. Apply $-\log_{\mathcal G}$ and use its compatibility with the torus logarithm. Formula \eqref{eq:mixed-branch-change} follows. The restriction of the additional term to $\ker R$ vanishes. This compares point coordinates with their corresponding normalized semistable embeddings; if the original branch is retained, the logarithm of the new uniformizer instead appears in the toric summand and cancels this change.
\end{proof}

\begin{corollary}[The logarithmic monodromy interpretation]\label{cor:logarithmicmonodromy}
The matrix $R$ in \eqref{eq:geometric-R-u0} represents Raynaud's integral pairing
$Y\times X^*(T)\to\mathbf Z$. Its reduction modulo $p^m$ is the geometric pairing in the logarithmic finite realization $M[p^m]$. Under the logarithmic $p$-divisible comparison, its rational realization is the operator in \eqref{eq:mixed-N-Hodge}.
\end{corollary}
\begin{proof}
Evaluation of a toric character on \eqref{eq:geometric-R-u0} gives $\nu_G(u(y))$ in that character direction, exactly Raynaud's definition in \cite[Section 4.3]{Ray94}. The finite-level logarithmic pairing is the reduction of this chart pairing, by \cite[Proposition 4.7]{WZ24}. Compatibility of Kato and Fontaine monodromy is \cite[Theorem F and Section 6.5]{BWZ26}; in our frames the explicit computation is already \Cref{thm:mixedrealization}. These are compatible descriptions of the same geometric monodromy.
\end{proof}

\subsection{Raynaud uniformization and unramified descent}

\begin{proposition}[Transfer through uniformization]\label{prop:uniformizationtransfer}
Let $K'/K$ be a finite unramified Galois extension, with group $\Omega$, and let $M$ be a $1$-motive over $K$. Suppose that a Raynaud replacement $\widetilde M=[\widetilde Y\to\widetilde G]$ over $K$ has the following properties: its Tate realization is identified with that of $M$; over $K'$ its lattice and torus are split; and its semi-abelian part has a model of the form \eqref{eq:good-model}. Then $V_M$ is semistable over $K$. Restriction induces isomorphisms
\begin{equation}\label{eq:unramified-descent}
H^1_\star(K,V_M)\xrightarrow{\sim}
H^1_\star(K',V_{\widetilde M})^\Omega,
\qquad \star\in\{e,f,\mathrm{st}\}.
\end{equation}
The normal forms, boundary, exponential quotient, and reciprocity of \Cref{thm:mixednormalform,thm:boundary,thm:semistablereciprocity} therefore give the corresponding rational statements over $K$ after taking $\Omega$-invariants. Arithmetic pairings descend by
\begin{equation}\label{eq:unramified-pairing-factor}
\langle x,y\rangle_K=[K':K]^{-1}
\langle\operatorname{res}x,\operatorname{res}y\rangle_{K'}.
\end{equation}
\end{proposition}
\begin{proof}
Raynaud's replacement preserves every Tate realization, including the $p$-adic one, by \cite[Theorem 4.2.2]{Ray94}. This is a statement about finite realizations and rigid uniformization; it does not use the prime-to-$p$ inertia formula of that paper's Section 4.6. Over $K'$, \Cref{thm:mixedrealization} proves semistability. Semistability and crystallinity descend through a finite unramified extension; this is the unramified base-change property of the period functors \cite[Proposition 9.3.1]{BC09}. It applies to the middle representation of an extension as well.

The rational Hochschild--Serre sequence gives $H^1(K,V_M)=H^1(K',V_M)^\Omega$, since averaging over the finite group is exact over $\Qp$. The period splitting conditions and the unramified descent just cited show that restriction identifies the finite and semistable subgroups. For the exponential subgroup one can use its defining $\Bcris^{\varphi=1}$ splitting condition, or the fundamental exact sequence and the same descent. All constructions in the preceding theorems are functorial, so their exact sequences and quotient pairings descend. Finally, restriction multiplies the local invariant in degree two by $[K':K]$, proving \eqref{eq:unramified-pairing-factor}.
\end{proof}

The hypothesis occurs for $1$-motives admitting the stated semistable Raynaud model: replacing the abelian quotient by its uniformization introduces additional toric and lattice directions and a good abelian quotient. Those new ranks must be used in the formulas. For a purely abelian variety $B$ with a split Raynaud presentation
\[
B^{\mathrm{an}}=\widetilde G^{\mathrm{an}}/\widetilde Y,
\qquad \widetilde T=\mathbf G_m^b,\quad \widetilde Y=\mathbf Z^b,
\]
suppose its valuation pairing $\widetilde R$ has full rank, as in the nondegenerate monodromy pairing of proper semistable uniformization. If the good abelian quotient has dimension $g_0$, then \Cref{thm:mixednormalform,thm:boundary} immediately give
\begin{equation}\label{eq:uniformized-abelian-local}
\He(K,V_p(B))=\Hf(K,V_p(B))=\Hst(K,V_p(B)),
\qquad \dim_{\Qp}=n(b+g_0)=n\dim B.
\end{equation}
Indeed $a=s=d=b$ and $h=0$. Thus nonzero coefficient monodromy need not produce a nonzero new-monodromy obstruction: here $\coker\widetilde R=0$. A Tate elliptic curve, or its product with a good abelian variety, gives an immediate example. The equality of the exponential, finite, and geometric subgroups for an abelian Tate module is classical \cite[Example 3.11]{BK90}; \Cref{cor:universalquotientduality} identifies the geometric subgroup with the semistable subgroup here. The comparison of Raynaud and Fontaine monodromy for abelian varieties is classical as well \cite{CI99}. Formula \eqref{eq:uniformized-abelian-local} recovers these facts through the specified uniformizing $1$-motive.

The invariants in \eqref{eq:unramified-descent} are rational. When $p\nmid[K':K]$, averaging works over $\Zp$ as well. If $p\mid[K':K]$, integral restriction and descent require their Hochschild--Serre terms; rational averaging does not provide an integral comparison. A ramified extension with potentially semistable descent data requires the corresponding filtered $(\varphi,N)$-module with its descent action, rather than merely taking invariant vectors in the split normal form.

\section{Monodromy, local conditions, and finite-level reciprocity}\label{sec:localconditions}

\subsection{A boundary for every semistable representation}

Recall the period splitting definitions in \eqref{eq:BK-localconditions} and the following paragraph. We first extract the monodromy obstruction without a motivic hypothesis.

For any semistable $V$, write $D=\Dst(V)$ and let $\bar\varphi$ be the induced Frobenius on $D/ND$. This is defined because $ND$ is Frobenius stable. Define the monodromy obstruction space and the crystalline Frobenius quotient by
\begin{equation}\label{eq:general-obstruction-spaces}
\mathcal M_N(V)=\ker(1-p\bar\varphi:D/ND\to D/ND),
\qquad
\mathcal U(V)=\frac{D^{N=0}}{(1-\varphi)D^{N=0}}.
\end{equation}
These are $\Qp$-vector spaces. The following statement allows any finite extension $K/\Qp$, with coefficient field $K_0$ as in \eqref{eq:cst}.

\begin{theorem}[Universal monodromy boundary]\label{thm:universalboundary}
There are canonical exact sequences
\begin{equation}\label{eq:universal-boundary}
0\longrightarrow\Hf(K,V)\longrightarrow\Hst(K,V)
\xrightarrow{\partial_N}\mathcal M_N(V)\longrightarrow0
\end{equation}
and
\begin{equation}\label{eq:universal-e-quotient}
0\longrightarrow\mathcal U(V)\longrightarrow
\Hst(K,V)/\He(K,V)\xrightarrow{\partial_N}\mathcal M_N(V)
\longrightarrow0.
\end{equation}
For a semistable extension $E_x$ and any lift $h\in\Dst(E_x)$ of $1$, the boundary is $[Nh]\in D/ND$. Moreover,
\begin{equation}\label{eq:universal-rankjump}
\operatorname{rank}N_{E_x}=\operatorname{rank}N_D+
\begin{cases}0,&\partial_Nx=0,\\1,&\partial_Nx\ne0.\end{cases}
\end{equation}
No splitting of \eqref{eq:universal-e-quotient} is chosen here.
\end{theorem}
\begin{proof}
Set $a=(1-\varphi)h$ and $b=Nh$. The extension identity is $Na=(1-p\varphi)b$. Thus $[b]\in\mathcal M_N(V)$, and changing $h$ changes $b$ by an element of $ND$. Conversely, for $[b]\in\mathcal M_N(V)$, the vector $(1-p\varphi)b$ lies in $ND$. Choose $a\in D$ with $Na=(1-p\varphi)b$. Then $(a,b,0)$ is a cocycle of \eqref{eq:cst}. By the Ext comparison of \cite[Theorem 1.1.4]{BLZ16}, it gives an extension of weakly admissible filtered modules, hence a semistable representation. This proves surjectivity without an additional admissibility assumption on the chosen Hodge entry.

The class belongs to $\Hf$ precisely when $1$ has an invariant lift over crystalline periods, equivalently a lift in $\Dcris(E_x)$. For semistable $E_x$ this is $\Dst(E_x)^{N=0}$. Such a lift exists exactly when $b$ can be changed to zero by adding $ND$. This proves \eqref{eq:universal-boundary}. Notice that this condition does not require $E_x$ to be crystalline when $V$ itself has monodromy.

For zero boundary, adjust $h$ so that $Nh=0$. Then $a\in D^{N=0}$, and further changes preserving this normalization come from $D^{N=0}$. The resulting finite-group presentation is
\[
\Hf(K,V)\simeq
\frac{D^{N=0}\oplus t_D}
 {\{((1-\varphi)v,\bar v):v\in D^{N=0}\}}.
\]
Quotienting by the Hodge-only classes gives $\Hf/\He=\mathcal U(V)$ and proves \eqref{eq:universal-e-quotient}. Finally, the image of $N_{E_x}$ is $ND+K_0b$. It has dimension $\operatorname{rank}N_D$ or one more, according to whether $[b]$ vanishes. This proves \eqref{eq:universal-rankjump}.
\end{proof}

\begin{proposition}[Dual obstruction spaces]\label{prop:obstruction-duality}
Put $V^\dagger=V^\vee(1)$. The period evaluation and $\Tr_{K_0/\Qp}$ induce a perfect pairing of $\Qp$-spaces
\begin{equation}\label{eq:obstruction-duality}
\mathcal M_N(V)\times\mathcal U(V^\dagger)\longrightarrow\Qp.
\end{equation}
In this assertion the coefficient of the standard basis of $\Dst(\Qp(1))$ is used before taking trace. This is an algebraic pairing of period spaces; the arithmetic normalization of local Tate duality remains \eqref{eq:reciprocity-convention}.
\end{proposition}
\begin{proof}
Write the scalar coefficient of period evaluation as $b_D(x,y)$. The relation $b_D(Nx,y)+b_D(x,N^\dagger y)=0$ identifies the dual of $D/ND$ with $(D^\dagger)^{N^\dagger=0}$. Frobenius satisfies
\[
b_D(\varphi x,\varphi^\dagger y)=p^{-1}\sigma(b_D(x,y)).
\]
After trace, the adjoint of $\varphi^\dagger$ is $(p\bar\varphi)^{-1}$. The annihilator of $(1-\varphi^\dagger)(D^\dagger)^{N^\dagger=0}$ is therefore $\ker(1-p\bar\varphi)$. Nondegeneracy of the original evaluation and of the field trace proves perfectness of \eqref{eq:obstruction-duality}.
\end{proof}

\begin{corollary}[Arithmetic duality of the universal quotients]\label{cor:universalquotientduality}
For any semistable $V$ over a finite extension $K/\Qp$, one has
\begin{equation}\label{eq:universal-geometric-condition}
\Hst(K,V)=\Hg(K,V).
\end{equation}
The local Tate pairing induces a perfect pairing
\begin{equation}\label{eq:universal-quotient-pairing}
\frac{\Hst(K,V)}{\He(K,V)}\ \times\
\frac{\Hst(K,V^\dagger)}{\He(K,V^\dagger)}
\longrightarrow\Qp.
\end{equation}
The subspaces $\Hf/\He=\mathcal U$ on the two sides are exact annihilators.
Consequently, this pairing induces perfect pairings between
$\mathcal U(V)$ and $\mathcal M_N(V^\dagger)$, and between
$\mathcal M_N(V)$ and $\mathcal U(V^\dagger)$, without choosing a splitting
of \eqref{eq:universal-e-quotient}.
\end{corollary}
\begin{proof}
The period splitting criterion for $\Hg$ says that the middle representation
of the extension is de Rham, since $V$ is de Rham
\cite[the remark following (3.7.4)]{BK90}. Berger's extension theorem
\cite[Theorem 6.2]{Ber02} says that a de Rham extension of two semistable
representations is semistable. This proves
\eqref{eq:universal-geometric-condition}; the reverse inclusion follows
from semistable comparison.

By Bloch--Kato orthogonality \cite[Proposition 3.8]{BK90}, the annihilator
of $\Hg(K,V^\dagger)$ in full local cohomology is $\He(K,V)$.
Since $\He\subset\Hst$, this is precisely the radical of the restricted
pairing. The same statement on the dual side proves
\eqref{eq:universal-quotient-pairing}. The finite subgroups are exact
annihilators in full cohomology by the same proposition. As both finite
subgroups lie in the respective semistable subgroups, their images in
the quotients remain exact annihilators. Use
\eqref{eq:universal-e-quotient} to identify those images with $\mathcal U$
and the remaining quotients with $\mathcal M_N$. Perfectness of the two
induced pairings follows from elementary duality of finite-dimensional
vector spaces. These arithmetic pairings and the algebraic pairing of
\Cref{prop:obstruction-duality} have separate definitions; their numerical
normalizations are not identified by this argument.
\end{proof}

\subsection{The canonical quotient for a mixed \texorpdfstring{$1$}{1}-motive}

We now return to the unramified field and the geometric hypotheses of \Cref{thm:mixedrealization}. They include the toric case $g=0$. The extra structure gives a canonical splitting at the level of the exponential quotient that is not part of \Cref{thm:universalboundary} for an arbitrary representation.

\begin{theorem}[Monodromy and exponential coordinates]\label{thm:boundary}
For $V=V_M$, the spaces in \eqref{eq:general-obstruction-spaces} are canonically identified, in the given frames, with
\[
\mathcal M_N(V)=\coker R,\qquad \mathcal U(V)=\ker R.
\]
In particular,
\begin{equation}\label{eq:boundaryexact}
0\longrightarrow\Hf(K,V)\longrightarrow\Hst(K,V)
\xrightarrow{\partial_N}\coker R\longrightarrow0,
\end{equation}
and
\begin{equation}\label{eq:mixed-he}
\He(K,V)\simeq\mathfrak g/\lambda(\ker R).
\end{equation}
The quotient $V\to Y\otimes\Qp$ defines a map
$\alpha_V:\Hst(K,V)\to\ker R$: its quotient character is $\lambda_0\alpha_V$. Together these maps induce
\begin{equation}\label{eq:canonical-quotient}
\Hst(K,V)/\He(K,V)\xrightarrow{\ (\alpha_V,\partial_N)\ }\ker R\oplus\coker R.
\end{equation}
This is an isomorphism, independent of $\vartheta$. For $h,d$ as in \Cref{thm:mixednormalform},
\begin{equation}\label{eq:mixed-local-dimensions}
\dim\He=n(a+g)-s+d+h,\qquad \dim\Hf=n(a+g)+h.
\end{equation}
\end{theorem}
\begin{proof}
In \eqref{eq:mixed-D}, quotienting by $ND=RK^s$ leaves toric Frobenius $p^{-1}\sigma$ on $K\otimes\coker R$, the middle block $\Phi_A$, and the lattice block $\sigma$. Only the first has $p\varphi$-fixed vectors, by \Cref{thm:mixedrealization,lem:frobenius}. Its fixed space is $\coker R$. Similarly, on $D^{N=0}=D_G\oplus K\ker R$, the operator $1-\varphi$ is bijective on $D_G$ and has cokernel $\ker R$ on the final summand, through trace. This identifies the two obstruction spaces. The boundary of a reduced cocycle $(u,b,z)$ is the fixed class $[b]$ modulo $R$; equation $Ru=(1-\sigma)b$ gives fixedness. Thus \eqref{eq:boundaryexact} is \eqref{eq:universal-boundary} in these coordinates.

In the reduced complex, Hodge-only cocycles are identified precisely by $z\mapsto z+\lambda y$ for $y\in\ker R\subset\Qp^s$. This proves \eqref{eq:mixed-he}. Projection to the lattice quotient sends $(u,b,z)$ to the unramified extension represented by $u$. Applying the convention of \Cref{lem:tau} gives
\[
\alpha_V=\epsilon^{-1}\Tr(u)\in\ker R.
\]
The normal form of \Cref{thm:threeblockcohomology} shows that $(\alpha_V,\partial_N)$ is surjective. If both coordinates vanish, take $u=0$ in normal form and write $b=Ry$ with $y\in\Qp^s$. Subtracting $\delta^0y$ leaves a Hodge-only class. This proves the kernel assertion and \eqref{eq:canonical-quotient}. Finally $\dim_{\Qp}\lambda(\ker R)=s-d-h$, giving the first dimension formula. The second follows from $\Hf/\He=\ker R$.
\end{proof}

For $G=\mathbf G_m^a$, the identification $\mathfrak g=K^a$ gives $\lambda=-L$. Thus \eqref{eq:mixed-he} specializes to the formula used in the toric calculations below,
\begin{equation}\label{eq:he-description}
\He(K,V_Q)\simeq K^a/L(\ker R).
\end{equation}
The general formula identifies the Hodge contribution with formal logarithms of the integral points $u_0(\ker R)$. By \Cref{cor:mixedbranch}, this subspace is unchanged by a change of uniformizer and the associated logarithm branch.

\begin{corollary}\label{cor:rankjump}
For a semistable extension $E_x$ of $\Qp$ by $V_M$, put $d=\operatorname{rank}R$. Then $N_{E_x}^2=0$ and
\[
\operatorname{rank}N_{E_x}=\begin{cases}d,&\partial_Nx=0,\\d+1,&\partial_Nx\ne0.\end{cases}
\]
Equivalently, $\dim_K\Dcris(E_x)-\dim_K\Dcris(V_M)$ is $1$ for a finite class and $0$ for a nonfinite semistable class.
\end{corollary}
\begin{proof}
The rank formula is \eqref{eq:universal-rankjump}. For square-zero monodromy, use a lift in \eqref{eq:extension-triple} and the splitting \eqref{eq:mixed-D}. Projecting $Na=(1-p\varphi)b$ to $D_A$ and $K^s$ forces $b_A=b_Y=0$, because $1-p\Phi_A$ and $1-p\sigma$ are bijective. Thus the new monodromy vector belongs to the toric summand, which $N$ kills. This proves $N_{E_x}^2=0$. The assertion about crystalline dimensions follows by taking kernels of $N$.
\end{proof}

\begin{proposition}[The boundary as a positive logarithmic coefficient]\label{prop:boundarypositive}
Let $x\in\Hst(K,V_M)$ under the hypotheses of \Cref{thm:mixedrealization}, and let $h\in\Dst(E_x)$ lift $1$. Put $b=Nh\in\Dst(V_M)$. Then $b$ belongs to the toric summand and there is $A_h\in\Rplus\otimes E_x$ such that
\begin{equation}\label{eq:boundary-positive-representative}
t h=A_h-\ell_p\,t b,
\qquad t b\in\Rplus\otimes V_M.
\end{equation}
The class $[b]\in D/ND$ is independent of the lift and equals $\partial_Nx$. In particular, the following conditions are equivalent:
\begin{enumerate}
\item $x\in\Hf(K,V_M)$;
\item there is a period lift $h$ of $1$ with $Nh=0$;
\item there is a period lift $h$ of $1$ with $t h\in\Rplus\otimes E_x$.
\end{enumerate}
Thus the nonremovable logarithmic coefficient, measured modulo the old monodromy image, is exactly the obstruction in \eqref{eq:boundaryexact}.
\end{proposition}
\begin{proof}
The proof of \Cref{cor:rankjump} shows that $b$ is toric and $Nb=0$. By \Cref{prop:generalpositive}, $t h\in\Rlog\otimes E_x$ and $t b\in\Rplus\otimes V_M$. On this ring $N=-\partial_{\ell_p}$ and $N(t)=0$. Therefore
\[
N(t h+\ell_p t b)=t b-t b=0.
\]
The kernel of $-\partial_{\ell_p}$ on $\Rplus[\ell_p]$ is $\Rplus$, which gives \eqref{eq:boundary-positive-representative}. Changing $h$ by $v\in D$ changes $b$ by $Nv$; \Cref{thm:universalboundary} identifies its class with the boundary and proves equivalence of (i) and (ii). Condition (ii) puts $h$ in $\Dcris(E_x)$, so \Cref{prop:generalpositive} gives (iii). Conversely (iii) implies $0=N(t h)=tNh$, hence (ii). This proves all assertions.
\end{proof}

\subsection{Cartier duality and the local pairing}

The Cartier dual $M^\vee$ has abelian quotient $A^\vee$, torus $\Hom(Y,\mathbf G_m)$, and lattice $X^*(T)$. Its finite realizations pair perfectly with those of $M$ into roots of unity; hence $V_p(M^\vee)=V_M^\dagger$. The construction uses the Poincar\'e biextension of $A\times A^\vee$, as in \cite[Section 1]{HS05} and \cite[Section 4.3]{Ray94}. Its semi-abelian part has a model of the form \eqref{eq:good-model}: the classifying points of $A^\vee$ and $A$ extend over $\OK$ by properness. The preceding theorems therefore apply to both motives.

Choose dual frames so that coefficient evaluation is the usual positive evaluation between the toric part on one side and the lattice quotient on the other. Monodromy compatibility then gives
\begin{equation}\label{eq:mixed-dual-R}
R^\dagger=-R^{\mathsf t}.
\end{equation}
This specifies the sign in the framed Tate realization; a different lattice frame for the geometric Cartier dual changes its displayed matrix correspondingly.

\begin{theorem}[Semistable reciprocity for $1$-motives]\label{thm:semistablereciprocity}
Let $M$ satisfy \eqref{eq:mixed-motive}--\eqref{eq:good-model}, and use the arithmetic normalization \eqref{eq:reciprocity-convention}. For $x\in\Hst(K,V_M)$ and $y\in\Hst(K,V_M^\dagger)$,
\begin{equation}\label{eq:semistable-reciprocity}
\boxed{\quad
\pair{x}{y}=\alpha_{V_M^\dagger}(y)(\partial_Nx)
             -\alpha_{V_M}(x)(\partial_{N^\dagger}y).
\quad}
\end{equation}
The evaluations use $(\coker R)^*=\ker R^{\mathsf t}$ and $(\coker R^{\mathsf t})^*=\ker R$. The radicals of this restricted pairing are exactly $\He(K,V_M)$ and $\He(K,V_M^\dagger)$. It induces a perfect pairing on the quotients in \eqref{eq:canonical-quotient}, of dimension $a+s-2d$ on either side.
\end{theorem}
\begin{proof}
Every boundary class in $\coker R$ is represented by a Kummer class of the toric subrepresentation, using a tuple of powers of $p$. Subtract such a class from $x$ to obtain an element of $\Hf$, by \Cref{thm:boundary}; do the same on the dual side. Evaluation between the two toric subrepresentations is zero: the torus of the Cartier dual consists of functionals on the lattice quotient and annihilates $V_p(G)$. Thus the two toric Kummer classes pair to zero. The two finite classes pair to zero by Bloch--Kato orthogonality \cite[Proposition 3.8]{BK90}.

If $x$ is the toric Kummer image of $a'=(a'_i)$ and $y$ is finite, functoriality of cup products projects $y$ to its lattice-quotient character $\lambda_0\alpha_{V_M^\dagger}(y)$. By \eqref{eq:reciprocity-convention},
\[
\pair{x}{y}=\sum_i\alpha_{V_M^\dagger}(y)_i\ord_p(a'_i)
            =\alpha_{V_M^\dagger}(y)(\partial_Nx).
\]
The expression factors through $\coker R$ because the character belongs to $\ker R^{\mathsf t}$. Reversing the arguments contributes the minus sign from graded commutativity in degree one. These four terms prove \eqref{eq:semistable-reciprocity}. Under \eqref{eq:canonical-quotient} and its dual version, the formula is the sum of two perfect evaluation pairings, one with negative sign. Its radicals and the dimension assertion follow.
\end{proof}

\begin{remark}\label{rem:meaning-reciprocity}
The theorem expresses the restriction of classical Tate duality using geometric monodromy and the unramified quotient. The good abelian part enlarges the exponential space through $\mathfrak g$ and $\lambda$. It creates no additional invariant summand in \eqref{eq:semistable-reciprocity}. The filtered abelian module and its Poincar\'e pairing remain necessary for the representation and for full local cohomology. In the toric case, \Cref{thm:actualperiods} supplies actual representatives for both arguments of this formula.
\end{remark}

\begin{example}[Monodromy created by an extension]\label{ex:unitextension}
Take $K=\Qp$, $a=s=1$, and $q=1+p$. Then $R=0$, $L=\log_p(1+p)\ne0$, and $V_q$ is crystalline. The long exact Kummer sequence gives $\dim H^1(K,V_q)=2$, while \Cref{thm:normalform} gives the same dimension for $\Hst$. Let $x_p$ be the Kummer image of $p$, and let $x_\ur$ be the class constructed with $\alpha=1$ and $a'=0$. They form a basis: their pairs $(\alpha_V,\partial_N)$ are $(0,1)$ and $(1,0)$. The period lift for $x_p$ is
\[
h_p=z_p-\ell_p e,\qquad Nh_p=e.
\]
Thus its middle extension is semistable with nonzero monodromy. The period lift for $x_\ur$ is
\[
h_\ur=z_\ur-\tau f-\rho e,\qquad
\theta\rho=\bar\tau L,\qquad Nh_\ur=0.
\]
Its middle extension is crystalline. For the dual parameter $q^{-1}$, choose the analogous classes $y_p,y_\ur$. \Cref{thm:semistablereciprocity} gives the full local pairing matrix
\[
\left(\pair{x_i}{y_j}\right)_{i,j\in\{p,\ur\}}
=\begin{pmatrix}0&1\\-1&0\end{pmatrix}.
\]
The coefficient representation has $N_D=0$, but the class $x_p$ has nonzero monodromy boundary. This distinguishes the monodromy of an extension from the coefficient-monodromy coordinate term in \eqref{eq:ZplusC} below.
\end{example}

\section{Abelian and mixed examples, and the full lifting obstruction}\label{sec:mixedexamples}

\subsection{A good abelian variety}

\begin{example}[The abelian exponential space]\label{ex:goodabelian}
Let $M=[0\to A]$, where $A/K$ has good reduction and dimension $g>0$. There are no outer weight blocks: $a=s=0$. By \Cref{thm:mixednormalform,thm:boundary},
\begin{equation}\label{eq:pure-abelian-local}
\He(K,V_p(A))=\Hf(K,V_p(A))=\Hst(K,V_p(A))
\simeq\operatorname{Lie}(A),\qquad \dim_{\Qp}=ng.
\end{equation}
The same holds for $A^\vee$. The restricted semistable pairing is identically zero, by \Cref{thm:semistablereciprocity}. In contrast,
\begin{equation}\label{eq:pure-abelian-full}
\dim_{\Qp}H^1(K,V_p(A))=2ng,
\end{equation}
and full local Tate duality pairs this space perfectly with $H^1(K,V_p(A^\vee))$.

To check the full dimension, $A(K)[p^\infty]$ is finite: a sufficiently small formal subgroup is torsion free and has finite index in the compact group $A(K)$. Thus $H^0(K,V_p(A))=0$, and the same argument for $A^\vee$, followed by Tate duality, gives $H^2(K,V_p(A))=0$. The local Euler characteristic formula \cite[Chapter I, Section 2]{Milne06} now gives \eqref{eq:pure-abelian-full}. This example shows why a semistable comparison alone does not compute full cohomology, even in good reduction.
\end{example}

\subsection{A nonsplit extension of an elliptic curve}

\begin{example}[Integral points and a monodromy change]\label{ex:mixedelliptic}
Let $E/K$ be a good-reduction elliptic curve. Choose a nontorsion point $P\in E^\vee(K)$ and let
\[
0\longrightarrow\mathbf G_m\longrightarrow G_P\longrightarrow E\longrightarrow0
\]
be the extension classified by $P$. The classification is the usual identification of $\operatorname{Ext}^1(E,\mathbf G_m)$ with $E^\vee(K)$ through the Poincar\'e bundle; see \cite[Section 1]{HS05}. The point $P$ extends over $\OK$, so this extension has a model as in \eqref{eq:good-model}. Choose $Q_0\in\mathcal G_P(\OK)$ whose image in $E(K)$ is nontorsion, and put
\[
M_m=[\mathbf Z\longrightarrow G_P],\qquad
1\longmapsto Q_0+i(p^m),\qquad m\in\mathbf Z.
\]
Then $a=s=g=1$ and \Cref{thm:mixedrealization} gives the concrete data
\begin{equation}\label{eq:elliptic-point-data}
R=(m),\qquad \mathfrak g=\operatorname{Lie}(G_P),\qquad
\lambda(1)=-\log_{\mathcal G_P}(Q_0).
\end{equation}
The last vector is nonzero, since its projection to $\operatorname{Lie}(E)$ is the logarithm of a nontorsion point. Hence $H^0(K,V_{M_m})=0$. The lattice generator of the Cartier dual projects to $P$ in $E^\vee$; its logarithm is nonzero. The same argument gives $H^0(K,V_{M_m}^\dagger)=0$. Thus both full degree-one spaces have dimension $4n$.

The semistable local conditions are
\begin{equation}\label{eq:elliptic-local-table}
\begin{array}{c|cc}
 &m=0&m\ne0\\ \hline
\operatorname{rank}N_{V_{M_m}}&0&1\\
\dim\Hst&2n+1&2n\\
\dim\Hf&2n&2n\\
\dim\He&2n-1&2n\\
\dim(\Hst/\He)&2&0.
\end{array}
\end{equation}
These entries follow directly from \Cref{thm:mixednormalform,thm:boundary}. For $m=0$, the coefficient representation is crystalline, while a semistable extension class with nonzero boundary creates rank-one monodromy. For $m\ne0$, every semistable class is finite and exponential, and the restricted semistable pairing is zero. In both columns the full $4n$-dimensional Tate pairing remains perfect. The distinction is rational: a nonzero integer $m$ is invertible over $\Qp$ even if $p\mid m$; integral effects of that divisibility are not suppressed by \eqref{eq:elliptic-local-table}.

For a concrete good ordinary curve take $K=\mathbf Q_5$ and $E:y^2=x^3+x+1$. Its discriminant is $-16\cdot31$, a $5$-adic unit, and a direct count gives $\#E(\mathbf F_5)=9$, so its Frobenius trace is $-3$ and the reduction is ordinary. Use its principal polarization to view $P$ on $E^\vee$. Choose $P=\exp_E(5^c)$ and $Q_0=\exp_{\mathcal G_P}(5^c v)$ for $c$ sufficiently large, where an integral tangent vector $v$ projects to $1\in\operatorname{Lie}(E)$ in fixed local coordinates. Formal exponentials converge for such $c$. These choices give
\[
\lambda(1)=-5^c v,
\qquad
(\dim H^1,\dim\Hst,\dim\Hf,\dim\He)
=\begin{cases}(4,3,2,1),&m=0,\\(4,2,2,2),&m\ne0.\end{cases}
\]
This is a nonsplit semi-abelian example: nontorsion of $P$ prevents its extension class from becoming split after an isogeny.
\end{example}

\subsection{Split products and actual division-point representatives}

\begin{proposition}[Adding a good abelian factor]\label{prop:splitproducts}
Let $M=M_Q\oplus[0\to A]$, with $A$ of good reduction. For $\star\in\{e,f,\mathrm{st}\}$ and for full local cohomology there are direct-sum decompositions
\begin{equation}\label{eq:splitproduct-cohom}
H^1_\star(K,V_M)=H^1_\star(K,V_Q)\oplus H^1_\star(K,V_p(A)),
\qquad
H^1(K,V_M)=H^1(K,V_Q)\oplus H^1(K,V_p(A)).
\end{equation}
Tate duality is the orthogonal direct sum of the toric and abelian pairings. The toric period representatives of \Cref{thm:actualperiods} combine with classical crystalline representatives for the abelian summand.

More concretely, every abelian semistable class is represented by the following Galois cocycle. For $z\in\operatorname{Lie}(A)$ choose $c\geq0$ so that $P=\exp_A(p^cz)$ is defined, and choose points $P_j\in A(\overline K)$ with $P_0=P$ and $[p]P_{j+1}=P_j$. Then
\begin{equation}\label{eq:abelian-division-cocycle}
g\longmapsto p^{-c}(gP_j-P_j)_{j\geq0}\ \in V_p(A)
\end{equation}
represents the class $\exp_{V_p(A)}(z)$. It is independent of these choices in cohomology.
\end{proposition}
\begin{proof}
The coefficient representation and the semistable complex are direct sums. The defining kernels for the local conditions commute with finite direct sums, proving \eqref{eq:splitproduct-cohom}. The evaluation between a representation and the dual of the other factor is zero, proving orthogonality. The abelian class in \eqref{eq:abelian-division-cocycle} is the Kummer class of $P$, divided by $p^c$. The formal-group description of Bloch--Kato's exponential \cite[Section V.2]{BerIntro02} identifies it with $\exp_{V_p(A)}(z)$. This exponential is an isomorphism by \eqref{eq:pure-abelian-local} and $\ker(1-\Phi_A)=0$. Thus these cocycles give all abelian semistable classes. Their middle representations are crystalline, by the good $p$-divisible construction in \Cref{thm:mixedrealization}. Changing division points adds a coboundary; changing $c$ gives the same exponential class.
\end{proof}

\begin{corollary}[An integral product extension of trace duality]\label{cor:integralproduct}
Let $T_Q$ be the toric basis-family lattice of \Cref{thm:integraltrace} and set $T=T_Q\oplus T_p(A)$. Assume
\[
A(K)[p^\infty]=A^\vee(K)[p^\infty]=0.
\]
Then $H^1(K,T_p(A))$ and $H^1(K,T_p(A^\vee))$ are free of rank $2ng$ and pair perfectly over $\Zp$. There are identifications
\begin{align}\label{eq:integralproduct-trace}
H^1(K,T)&=\pgl_{\Zp}(W)\oplus H^1(K,T_p(A)),\\
H^1(K,T^\dagger)&=\slie_{\Zp}(W)\oplus H^1(K,T_p(A^\vee)),\notag
\end{align}
under which
\[
\pair{([B],x)}{(C,y)}=\tr(CB)+\langle x,y\rangle_{K,A}.
\]
If a polarization of degree prime to $p$ is used to identify the two abelian Tate lattices, the trace-zero-to-projective map on the toric factor and the identity on the abelian factor have cokernel $\Zp/r\Zp$. The abelian factor adds no lattice defect to that map.
\end{corollary}
\begin{proof}
The torsion in $H^1(K,T_p(A))$ is $A(K)[p^\infty]$, by the cohomology sequence for $0\to T_p(A)\to V_p(A)\to A[p^\infty]\to0$ and $H^0(K,V_p(A))=0$. Integral local duality identifies $H^2(K,T_p(A))$ with the dual of $A^\vee(K)[p^\infty]$. Both groups therefore vanish under the hypotheses. The same statements hold after interchanging $A$ and $A^\vee$. For every $j$, reduction identifies $H^1(K,T_p(A))/p^j$ with $H^1(K,A[p^j])$; this follows from the vanishing of $H^2(K,T_p(A))[p^j]$. Finite local duality gives perfect pairings on these reductions, so the integral pairing is perfect. Its rank is \eqref{eq:pure-abelian-full}. Apply the direct-sum decomposition and \Cref{thm:integraltrace}; the final cokernel is \Cref{prop:tracedefect} on the toric factor.
\end{proof}

\subsection{The obstruction replacing a trace-zero equation}

The preceding product formula is a direct sum. For a nonsplit mixed motive, full dual cohomology instead has a relative lifting description. Put $U=V_p(G)$, $Y_{\Qp}=Y\otimes\Qp$, and let
\[
\xi_M:Y_{\Qp}\longrightarrow H^1(K,U),\qquad
y\longmapsto\kappa_G(u(y))
\]
be its Kummer extension map. Let $U^\dagger=U^\vee(1)$ and $Y_{\Qp}^*=\Hom_{\Qp}(Y_{\Qp},\Qp)$. Define
\begin{equation}\label{eq:mixed-full-obstruction}
\delta_M:H^1(K,U^\dagger)\longrightarrow Y_{\Qp}^*,
\qquad
\delta_M(z)(y)=\inv_K(\xi_M(y)\cup z).
\end{equation}
Thus the obstruction is the adjoint of the Kummer map for full local Tate duality, including its abelian pairing.

\begin{proposition}[Full dual lifting criterion]\label{prop:fullmixedobstruction}
A class $z\in H^1(K,U^\dagger)$ lifts to $H^1(K,V_M^\dagger)$ if and only if $\delta_M(z)=0$. There is an exact sequence
\begin{equation}\label{eq:mixed-dual-exact}
0\longrightarrow
\frac{H^1(K,Y_{\Qp}^*(1))}{\operatorname{im}H^0(K,U^\dagger)}
\longrightarrow H^1(K,V_M^\dagger)
\longrightarrow\ker\delta_M\longrightarrow0,
\end{equation}
where the denominator denotes the connecting image. Furthermore,
\begin{equation}\label{eq:mixed-obstruction-rank}
\operatorname{rank}\delta_M=s-\dim H^0(K,V_M)=s-h.
\end{equation}
The map $\delta_M$ is minus the connecting map for the framed dual extension, after identifying $H^2(K,Y_{\Qp}^*(1))$ with $Y_{\Qp}^*$ by local invariants.
\end{proposition}
\begin{proof}
Dualizing $0\to U\to V_M\to Y_{\Qp}\to0$ gives
\[
0\longrightarrow Y_{\Qp}^*(1)\longrightarrow V_M^\dagger
\longrightarrow U^\dagger\longrightarrow0.
\]
Its extension class is minus the transpose of the original one in the evaluation frames used in \eqref{eq:mixed-dual-R}. The connecting map in degree one is cup product with this dual extension class. Evaluation on $y$ is consequently $-\inv_K(\xi_M(y)\cup z)$. The long exact sequence proves the lifting criterion and \eqref{eq:mixed-dual-exact}. To obtain the rank, note that $H^0(K,U)=0$: the $p$-primary torsion in $G(K)$ is finite, since this holds for its torus and its abelian quotient. The original long exact sequence therefore identifies $\ker\xi_M$ with $H^0(K,V_M)$. Since full local Tate duality for $U$ is perfect, a map and its adjoint have equal rank. This gives \eqref{eq:mixed-obstruction-rank}.
\end{proof}

For the toric basis family, $s=1$, $U=\Qp(1)^r$, and $\xi_M(1)=(\kappa(q_i))_i$. Equation \eqref{eq:mixed-full-obstruction} then becomes $B\mapsto\tr(B)$ in the matrix coordinates of \Cref{thm:integraltrace}. In a general mixed motive it is a vector of $s$ cup products with the geometric Kummer classes. \Cref{thm:relativeHerrlift} below turns this criterion into actual dual Herr representatives by solving the resulting Tate-valued equations component by component.

\section{An integral reduced complex}\label{sec:integralcomplex}

The following complex is defined over $\Zp$ without dividing trace by $n$:
\begin{equation}\label{eq:integralcomplex}
\mathcal C_Q^\bullet:\quad
\OK^s\xrightarrow{\delta^0}
\OK^s\oplus\OK^a\oplus p\OK^a
\xrightarrow{\delta^1}\OK^a,
\end{equation}
with the same formulas as \eqref{eq:reducedcomplex}. The use of $p\OK^a$ in the Hodge coordinate is compatible with $L\in M_{a,s}(p\OK)$. Since $R$ has rational integer entries and commutes with $\sigma$, this is a complex.

\begin{proposition}\label{prop:integralcohom}
Its cohomology is described by
\[
H^0(\mathcal C_Q)=\ker(R,L)\cap\Zp^s,\qquad
H^2(\mathcal C_Q)=\Zp^a/R\Zp^s,
\]
and the exact sequence
\begin{equation}\label{eq:integralnormal}
0\longrightarrow
\frac{\Zp^a\oplus p\OK^a}{(R,-L)\Zp^s}
\longrightarrow H^1(\mathcal C_Q)
\xrightarrow{\Tr(u)}\ker(R:\Zp^s\to\Zp^a)
\longrightarrow0.
\end{equation}
A trace-one $\vartheta\in\OK$ splits this sequence. Rationalization of \eqref{eq:integralcomplex} is quasi-isomorphic to $C_{\mathrm{st}}^\bullet(\Dst(V_Q))$.
\end{proposition}
\begin{proof}
The proof of \Cref{thm:normalform} now applies integrally using \Cref{lem:frobenius}. For completeness, the cocycle equation implies $R\Tr(u)=0$. Subtract $\delta^0y$ with $(1-\sigma)y=u-\vartheta\Tr(u)$, where an integral solution exists by that lemma. The second coordinate becomes $\sigma$-fixed and hence belongs to $\Zp^a$. Remaining gauge changes use $y\in\Zp^s$, giving the stated quotient. The splitting is $(\vartheta\alpha,0,0)$. In degree two, trace identifies the quotient by $(1-\sigma)\OK^a$ with $\Zp^a$, and the image of $R\OK^s$ with $R\Zp^s$ because trace is surjective. The degree-zero assertion is immediate. Rationalization follows from the explicit contraction in \Cref{thm:normalform}.
\end{proof}

This proposition describes the indicated integral complex for every $Q$. An identification with the integral Galois lattice requires additional arithmetic input. We prove such an identification for the basis family in \Cref{thm:integralcomparison}; no general integral comparison theorem is inferred just from rationalization.

\section{The basis family and the full integral Tate pairing}\label{sec:tracefamily}

Let $r=n+1$ and fix a basis $q_0=p,q_1,\ldots,q_n$ of $W$ such that $q_i\in1+p\OK$ for $i>0$. We index all $r\times r$ matrices from $0$ to $n$. Consider
\[
M=[\mathbf Z\to\mathbf G_m^r],\quad1\mapsto(q_0,\ldots,q_n),
\qquad T=T_p(M),\quad V=T\otimes\Qp.
\]
Here $R=e_0$ is a column and $L=(0,L_1,\ldots,L_n)^{\mathsf t}$ with $L_i=\log_p q_i$. Define
\[
\pgl_O(W)=\End_O(W)/O I,\qquad
\slie_O(W)=\ker(\tr:\End_O(W)\to O),\quad O=\Zp.
\]

\begin{theorem}[Full integral trace duality]\label{thm:integraltrace}
There are $Q$-dependent integral isomorphisms
\begin{equation}\label{eq:trace-identifications}
H^1(K,T)\simeq\pgl_O(W),\qquad
H^1(K,T^\dagger)\simeq\slie_O(W).
\end{equation}
Moreover $H^0(K,T)=H^2(K,T)=H^0(K,T^\dagger)=H^2(K,T^\dagger)=0$. For $A\in\End_O(W)$ and $B\in\slie_O(W)$,
\begin{equation}\label{eq:integraltrace}
\pair{[A]}{B}=\tr(BA).
\end{equation}
The pairing is perfect over $O$.
\end{theorem}
\begin{proof}
The long exact sequence for $0\to O(1)^r\to T\to O\to0$ begins with the injective map $O\to W^r$, $1\mapsto(q_i)$. Its next relevant boundary is
\[
W^*\longrightarrow O^r,\qquad\lambda\longmapsto(\lambda(q_i))_i,
\]
up to a harmless common sign from cup order. It is an isomorphism because $(q_i)$ is a basis of $W$. Consequently $H^1(K,T)=W^r/O(q_i)$ and $H^0(K,T)=H^2(K,T)=0$. Send a tuple $(a_i)$ to the endomorphism $A$ with $A(q_i)=a_i$. The tuple $(q_i)$ is $I$, which gives the first identification.

The dual extension is $0\to O(1)\to T^\dagger\to O^r\to0$ with Kummer row $-(q_i)^{\mathsf t}$. Its boundary $O^r\to W$ is an isomorphism. The next boundary on $(W^*)^r$ is, up to a global sign, $(\lambda_i)\mapsto\sum_i\lambda_i(q_i)$. It is surjective. Thus the dual $H^0,H^2$ vanish and
\[
H^1(K,T^\dagger)=\{(\lambda_i)\in(W^*)^r:\ \sum_i\lambda_i(q_i)=0\}.
\]
Set $B=\sum_iq_i\otimes\lambda_i$. Its trace is the displayed sum, so this gives the second identification.

Functoriality of cup product, followed by \eqref{eq:reciprocity-convention}, gives
\[
\pair{[A]}{B}=\sum_i\lambda_i(Aq_i)=\tr(BA).
\]
This is unchanged when a scalar is added to $A$.

The trace form on $\End_O(W)$ is perfect, with matrix-unit bases $E_{ij}$ and $E_{ji}$. The scalar line $OI$ is primitive, and its annihilator is exactly $\slie_O(W)$. Dualizing the split sequence of free modules $0\to OI\to\End_O(W)\to\pgl_O(W)\to0$ proves perfectness of the quotient--annihilator pairing. In particular no division by $r$ is used.
\end{proof}

\begin{theorem}[Integral valuation--logarithm comparison]\label{thm:integralcomparison}
For this basis family, the following explicit map is an isomorphism:
\begin{equation}\label{eq:integralcomparison}
H^1(K,T)\xrightarrow{\sim}H^1(\mathcal C_Q),\qquad
[(a_i)]\longmapsto[(0,(\ord_p a_i)_i,-(\log_p a_i)_i)].
\end{equation}
It agrees after rationalization with the period comparison of \Cref{thm:actualperiods}. Every rational class in $H^1(K,V)$ is semistable.
\end{theorem}
\begin{proof}
The column $R=e_0$ is injective over $O$, so \Cref{prop:integralcohom} gives
\[
H^1(\mathcal C_Q)=(O^r\oplus p\OK^r)/O(R,-L).
\]
Equation \eqref{eq:val-log} and the first long exact sequence calculation in \Cref{thm:integraltrace} give exactly \eqref{eq:integralcomparison}. For integral $A$, choose $a_i=\sum_k A_{ki}q_k\in W$ and $\omega_i=\sum_k A_{ki}\ell(q_k)$. The Galois extension and invariant period
\[
h_A=z_A-\sum_i\omega_i e_i
\]
are \eqref{eq:actualextension} and \eqref{eq:actualh} with $\alpha=0$. Their complex representative is the right side of \eqref{eq:integralcomparison}. This proves compatibility and semistability of every rational class.
\end{proof}

Changing $A$ to $A+cI$ changes the reduced representative by $\delta^0(cf)$, where $f$ is the original lattice period. This explicitly verifies scalar ambiguity at the period level.

\begin{proposition}[Exact trace-zero defect]\label{prop:tracedefect}
The natural map $\slie_O(W)\to\pgl_O(W)$ fits into
\begin{equation}\label{eq:trace-defect}
0\longrightarrow\slie_O(W)\longrightarrow\pgl_O(W)
\xrightarrow{\ \tr\bmod r\ }O/rO\longrightarrow0.
\end{equation}
The discriminant ideal of the self trace form on $\slie_O(W)$ is $(r)$. Over $\Fp$, its radical is zero if $p\nmid r$ and is $\Fp I$ if $p\mid r$.
\end{proposition}
\begin{proof}
A scalar matrix of trace zero over $O$ is zero, so the first map is injective. Trace modulo $r$ is well defined on the projective quotient and is surjective. If $\tr(A)=rc$, then $A-cI$ is trace zero and represents the same projective class, proving exactness.

Choose the off-diagonal matrix units and $H_i=E_{ii}-E_{nn}$ for $0\leq i<n$ as a basis of $\slie_O(W)$. Off-diagonal units form mutually dual pairs under trace. The diagonal Gram matrix is $I_n+J_n$, where $J_n$ is the all-ones matrix. Its determinant is $n+1=r$, proving the discriminant assertion. Over a field the annihilator of trace-zero matrices inside all matrices is the scalar line. Intersecting that line with trace-zero matrices gives the stated radical.
\end{proof}

The rational operation $[A]\mapsto A-r^{-1}\tr(A)I$ therefore identifies the two rational spaces, but preserves these lattices precisely when $r$ is a unit of $O$. Equation \eqref{eq:trace-defect} is an arithmetic constraint on any proposed integral self-identification; it is not a computation of all regulator denominators.

\section{Filtered trace duality and a concrete example}\label{sec:filteredtrace}

Write $W=Oq_0\oplus U$, where $U$ is the principal-unit summand with basis $q_1,\ldots,q_n$. Every projective matrix has a unique representative
\begin{equation}\label{eq:blockA}
A=\begin{pmatrix}0&\alpha\\ u&D_0\end{pmatrix},
\qquad\alpha\in U^*,\quad u\in U,\quad D_0\in\End_O(U).
\end{equation}
Every trace-zero matrix has a unique form
\begin{equation}\label{eq:blockB}
B=\begin{pmatrix}-\tr(E)&\beta\\ v&E\end{pmatrix},
\qquad\beta\in U^*,\quad v\in U,\quad E\in\End_O(U).
\end{equation}

\begin{theorem}\label{thm:filteredtrace}
With the integral identifications of \Cref{thm:integraltrace},
\begin{equation}\label{eq:blockpair}
\pair{[A]}{B}=\alpha(v)+\beta(u)+\tr(ED_0).
\end{equation}
Define integral local conditions by intersection with the corresponding rational subspaces. Then
\begin{align*}
\Hf(K,T)=\He(K,T)&=\{[A]:\alpha=0\},\\
\Hf(K,T^\dagger)=\Hst(K,T^\dagger)&=
\left\{\begin{pmatrix}0&0\\ v&0\end{pmatrix}:v\in U\right\},
\qquad \He(K,T^\dagger)=0.
\end{align*}
Here $\Hst(K,T^\dagger)$ means the intersection with $\Hst(K,V^\dagger)$. The monodromy boundary is $\partial_N[A]=\alpha$ after the basis identification $\coker R\simeq U^*$. Hence it gives a perfect integral pairing
\begin{equation}\label{eq:boundarypair}
\frac{H^1(K,T)}{\Hf(K,T)}\times\Hf(K,T^\dagger)\longrightarrow O,
\qquad (\alpha,v)\longmapsto\alpha(v).
\end{equation}
\end{theorem}
\begin{proof}
Multiplication of the block matrices gives \eqref{eq:blockpair}. If $A(q_i)=\sum_k A_{ki}q_k$, the valuation vector of the corresponding Kummer tuple is $(A_{0i})_i$. Its image modulo $R O=Oe_0$ is $(A_{01},\ldots,A_{0n})$, namely $\alpha$. Since $\ker R=0$, \Cref{thm:boundary} gives $\He(K,V)=\Hf(K,V)=\ker\partial_N$. This proves the first equality after intersecting integrally.

For the dual matrix $R^\dagger=-e_0^{\mathsf t}$, the map
\[
(R^\dagger,-L^\dagger):\Qp^r\longrightarrow\Qp\oplus K
\]
is an isomorphism, by \eqref{eq:val-log} and the basis property of $q_i$. Thus its Kummer term in \eqref{eq:h1-exact} is zero. Its semistable classes consist precisely of unramified quotient characters $\lambda_i=v_i\lambda_0$ with $v_0=0$. In the trace-zero description these are $B=\sum_{i>0}v_iE_{i0}$. Moreover $\coker R^\dagger=0$, so all these classes are finite. Equation \eqref{eq:he-description} gives $\He(K,V^\dagger)=0$ because $L_1,\ldots,L_n$ span $K$ over $\Qp$. Intersecting with the integral character lattice gives exactly $v\in U$. Finally \eqref{eq:blockpair} on this subspace is $\alpha(v)$, proving perfectness of \eqref{eq:boundarypair}.
\end{proof}

\begin{corollary}[Dimension obstruction]\label{cor:dimensionobstruction}
For this family,
\[
\begin{array}{c|cc}
 & V&V^\dagger\\ \hline
\dim H^1(K,-)&r^2-1&r^2-1\\
\dim\Hst(K,-)&r^2-1&r-1\\
\dim\Hf(K,-)&r(r-1)&r-1\\
\dim\He(K,-)&r(r-1)&0\\
\dim H^2(C_{\mathrm{st}}^\bullet)&r-1&0.
\end{array}
\]
The semistable extension complex of $V^\dagger$ cannot represent all of $H^1(K,V^\dagger)$. In particular the restricted semistable pairing is degenerate on the $V$ side with radical $\Hf(K,V)$, although the full pairing is perfect.
\end{corollary}
\begin{proof}
The full dimensions follow by rationalizing \Cref{thm:integraltrace}. The semistable and local-condition dimensions follow from \Cref{thm:normalform,thm:filteredtrace}. The degree-two entries are cokernels of $R$ and $-R^{\mathsf t}$. Since $r=n+1\geq2$, $r^2-1>r-1$, which proves the obstruction.
\end{proof}

\begin{example}[An integral defect at $p=3$]\label{ex:three}
Let $K=\mathbf Q_3(\zeta)$ with $\zeta^2=2$. This is unramified quadratic since $X^2-2$ is irreducible modulo $3$. Choose
\[
q_0=3,\qquad q_1=4,\qquad q_2=1+3\zeta.
\]
The logarithm series gives
\[
\frac{\log 4}{3}\equiv1\pmod{3\OK},\qquad
\frac{\log(1+3\zeta)}3\equiv\zeta\pmod{3\OK}.
\]
Thus $\log(q_1),\log(q_2)$ are a $\mathbf Z_3$-basis of $3\OK$ by Nakayama's lemma. The cohomology lattices have rank eight. A basis of the projective lattice consists of the six $E_{ij}$ with $i\ne j$, followed by $[E_{00}],[E_{11}]$. Its trace-dual basis is the six transposed units, followed by $E_{00}-E_{22},E_{11}-E_{22}$. The pairing matrix for these two bases is the identity. In contrast, the self trace form on $\slie_3(\mathbf Z_3)$ has discriminant ideal $(3)$ and the map to $\pgl_3(\mathbf Z_3)$ has cokernel $\mathbf Z_3/3\mathbf Z_3$.

For example $A=E_{01}$ has monodromy boundary $(1,0)$, while $B=E_{10}$ is the unramified semistable dual class corresponding to that first unit direction. Their Tate pairing is $1$. The period lift for $A$ uses $\omega_1=\ell_p$ and all other $\omega_i=0$, so $Nh_A=e_1$ and its Hodge correction is zero. This computes the boundary value directly from an actual logarithmic period.
\end{example}

\section{Logarithmic descent and its cup-product coordinates}\label{sec:logdescent}

\subsection{Horizontalization and a contracting homotopy}

Let $A$ be a commutative $\Qp$-algebra, let $D_A$ be a finite free $A$-module, and let $N_D$ be a nilpotent $A$-linear endomorphism. Introduce an independent polynomial variable $\ell$ and set
\[
S=A[\ell],\qquad \Delta=-\partial_\ell+N_D.
\]
This subsection is an algebraic calculation. When $S$ is embedded in logarithmic periods, $\Delta$ is the sum of coefficient-ring monodromy and module monodromy.

\begin{proposition}\label{prop:horizontal}
The map
\begin{equation}\label{eq:horizontal}
\mathcal H_\ell:D_A\longrightarrow\ker\Delta,\qquad
x\longmapsto e^{\ell N_D}x
\end{equation}
is an isomorphism, inverse to evaluation at $\ell=0$. Define
\begin{equation}\label{eq:homotopy}
\mathcal K(F)=-e^{\ell N_D}\int_0^\ell e^{-sN_D}F(s)\,\dd s.
\end{equation}
Then $\Delta\mathcal K=1$ and
$\mathcal K\Delta=1-\mathcal H_\ell\ev_0$.
If $N_D^2=0$, coefficient extraction gives an exact sequence
\begin{equation}\label{eq:boundary-horizontal}
0\longrightarrow\ker N_D\longrightarrow\ker\Delta
\xrightarrow{[\ell]}\im N_D\longrightarrow0.
\end{equation}
\end{proposition}
\begin{proof}
All exponentials are finite polynomials. Differentiation gives
$\Delta(e^{\ell N_D}x)=0$. Conversely $\Delta F=0$ implies
$\partial_\ell(e^{-\ell N_D}F)=0$, whence $F=e^{\ell N_D}F(0)$. The same conjugation gives
$\Delta F=-e^{\ell N_D}\partial_\ell(e^{-\ell N_D}F)$, and the two homotopy identities follow from the fundamental theorem of polynomial integration. If $N_D^2=0$, every horizontal polynomial is $x+\ell N_Dx$. Its $\ell$-coefficient is $N_Dx$, with kernel $\ker N_D$ and image $\im N_D$.
\end{proof}

If $N_D\varphi_D=p\varphi_DN_D$ and $\varphi(\ell)=p\ell$, the last map in \eqref{eq:boundary-horizontal} is Frobenius compatible when $\im N_D$ is equipped with $p\varphi_D$. The factor $p$ follows directly from $N_D\varphi_D=p\varphi_DN_D$. This polynomial boundary and the arithmetic boundary of \Cref{thm:boundary} express the same lift-changing mechanism at different levels. The coefficient lattice that connects the two descriptions is constructed in \Cref{prop:integralRobba}.

\subsection{Affine Frobenius and group actions}

Suppose $A$ has commuting endomorphisms $\varphi_A,\gamma_A$, and the logarithmic variable transforms by
\begin{equation}\label{eq:affine}
\varphi(\ell)=p\ell+\eta_\varphi,\qquad
\gamma(\ell)=\ell+\eta_\gamma,
\end{equation}
where $\eta_\varphi,\eta_\gamma\in A$. Assume the untwisted semilinear operators $\Phi_0,G_0$ on $D_A$ commute, satisfy $N_D\Phi_0=p\Phi_0N_D$ and $N_DG_0=G_0N_D$, and extend these coefficient actions. Commutation on $S$ means
\begin{equation}\label{eq:etaidentity}
\eta_\varphi+\varphi_A(\eta_\gamma)
=p\eta_\gamma+\gamma_A(\eta_\varphi).
\end{equation}

\begin{proposition}\label{prop:affinedescent}
Transport of the actions through \eqref{eq:horizontal} gives
\begin{equation}\label{eq:descended-actions}
F=e^{\eta_\varphi N_D/p}\Phi_0,\qquad
G=e^{\eta_\gamma N_D}G_0.
\end{equation}
These two semilinear operators commute. If $\ell'=\ell+c$ with $c\in A$, then
\begin{equation}\label{eq:etachange}
\eta'_\varphi=\eta_\varphi+\varphi_A(c)-pc,\qquad
\eta'_\gamma=\eta_\gamma+\gamma_A(c)-c,
\end{equation}
and $U=e^{-cN_D}$ identifies the descended coordinates by
$F'=UFU^{-1}$ and $G'=UGU^{-1}$.
\end{proposition}
\begin{proof}
Since $\Phi_0N_D=p^{-1}N_D\Phi_0$,
\[
\varphi(\mathcal H_\ell x)
=e^{(p\ell+\eta_\varphi)N_D/p}\Phi_0x
=\mathcal H_\ell e^{\eta_\varphi N_D/p}\Phi_0x.
\]
The calculation for $\gamma$ uses its commutation with $N_D$ and gives the second formula. The original actions commute by \eqref{eq:etaidentity}, so their transports commute. The formulas for $\eta'$ follow by substitution. Finally $\mathcal H_{\ell'}=\mathcal H_\ell e^{cN_D}$, proving conjugacy by $U$; this computation includes the semilinear action on $c$.
\end{proof}

\subsection{A cup formula and the apparent correction}

Let $D_A^\dagger$ carry the dual data, and let
$b:D_A\times D_A^\dagger\to L_A$ be a perfect pairing into a rank-one module with zero monodromy. Assume $b(N_Dx,y)+b(x,N_D^\dagger y)=0$ and compatibility with Frobenius and $\gamma$. For commuting actions $F,G$, use the cochain convention
\begin{equation}\label{eq:herr-formal}
C^\bullet(D_A):\quad D_A\xrightarrow{(F-1,G-1)}D_A^2
\xrightarrow{(1-G,F-1)}D_A.
\end{equation}
For $x=(u,v)$ and $y=(u',v')$ of degree one set
\begin{equation}\label{eq:cupformula}
Z(x,y)=b(u,F^\dagger v')-b(v,G^\dagger u').
\end{equation}

\begin{theorem}\label{thm:cupcorrection}
Formula \eqref{eq:cupformula} descends to a pairing in degree-two cohomology. If $(N_D^\dagger)^2=0$, it is
\begin{align}\label{eq:ZplusC}
Z&=Z_0+\mathcal C_N,\\
Z_0&=b(u,\Phi_0^\dagger v')-b(v,G_0^\dagger u'),\notag\\
\mathcal C_N&=\frac{\eta_\varphi}{p}b(u,N_D^\dagger\Phi_0^\dagger v')
-\eta_\gamma b(v,N_D^\dagger G_0^\dagger u').\notag
\end{align}
Changing $\ell$ to $\ell+c$ and transporting both arguments by $U,U^\dagger$ leaves the total $Z$ unchanged, even before passage to cohomology. The summand $\mathcal C_N$ is consequently a coordinate term, rather than an independently defined invariant.
\end{theorem}
\begin{proof}
The cup formula is the usual ordered cup product for a two-generator commuting-action complex. One can check descent directly. If $x=((F-1)a,(G-1)a)$ and $y$ is a cocycle, pairing equivariance gives
\[
Z(x,y)=(F_L-1)b(a,v')-(G_L-1)b(a,u'),
\]
because $(1-G^\dagger)u'+(F^\dagger-1)v'=0$. This is a boundary in $C^2(L_A)$. If $y$ is a coboundary and $x$ is a cocycle, the same calculation gives
\[
Z(x,y)=-(F_L-1)b(v,G^\dagger a')
+(G_L-1)b(u,F^\dagger a').
\]
Hence both choices of representatives are harmless. Expanding the exponentials in \eqref{eq:descended-actions} to first order in $N_D^\dagger$ gives \eqref{eq:ZplusC}.

Monodromy adjointness implies
\[
b(e^{sN_D}x,e^{sN_D^\dagger}y)=b(x,y),
\]
as follows by differentiating this polynomial in $s$. Thus $U$ and $U^\dagger$ preserve the pairing. Conjugacy in \Cref{prop:affinedescent} now gives
\[
b(Uu,F'^\dagger U^\dagger v')=b(Uu,U^\dagger F^\dagger v')
=b(u,F^\dagger v'),
\]
and the second term is identical as well. This proves the claimed invariance.
\end{proof}

For the mixed motives of \Cref{thm:mixedrealization}, $N^2=0$, so \eqref{eq:ZplusC} applies to their logarithmic period data as well. Its monodromy term only involves the outer weight blocks: $N^\dagger$ maps the dual lattice block into the dual torus, and pairing with that torus reads the original lattice quotient. The abelian graded pairing contributes to $Z_0$ through the Poincar\'e realization. Consequently a change of logarithmic coordinate changes the displayed correction without changing the full pairing, just as the point-branch change in \Cref{cor:mixedbranch} leaves the exponential quotient invariant. This algebraic transport does not replace the additional coefficient input required in \Cref{thm:relativeHerrlift}.

For the basis family, write $Nf=e_0$ and choose dual vectors $b_0,\ldots,b_n,a$ with $b(e_i,b_j)=\delta_{ij}e_\Tate$ and $b(f,a)=e_\Tate$. Then $N^\dagger b_0=-a$ and $N^\dagger b_i=0$ for $i>0$. If $u_f,v_f$ denote the $f$-coordinates and $c_0,d_0$ the $b_0$-coordinates of $u',v'$, respectively, \eqref{eq:ZplusC} reads
\begin{equation}\label{eq:scalarcorrection}
\mathcal C_N=
\left(-\frac{\eta_\varphi}{p}u_f\varphi_A(d_0)
+\eta_\gamma v_f\gamma_A(c_0)\right)e_\Tate.
\end{equation}
The full coefficient $Z$ belongs to the actual Tate complex for the representatives constructed in \Cref{thm:fulldualrepresentatives}. Its normalized residue is evaluated in \Cref{thm:fullresiduereciprocity}. The summands in \eqref{eq:ZplusC} depend on logarithmic coordinates and can require localization separately; the residue theorem applies to their full cup product.

\subsection{Cyclotomic coordinates and normalization}

Let $X$ be the cyclotomic Robba variable, with
$\varphi(X)=(1+X)^p-1$ and $\gamma(X)=(1+X)^{\chi(\gamma)}-1$. A Robba ring consists of Laurent series converging on some annulus $r<|X|<1$, with $r<1$ allowed to depend on the series. In Berger's corrected normalization,
\[
N(\log X)=-\frac{p}{p-1}.
\]
This is the normalization specified in \cite[Appendix B]{Ber08}, correcting the earlier convention in \cite{Ber02}. Thus a variable with $N\ell=-1$ is $\ell=(p-1)\log X/p$, and its affine terms are
\begin{equation}\label{eq:cyclotomicetas}
\eta_\varphi=\frac{p-1}{p}\log\frac{(1+X)^p-1}{X^p},\qquad
\eta_\gamma=\frac{p-1}{p}\log\frac{(1+X)^{\chi(\gamma)}-1}{X}.
\end{equation}
Here $\chi(\gamma)\in1+p\Zp$. The first logarithm converges uniformly on every closed subannulus of
$0<v_p(X)<1/(p-1)$: each nonconstant term in its argument minus one has positive valuation there. The second converges on such subannuli because its argument is $1$ modulo the maximal ideal. These statements concern convergence of the coefficient functions; an action on a fixed radius requires the usual change of annulus under Frobenius.

For \Cref{ex:three}, take $\chi(\gamma)=4$. The functions are
\[
\eta_\varphi=\frac23\log\left(1+\frac3X+\frac3{X^2}\right),\qquad
\eta_\gamma=\frac23\log(4+6X+4X^2+X^3).
\]
They converge uniformly on closed subannuli of $0<v_3(X)<1/2$. These expressions can be substituted into \eqref{eq:scalarcorrection} using the actual module and representatives of \Cref{sec:explicitrobba,sec:fullresidue}.

On the classical cyclotomic Robba ring $\mathcal R_K$, the finite-free Frobenius extension has degree $p$ and defines
\[
\psi_{\mathrm{cyc}}(f)=\varphi^{-1}\left(\frac1p
\Tr_{\mathcal R_K/\varphi(\mathcal R_K)}(f)\right),\qquad
\psi_{\mathrm{cyc}}\varphi=1.
\]
The last identity is the trace identity for an element of the base subring. This familiar operator \cite[Section 4.4.2]{FO} explains why a decompleted finite-free Frobenius extension supplies more information than the automorphism of the perfect ambient ring. The coefficient traces used throughout the multi-Kummer tower are specified in \Cref{subsec:towerperiods}.

\section{Explicit Robba modules and integral coefficient lattices}\label{sec:explicitrobba}

We now pass from the ambient periods of \Cref{sec:actualperiods} to the cyclotomic coefficient ring. The construction has two stages. Interpolation at cyclotomic points gives a basis satisfying Berger's reconstruction conditions. A change of representative in the rank-one Kummer classes then gives integral matrices. The comparison theorem is Berger's; the formulas below make its conditions explicit for our toric module.

\subsection{Coefficient rings and the rank-one inputs}

Write $\mathcal R_K$ for the Robba ring defined in \Cref{sec:logdescent}, and write $\mathcal R_K^+$ for the ring of power series converging on $|X|<1$. These are decompleted rings. In particular, $\mathcal R_K^+$ and the perfect ambient ring $\Rplus$ have different meanings. Frobenius acts on $K$ by $\sigma$. Put
\[
\begin{aligned}
 \mathcal E_K^\dagger&=\{f\in\mathcal R_K:\text{the Laurent coefficients of }f\text{ are bounded}\},\\
 \mathcal A_K^\dagger&=\left\{\sum_{j\in\mathbf Z}a_jX^j\in\mathcal E_K^\dagger:a_j\in\OK\right\},
\end{aligned}
\]
and let $\mathcal A_K$ be the $p$-adic completion of $\mathcal A_K^\dagger$. Thus the negative coefficients of an element of $\mathcal A_K$ tend to zero $p$-adically, but need not satisfy an overconvergence estimate. Set $t=\log(1+X)$ and $c=(p-1)/p$.

Let $\Delta_K\simeq\mu_{p-1}$ be the torsion subgroup of $\Gamma_K$. An integral \emph{\'etale} $(\varphi,\Gamma_K)$-module over $\mathcal A_K$ is a finite free module with continuous commuting semilinear actions for which the linearized Frobenius is an isomorphism. Fontaine's coefficient equivalence attaches such a module $D_K(T)$ to every Galois lattice $T$. Its overconvergent rational form is denoted $D_K^\dagger(T\otimes\Qp)$; scalar extension to $\mathcal R_K$ is denoted $D_{\mathrm{rig}}^\dagger(T\otimes\Qp)$. The equivalence and overconvergence are those of \cite[Section I.3 and Proposition III.3.2]{CC99}. On the $\Delta_K$-fixed part, the complex \eqref{eq:herr-formal} with $F=\varphi$ and $G=\gamma$ computes full Galois cohomology and its cup products \cite{Her01,Liu08}. We use these comparisons throughout \Cref{sec:explicitrobba,sec:fullresidue}.

Fix $\gamma\in\Gamma_K$ with $\chi(\gamma)$ generating $1+p\Zp$. For the Tate coefficient module $\mathcal R_K(1)=\mathcal R_K v$, the actions are
\[
 \varphi(v)=v,\qquad \delta(v)=\chi(\delta)v,\qquad \gamma(v)=\chi(\gamma)v.
\]
Denote averaging for this action by
\begin{equation}\label{eq:tate-average}
 \operatorname{Av}_1(f)=\frac1{p-1}\sum_{\delta\in\Delta_K}\chi(\delta)\delta(f),
 \qquad h(X)=X^{-1}+\tfrac12,\qquad h_\Delta=\operatorname{Av}_1(h).
\end{equation}
Averaging without the factor $\chi(\delta)$ will be denoted $\operatorname{Av}_0$. Both operators preserve integral coefficients. Throughout this section, a Tate-valued Herr cocycle is a pair $(u,v)$ satisfying
\begin{equation}\label{eq:tate-cocycle}
 (\varphi-1)v=(\chi(\gamma)\gamma-1)u.
\end{equation}
The first entry is the Frobenius entry. This order agrees with \eqref{eq:herr-formal}.

We use two established coefficient-ring results. First, $\gamma-1$ is invertible on the $\psi=0$ part of an integral \emph{\'etale} coefficient module; its inverse preserves overconvergence. This is \cite[Proposition I.5.1 and Proposition III.3.2(iii)]{CC99}. Applying this statement to the Tate module gives an inverse for $\chi(\gamma)\gamma-1$ on its scalar coefficients with $\psi=0$. Second, the rank-one residue calculation of Colmez and Benois identifies $H^1(\mathcal R_{\Qp}(1))$ by the two residues of a cocycle; see \cite[Sections 2.3--2.7]{Col08} and \cite[Corollary 1.5.5]{Ben11}. Our use of these results includes their stated coefficient modules and the order of the two entries.

\subsection{An interpolation formula for the unit filtration}

For $j\geq0$ set
\begin{equation}\label{eq:beta-interpolation}
 \beta_j(X)=\frac{t}{\varphi^j(X)}-p^{-j}+\frac t2.
\end{equation}
The quotient is analytic on the open unit disc: every zero of its denominator is a simple zero of $t$. If $\pi_m=\zeta_{p^m}-1$ for a primitive root $\zeta_{p^m}$, then
\begin{equation}\label{eq:beta-values}
 \beta_j(\pi_m)=\begin{cases}-p^{-j}&j<m,\\0&j\geq m,\end{cases}
 \qquad \beta_j(0)=\beta_j'(0)=0,
 \qquad \varphi(\beta_j)=p\beta_{j+1}.
\end{equation}

\begin{lemma}\label{lem:unitinterpolation}
For $L_0\in K$, put $C_{L_0}=L_0-p^{-1}\sigma(L_0)$. The series
\begin{equation}\label{eq:unit-F}
 F_{L_0}=L_0+\sum_{j\geq0}\sigma^j(C_{L_0})\beta_j
\end{equation}
converges in $\mathcal R_K^+$. It satisfies
\begin{equation}\label{eq:unit-F-identities}
 F_{L_0}(\pi_m)=p^{-m}\sigma^m(L_0),\qquad
 (1-p^{-1}\varphi)F_{L_0}=C_{L_0}t h,\qquad F_{L_0}'(0)=0.
\end{equation}
If $L_0\in p\OK$, the cocycle
\begin{equation}\label{eq:unit-cocycle}
 u_{L_0}=\operatorname{Av}_1(-C_{L_0}h),\qquad
 v_{L_0}=\operatorname{Av}_1\left(\frac{(\gamma-1)F_{L_0}}t\right)
\end{equation}
belongs to $(\mathcal A_K^\dagger)^2$ and represents $\kappa(\exp L_0)$.
\end{lemma}
\begin{proof}
Fix a closed disc $|X|\leq\rho<1$. For large $j$, one has $\varphi^j(X)=\exp(p^jt)-1$ on this disc. Expanding $z/(e^z-1)=1-z/2+O(z^2)$ gives $\|\beta_j\|_\rho=O(p^{-j})$. The coefficients $\sigma^j(C_{L_0})$ have constant norm, so \eqref{eq:unit-F} converges there. The values in \eqref{eq:beta-values} telescope to the first identity of \eqref{eq:unit-F-identities}. Shifting the series by Frobenius gives the second identity. The derivative assertion follows from \eqref{eq:beta-values}.

The numerator $(\gamma-1)F_{L_0}$ vanishes at every nonzero zero of $t$, since $\gamma$ permutes the primitive roots at each level. It vanishes to order at least two at zero. Consequently its quotient by $t$ is analytic and has constant term zero. Before averaging, write $u=-C_{L_0}h$ and $v=(\gamma-1)F_{L_0}/t$. The Frobenius identity gives \eqref{eq:tate-cocycle}. If $L_0\in p\OK$, then $C_{L_0}\in\OK$, and
\[
 r=(\chi(\gamma)\gamma-1)u\in X\OK[[X]].
\]
Indeed, the pole and the constant term cancel in $\chi(\gamma)\gamma(h)-h$. The series $-\sum_{m\geq0}\varphi^m(r)$ converges on every closed subdisc, has integral coefficients, and solves $(\varphi-1)v=r$. The difference from the previously constructed $v$ is a Frobenius-fixed analytic function with zero constant term, hence is zero. This proves integrality, including after averaging.

For completeness, the identification of the class can be read directly from the filtered reconstruction. In the crystalline extension with filtration $\Fil^0=K(f+L_0e)$, put $v=te$ and $w=f+F_{L_0}e$. At the $m$th de Rham embedding, the coefficient of $e$ in $w$ becomes $p^m\sigma^{-m}(F_{L_0}(\pi_m))=L_0$. Thus $v,w$ give Berger's reconstructed module. The crystalline Kummer calculation in \Cref{prop:toricmodule} identifies its extension with $\kappa(\exp L_0)$. Averaging $w$ changes its lift by an element of $\mathcal R_Kv$ and produces exactly \eqref{eq:unit-cocycle}. We use Berger's comparison \cite[Definition II.2.2 and Theorem V.2.3(2)]{Ber08} here, rather than infer a Galois comparison from interpolation alone.
\end{proof}

\subsection{The valuation interpolation and the reconstructed basis}

We give a convergent interpolation for the logarithmic part as well. Put
\[
 Q_m(X)=\Phi_{p^m}(1+X),\quad U_m=Q_m/p,\quad
 E_m=(1-U_m^2)^2\prod_{\nu>m}\bigl(1-(1-U_\nu)^2\bigr)^2.
\]
Here $\Phi_{p^m}$ is the cyclotomic polynomial. These selectors are the functions $p^mG_{m,2}$ of \cite[Definition 3.1]{CGY25}. We recall the convergence argument in the form needed here. On a fixed closed subdisc, $\|U_m-1\|_\rho=O(p^{-m})$. The infinite product therefore converges, and $\|E_m\|_\rho\leq C_\rho p^{-2m}$ for large $m$. At the primitive cyclotomic points, $E_m$ is one at level $m$ and zero at every other level. In fact these congruences hold modulo $Q_j^2$.

Let $P_m\in\Qp[X]$ be the polynomial of degree less than $e_m=p^{m-1}(p-1)$ whose value at $\pi_m$ is $\log_p(\pi_m)$. The branch $\log_p(p)=0$ makes these values Galois equivariant. The polynomial is therefore defined over $\Qp$.

\begin{lemma}\label{lem:valuationinterpolation}
The series
\begin{equation}\label{eq:valuation-F}
 F_0=-c\sum_{m\geq1}E_mP_m
\end{equation}
converges in $\mathcal R_{\Qp}^+$ and satisfies
\[
 F_0(\pi_m)=-c\log_p(\pi_m),\qquad F_0(0)=F_0'(0)=0.
\]
\end{lemma}
\begin{proof}
The ratio $\pi_{m-1}/\pi_m^p$ is a principal unit for $m\geq2$. Its logarithm has norm at most one: each nonleading term of the ratio minus one has valuation at least $1-(p-1)/e_m>1/(p-1)$. For $m=1$, the element $\pi_1^{p-1}/p$ is a unit in the tamely ramified field $\Qp(\pi_1)$. After removing its prime-to-$p$ torsion part, this unit differs from one by an element of valuation at least $1/(p-1)$. Its logarithm therefore has norm at most one. The identity $(p-1)\log_p\pi_1=\log_p(\pi_1^{p-1}/p)$ proves $|\log_p\pi_1|\leq1$. Induction in
\[
 p\log_p\pi_m=\log_p\pi_{m-1}-\log_p(\pi_{m-1}/\pi_m^p)
\]
gives $|\log_p\pi_m|\leq p^{m-1}$. The basis $1,\pi_m,\ldots,\pi_m^{e_m-1}$ of $\Qp(\pi_m)$ has terms with distinct valuations modulo $\mathbf Z$. It follows that every coefficient of $P_m$ has norm at most $p^m$. Consequently $\|E_mP_m\|_\rho=O(p^{-m})$. This proves convergence. The selector congruences give the values. Each $E_m$ is divisible by $X^2$, proving the two assertions at zero.
\end{proof}

Use $\ell=c\log X$ and $\eta_\varphi,\eta_\gamma$ from \eqref{eq:cyclotomicetas}. In the logarithmic tensor product with $D=\Dst(V_Q)$ define
\begin{equation}\label{eq:actual-Robba-basis}
 v_i=te_i,\qquad
 w_j=f_j+\sum_i\bigl(R_{ij}(\ell+F_0)+F_{L_{ij}}\bigr)e_i.
\end{equation}
The symbols $v_i,w_j$ now denote a Robba basis; the comparison identifies their toric and quotient frames with those of \eqref{eq:galois-toric}.

\begin{theorem}\label{thm:Robbabasis}
The vectors in \eqref{eq:actual-Robba-basis} form an $\mathcal R_K$-basis of $D_{\mathrm{rig}}^\dagger(V_Q)$. In that basis,
\begin{align}\label{eq:Robba-matrices}
 \varphi(v_i)&=v_i,&\gamma(v_i)&=\chi(\gamma)v_i,\notag\\
 \varphi(w_j)&=w_j+\sum_i u_{ij}v_i,&
 \gamma(w_j)&=w_j+\sum_i v_{ij}v_i,
\end{align}
where, before averaging the quotient lifts,
\begin{align}\label{eq:valuation-raw}
 u_{ij}&=R_{ij}u_p^{\mathrm{an}}-C_{L_{ij}}h,&
 v_{ij}&=R_{ij}v_p^{\mathrm{an}}+t^{-1}(\gamma-1)F_{L_{ij}},\notag\\
 u_p^{\mathrm{an}}&=t^{-1}\bigl(\eta_\varphi/p+(p^{-1}\varphi-1)F_0\bigr),&
 v_p^{\mathrm{an}}&=t^{-1}\bigl(\eta_\gamma+(\gamma-1)F_0\bigr).
\end{align}
All these coefficients belong to $\mathcal R_K$. Replacing $w_j$ by its $\Delta_K$-average replaces every coefficient by $\operatorname{Av}_1$ and makes $w_j$ fixed by $\Delta_K$.
\end{theorem}
\begin{proof}
The vectors are horizontal for the sum of ring and module monodromy, because $N\ell=-1$ and $Nf_j=\sum_iR_{ij}e_i$. They form a basis after inverting $t$. To check the integral powers of $t$ required by the filtration, use the embedding
\[
 \iota_m(X)=\zeta_{p^m}\exp(t/p^m)-1,\quad
 \iota_m|_K=\sigma^{-m},\quad
 \iota_m(\log X)=\log\bigl(\zeta_{p^m}\exp(t/p^m)-1\bigr).
\]
In Berger's transported filtration, $e_i$ is multiplied by $p^m$, while $f_j$ is unchanged. \Cref{lem:unitinterpolation,lem:valuationinterpolation} show that the constant term of the coefficient of $e_i$ in the image of $w_j$ is $L_{ij}$. The image of $v_i$ is $te_i$. Hence the image lattice is precisely $\Fil^0(K_m((t))\otimes D)$ at every sufficiently deep embedding.

Here is the elementary divisibility point that makes this a basis statement. On an annulus the zeros of $t$ are the simple zeros $\pi_m$ and their conjugates. The product $t=X\prod_{m\geq1}(Q_m/p)$ converges on the open unit disc. A function in $\mathcal R_K[1/t]$ has no poles exactly when its poles at all sufficiently deep such points cancel; a function in $\mathcal R_K$ is divisible by $t$ exactly when it vanishes at those points. The finitely many earlier points can be excluded by shrinking the annulus. Apply this observation to the coefficients of a horizontal element in the basis \eqref{eq:actual-Robba-basis}. The reconstructed filtration forces its quotient coefficients to be regular and its remaining toric coefficients to be divisible by $t$. Thus the coefficients in the displayed basis belong to $\mathcal R_K$.

Berger's reconstruction and comparison \cite[Definition II.2.2, Theorem II.2.6, and Theorem V.2.3(2)]{Ber08} now identify this module with $D_{\mathrm{rig}}^\dagger(V_Q)$. Applying the actions gives \eqref{eq:Robba-matrices}--\eqref{eq:valuation-raw}. Alternatively, the apparent quotients by $t$ are regular by the same cyclotomic-point calculation. Averaging a quotient lift changes it by a toric vector and preserves a basis. Since the actions commute with $\Delta_K$, its new coefficients are their Tate averages.
\end{proof}

\subsection{Integral matrices and the trace operator}

The analytic basis above records the filtration directly. To normalize its valuation component integrally, we give a direct weight-one argument. Extend the actions to $\mathcal R_{\Qp}[\log X]$ by the affine formulas of \Cref{prop:affinedescent}; for $\delta\in\Delta_K$, use $\delta(\log X)=\log X+\log(\delta(X)/X)$ with $\log\chi(\delta)=0$. Define the logarithmic averaging operator on $\log X$ by
\begin{equation}\label{eq:logarithmic-average}
 J_\gamma(\log X)=\log X+
 \sum_{k\geq1}\frac{(-1)^k}{k+1}(\gamma-1)^k\log X,
 \qquad A_p=c\operatorname{Av}_0\bigl(J_\gamma(\log X)\bigr).
\end{equation}
The notation $J_\gamma=\log\gamma/(\gamma-1)$ denotes this power series. It does not require an inverse of $\gamma-1$ on constants. This logarithmic averaging construction occurs in Benois's treatment, as written in \cite[Section 3.2, equation (16)]{CGY25}. We prove the properties used here directly in weight one.

\begin{lemma}[The normalized valuation period]\label{lem:valuationaverage}
The series \eqref{eq:logarithmic-average} converges after its formal logarithmic term is removed. More precisely,
\[
 A_p=c\log X+F_\gamma,\qquad F_\gamma\in\mathcal R_{\Qp}^+,
 \qquad N A_p=-1,\qquad A_p(\pi_m)=0\quad(m\geq1).
\]
For the normalized trace extended to the logarithmic coefficient ring,
\begin{equation}\label{eq:average-psi}
 \psi(A_p)=p^{-1}A_p,\qquad
 (\gamma-1)A_p=c\log\chi(\gamma)\,t h_\Delta.
\end{equation}
Consequently
\begin{equation}\label{eq:valuation-average-cochain}
 u_p'=t^{-1}(p^{-1}\varphi-1)A_p\in\mathcal R_{\Qp},\qquad
 v_p'=t^{-1}(\gamma-1)A_p=c\log\chi(\gamma)h_\Delta
\end{equation}
is a Tate Herr cocycle representing $\kappa(p)$, and $\psi(u_p')=0$.
The change from the averaged selector period to $A_p$ is the explicit gauge
\begin{equation}\label{eq:valuation-gauge}
 d_\gamma=\frac{A_p-\operatorname{Av}_0(c\log X+F_0)}t,
 \qquad d_\gamma+\frac{c\log\chi(\gamma)}{2X}\in\mathcal R_{\Qp}^+.
\end{equation}
\end{lemma}
\begin{proof}
Fix $\rho<1$. The function $H_\gamma=\gamma\log X-\log X=\log(\gamma(X)/X)$ is analytic on the disc. On analytic functions there, the substitution $\gamma$ satisfies
\[
 \|(\gamma-1)f\|_\rho\leq q_\rho\|f\|_\rho,
 \qquad q_\rho=\max(p^{-1},\rho^{p-1})<1.
\]
Indeed, $\|\gamma(X)-X\|_\rho\leq q_\rho\rho$, and the estimate follows first for monomials, then for convergent series. The $k$th summand of \eqref{eq:logarithmic-average}, for $k\geq1$, is therefore bounded by $|k+1|_p^{-1}q_\rho^{k-1}\|H_\gamma\|_\rho$. This proves convergence. The functions $\delta\log X-\log X$ are analytic as well, so averaging gives $F_\gamma\in\mathcal R_{\Qp}^+$. Monodromy commutes with these operators, acts on $\log X$ as $-1/c$, and kills analytic coefficients. Hence $N A_p=-1$.

At the primitive $p^m$th roots, $\gamma$ acts through a cyclic group of order $d_m=p^{m-1}$. Evaluation of $J_\gamma$ on this finite permutation representation is the averaging projector $d_m^{-1}\sum_{j=0}^{d_m-1}\gamma^j$. To verify this assertion, diagonalize the finite-order operator: the series $\log z/(z-1)$ has value zero at a nontrivial $p$-power root of unity and value one at $z=1$. Combining this projector with $\operatorname{Av}_0$ gives
\[
 A_p(\pi_m)=\frac{c}{p^{m-1}(p-1)}
 \sum_{a\in(\mathbf Z/p^m\mathbf Z)^\times}\log_p(\zeta_{p^m}^{a}-1)
 =\frac{c\log_p p}{p^{m-1}(p-1)}=0.
\]
We used $\Phi_{p^m}(1)=p$ and the chosen logarithm branch. Thus $F_\gamma$ has exactly the interpolation values of $F_0$ in \Cref{lem:valuationinterpolation}.

The trace extends to $\mathcal R_{\Qp}[\log X]$ using the same degree-$p$ Frobenius extension. Indeed, $p\log X=\varphi(\log X)-\eta_\varphi/c$, so this polynomial extension is the base change of the original Frobenius extension and retains degree $p$. On an outer annulus sufficiently close to the outer boundary, the conjugates of $\log X$ are $\log(\zeta(1+X)-1)$ for $\zeta^p=1$. Their sum is $\log\varphi(X)$, since their arguments have product $\varphi(X)$ and the logarithms of roots of unity vanish. Hence
\[
 \varphi\psi(\log X)=p^{-1}\varphi(\log X),
 \qquad \psi(\log X)=p^{-1}\log X.
\]
The trace commutes with $\gamma$ and $\Delta_K$ and is continuous on annuli where it is defined. Applying it to the convergent series proves the first identity of \eqref{eq:average-psi}.

The power-series identity $(\gamma-1)J_\gamma=\log\gamma$, together with differentiation of the cyclotomic substitution, gives
\[
 (\log\gamma)(\log X)=\log\chi(\gamma)\,t\frac{1+X}{X}.
\]
Since $\operatorname{Av}_0(t)=0$ and $\operatorname{Av}_0(tf)=t\operatorname{Av}_1(f)$, this proves the second identity of \eqref{eq:average-psi}.

The difference $A_p-\operatorname{Av}_0(c\log X+F_0)$ is analytic and vanishes at all nonzero cyclotomic zeros of $t$. Its quotient $d_\gamma$ therefore extends analytically across those points and has at most a simple pole at zero. The constant term of $F_\gamma$ is $-c\log\chi(\gamma)/2$: the $k=1$ summand of \eqref{eq:logarithmic-average} contributes this term, and every later summand has constant term zero. Averaging the logarithmic terms for $\Delta_K$ introduces no constant because the Teichm\"uller roots have logarithm zero. The averaged selector period has analytic constant term zero. Since $t=X+O(X^2)$, these facts prove \eqref{eq:valuation-gauge}.

Writing $(\bar u_p^{\mathrm{an}},\bar v_p^{\mathrm{an}})$ for the averaged analytic pair, the change of lift is explicitly
\[
 u_p'=\bar u_p^{\mathrm{an}}+(\varphi-1)d_\gamma,\qquad
 v_p'=\bar v_p^{\mathrm{an}}+(\chi(\gamma)\gamma-1)d_\gamma.
\]
This proves both the coefficient-ring membership of \eqref{eq:valuation-average-cochain} and its identification with $\kappa(p)$ by \Cref{thm:Robbabasis}. Finally, $\psi(t^{-1}b)=p t^{-1}\psi(b)$ and the first identity of \eqref{eq:average-psi} give
\[
 \psi(u_p')=p t^{-1}\bigl(p^{-1}A_p-\psi(A_p)\bigr)=0.
\]
The cocycle equation follows from the commuting actions, completing the proof.
\end{proof}

For integral calculations define
\begin{equation}\label{eq:valuation-integral}
 v_p=c\log\chi(\gamma)\,h_\Delta,\qquad
 u_p=(\chi(\gamma)\gamma-1)^{-1}(\varphi-1)v_p,
 \qquad \psi(u_p)=0.
\end{equation}
The inverse is taken on the $\psi=0$ space, not on all Robba functions.

\begin{proposition}\label{prop:integralRobba}
The coefficients in \eqref{eq:valuation-integral} belong to $\mathcal A_{\Qp}^\dagger$ and represent the same class as the averaged pair in \eqref{eq:valuation-raw}. Consequently the matrices
\begin{equation}\label{eq:integral-Kummer-matrices}
 u_{ij}=R_{ij}u_p+u_{L_{ij}},\qquad
 v_{ij}=R_{ij}v_p+v_{L_{ij}}
\end{equation}
define an integral \emph{\'etale} coefficient module for $T_Q$, with quotient lifts fixed by $\Delta_K$. Its scalar extension to $\mathcal R_K$ is the module in \Cref{thm:Robbabasis}.
\end{proposition}
\begin{proof}
The trace definition of $\psi$ gives $\psi(h)=h$, hence $\psi(v_p)=v_p$. Since $c\log\chi(\gamma)\in\Zp$, $v_p$ has integral coefficients. The element $(\varphi-1)v_p$ is integral, overconvergent, and killed by $\psi$. The integral inverse theorem cited above therefore gives $u_p\in\mathcal A_{\Qp}^\dagger$ and proves \eqref{eq:tate-cocycle}.

The pair $(u_p',v_p')$ of \Cref{lem:valuationaverage} has the same second entry and satisfies the same cocycle equation as the integral pair. Hence $(\chi(\gamma)\gamma-1)(u_p'-u_p)=0$. If $w\in\mathcal R_{\Qp}$ is killed by this operator, then $tw$ is $\gamma$-invariant and is constant. Indeed, with $\partial=(1+X)\dd/\dd X$, differentiation of the locally analytic cyclotomic action gives $t\partial(tw)=0$. The constant must be zero because a nonzero multiple of $t^{-1}$ has poles at arbitrarily deep cyclotomic points and does not belong to the Robba ring. Thus $w=0$, and $u_p'=u_p$. The lemma proves that the integral pair represents $\kappa(p)$.

For later normalization, write
\[
 \operatorname{res}(f)=\operatorname{Res}_{X=0}\left(f\frac{\dd X}{1+X}\right).
\]
The trace formula gives $\operatorname{res}\psi(f)=\operatorname{res}(f)$ over $\Qp$, while $\operatorname{res}(h_\Delta)=1$. Consequently \Cref{lem:valuationaverage} gives
\begin{equation}\label{eq:valuation-residues}
 \operatorname{res}(u_p)=0,\qquad
 \operatorname{res}(v_p)=c\log\chi(\gamma).
\end{equation}
Residues of a Tate coboundary vanish over $\Qp$ by the change-of-variable identities in \eqref{eq:residue-transform}. Thus the averaged analytic pair has the same residues. These coordinates agree with the rank-one residue description in \cite[Corollary 1.5.5]{Ben11}. Base change gives the assertion over $K$.

The unit entries are integral by \Cref{lem:unitinterpolation}. The Frobenius matrices in \eqref{eq:integral-Kummer-matrices} and their inverses have integral entries and determinant one. For clarity, the action of the whole procyclic group is determined by these entries. If $G_1=\chi(\gamma)\gamma$ on scalar Tate coefficients and $z\in\Zp$, then the scalar off-diagonal entry of $\gamma^z$ is
\[
 v_{ij}(z)=\frac{G_1^z-1}{G_1-1}v_{ij}
 =\sum_{k\geq0}\binom z{k+1}(G_1-1)^k v_{ij}.
\]
This series is the continuous Iwasawa-algebra action on the integral Tate module with its weak coefficient topology. The known continuous action on that module gives convergence; for integer $z\geq0$ it is the finite sum $\sum_{k=0}^{z-1}G_1^kv_{ij}$. The cocycle relation gives commutation with Frobenius, and the Tate averages give commutation with $\Delta_K$. Thus the matrices define an integral \emph{\'etale} coefficient module over $\mathcal A_K$.

The integral coefficient equivalence \cite[Section I.3]{CC99} attaches a Galois lattice to this module. Its rational framed extension has Kummer entries $\kappa(q_{ij})$. Integral framed extensions inject into rational ones here: the kernel is controlled by $H^0(K,(\Qp/\Zp)(1))$, which is zero because $K$ has no nontrivial $p$-power roots of unity. The integral representation is therefore $T_Q$, proving the claimed lattice identification.
\end{proof}

Let $D_Q$ denote this integral coefficient module or its scalar extension, according to the coefficient ring. Its left inverse to Frobenius is completely explicit. For coefficients $a_i,b_j$,
\begin{equation}\label{eq:psi-toric}
 \psi_{D_Q}\left(\sum_i a_iv_i+\sum_jb_jw_j\right)
 =\sum_i\left(\psi(a_i)-\sum_j\psi(b_ju_{ij})\right)v_i
   +\sum_j\psi(b_j)w_j.
\end{equation}
Indeed, rewrite the argument in the Frobenius basis $\varphi(v_i),\varphi(w_j)$ and apply the scalar normalized trace. This proves \eqref{eq:psi-toric}, including integrality and $\psi_{D_Q}\varphi=1$.

Choose dual basis vectors $b_i,a_j$ with $b(v_i,b_k)=\delta_{ik}v_{\Tate}$ and $b(w_j,a_k)=\delta_{jk}v_{\Tate}$. The integral dual matrices are
\begin{equation}\label{eq:Robba-dual-matrices}
 \varphi(b_i)=b_i-\sum_j u_{ij}a_j,\quad
 \gamma(b_i)=b_i-\sum_j v_{ij}a_j,\quad
 \varphi(a_j)=a_j,\quad\gamma(a_j)=\chi(\gamma)a_j.
\end{equation}
These follow from pairing equivariance. They will be used to construct every dual class in the basis family, including the classes excluded from the semistable extension complex by \Cref{cor:dimensionobstruction}.

\section{Full dual cochains and a direct residue reciprocity formula}\label{sec:fullresidue}

The integral matrices of \Cref{prop:integralRobba} describe both coefficient modules. To describe their degree-one classes, we need one more operation: solving a Tate-valued degree-two equation when its invariant vanishes. We give a solution inside the Robba ring. It supplies the missing scalar coordinate of a full dual cocycle and makes the trace-zero obstruction effective.

\subsection{The Tate residue and a solution of the degree-two equation}

For $f\in\mathcal R_K$ define
\[
 \operatorname{res}(f)=\operatorname{Res}_{X=0}\left(f\frac{\dd X}{1+X}\right)\in K.
\]
For a Laurent series on an annulus this means its Laurent coefficient of degree $-1$, not an assumption that the series extends meromorphically to zero. The change-of-variable formulas are
\begin{equation}\label{eq:residue-transform}
 \operatorname{res}(\varphi f)=\sigma\operatorname{res}(f),\qquad
 \operatorname{res}(\chi(\gamma)\gamma f)=\operatorname{res}(f),\qquad
 \operatorname{res}(\operatorname{Av}_1 f)=\operatorname{res}(f).
\end{equation}
They follow from Frobenius substitution and the differential $\dd X/(1+X)$. The residue description of $H^2(\mathcal R_K(1))$ is the standard Tate-coefficient calculation of Herr and Liu \cite{Her01,Liu08}. In particular, the trace of the residue kills both differentials.

Put $\mathcal R_K^{\Delta_K,1}=\operatorname{Av}_1\mathcal R_K$; this is the scalar coefficient space for the $\Delta_K$-fixed Tate module.

\begin{theorem}\label{thm:explicitprimitive}
There is a specified $\Qp$-linear operator
\[
 \mathcal T_\gamma:\{f\in\mathcal R_K^{\Delta_K,1}:\Tr_{K/\Qp}\operatorname{res}(f)=0\}
 \longrightarrow (\mathcal R_K^{\Delta_K,1})^2
\]
such that, if $\mathcal T_\gamma(f)=(a,b)$, then
\begin{equation}\label{eq:primitive-equation}
 (1-\chi(\gamma)\gamma)a+(\varphi-1)b=f.
\end{equation}
The operator is given by \eqref{eq:primitive-construction}--\eqref{eq:negative-resolvent} below. Its only inverse from coefficient-ring theory is $(1-\gamma)^{-1}$ on the bounded overconvergent $\psi=0$ space. All other operations are coefficientwise integration, finite linear algebra over $K$, and convergent series.
\end{theorem}
\begin{proof}
We give the construction and its convergence. Fix the trace-one element $\vartheta\in\OK$ from \Cref{lem:frobenius}. If $r=\operatorname{res}(f)$, then $\Tr r=0$. An explicit additive Hilbert 90 solution is
\[
 c_r=-\sum_{j=0}^{n-1}\sigma^j(\vartheta)\sum_{k=0}^{j-1}\sigma^k(r),
 \qquad (\sigma-1)c_r=r.
\]
The identity follows by shifting both sums, using $\Tr\vartheta=1$ and $\Tr r=0$; the inner sum is zero for $j=0$. Put
\begin{equation}\label{eq:primitive-construction}
 h_*=\operatorname{Av}_1\left(\frac{1+X}{X}\right),\qquad
 f_0=f-(\varphi-1)(c_rh_*),\qquad
 g=\int f_0\frac{\dd X}{1+X}.
\end{equation}
The last integral is the Laurent antiderivative with constant term zero. It exists because $\operatorname{res}(f_0)=0$. Dividing Laurent coefficients by nonzero integers preserves convergence on every slightly smaller annulus. Thus $g\in\mathcal R_K$ and $\partial g=f_0$, where $\partial=(1+X)\dd/\dd X$.

It remains to solve
\begin{equation}\label{eq:primitive-untwisted}
 (1-\gamma)a_0+(p^{-1}\varphi-1)b_0=g.
\end{equation}
This is useful because $\partial\gamma=\chi(\gamma)\gamma\partial$ and $\partial(p^{-1}\varphi)=\varphi\partial$. Decompose $g=g^++g^-$ into its nonnegative and negative Laurent parts. Write
\[
 A_0=g^+(0),\qquad A_1=(g^+)'(0),\qquad q=g^+-A_0-A_1t.
\]
The function $q$ is analytic on the disc and has a double zero at zero. Set
\begin{equation}\label{eq:positive-resolvent}
 a_+=\frac{A_1t}{1-\chi(\gamma)},\qquad
 b_+=(p^{-1}\sigma-1)^{-1}A_0-\sum_{m\geq0}(p^{-1}\varphi)^m q.
\end{equation}
The constant inverse exists because $\sigma^n=1$ and $p^{-n}\ne1$. On each closed disc, the $m$th summand has norm $O(p^{-m})$: after enough iterations $\varphi^m(X)=\exp(p^mt)-1$, and the double zero of $q$ compensates for $p^{-m}$ as a scalar. Telescoping proves that these terms solve \eqref{eq:primitive-untwisted} for $g^+$.

For the negative part put
\begin{equation}\label{eq:negative-resolvent}
 b_- =\sum_{j\geq1}p^j\psi^j(g^-),\quad
 c_-=g^--(p^{-1}\varphi-1)b_-,\quad
 a_-=(1-\gamma)^{-1}c_-.
\end{equation}
The negative Laurent part has bounded coefficients. The trace formula shows that $\psi$ preserves negative Laurent series and has integral Gauss norm at most one. For a fixed inner radius $r$ sufficiently close to one, its annular estimate is
\[
 \|\psi h\|_r\leq p\|h\|_{r^{1/p}}.
\]
The estimate follows by writing $\varphi\psi h=p^{-1}\sum_{\zeta^p=1}h(\zeta(1+X)-1)$ on that annulus. Negative Laurent series have decreasing norms as the radius increases. It follows that $p^j\psi^j(g^-)$ is uniformly bounded at radius $r$, while its bounded-coefficient norm is at most $p^{-j}\|g^-\|_1$. Log-convexity of annular norms gives geometric decay at any radius $\rho$ with $r<\rho<1$. The series for $b_-$ therefore converges in $\mathcal E_K^\dagger$ after a possible increase of the inner radius.

Telescoping gives $(1-p\psi)b_-=p\psi(g^-)$ and consequently $\psi(c_-)=0$. The inverse in \eqref{eq:negative-resolvent} exists in $\mathcal E_K^\dagger$ by the overconvergent inverse theorem \cite[Proposition I.5.1 and Proposition III.3.2(iii)]{CC99}, applied to the trivial coefficient module. Thus $a_-,b_-$ solve \eqref{eq:primitive-untwisted} for $g^-$. This step uses the trivial module, whose $\gamma$-action is untwisted; no \emph{\'etale} assertion about a module with Frobenius $p^{-1}\varphi$ is needed.

Finally set
\[
 \mathcal T_\gamma(f)=
 \left(\operatorname{Av}_1\partial(a_++a_-),\quad
 c_rh_*+\operatorname{Av}_1\partial(b_++b_-)\right).
\]
Apply $\partial$ to \eqref{eq:primitive-untwisted}, average, and restore the subtracted term in \eqref{eq:primitive-construction}. The result is \eqref{eq:primitive-equation}. Every operation is $\Qp$-linear once the trace-one choice is fixed. This completes the construction. The derivative reduction is the one underlying \cite[Proposition 2.12]{Liu08}; the two displayed resolvents specify its primitives for the present coefficients.
\end{proof}

The preceding theorem and the nonzero valuation residue prove directly that
\[
 H^2(\mathcal R_K(1))\xrightarrow{\ \Tr\operatorname{res}\ }\Qp
\]
is an isomorphism. To normalize this isomorphism arithmetically, use the unramified character $\lambda_0$ of \eqref{eq:epsilon}. Its Herr representative is $(-\epsilon\vartheta,0)$. This follows from the unramified period equation $\varphi(\tau)-\tau=\epsilon\vartheta$ in \eqref{eq:tau}. The cup product with the valuation pair \eqref{eq:valuation-integral} has residue trace $\epsilon c\log\chi(\gamma)$. Since \eqref{eq:reciprocity-convention} gives $\kappa(p)\cup\lambda_0$ invariant one, the normalized Tate invariant in these coordinates is
\begin{equation}\label{eq:normalized-residue}
 \operatorname{Inv}_{K,\gamma}([f v_{\Tate}])
 =\frac{\Tr_{K/\Qp}\operatorname{res}(f)}{\epsilon c\log\chi(\gamma)}.
\end{equation}
This calibration fixes both the Frobenius sign and the generator factor. Formula \eqref{eq:normalized-residue} is the standard Herr invariant in the cochain convention \eqref{eq:herr-formal}.

\subsection{Constructing every trace-zero dual class}

Return to the basis family of \Cref{sec:tracefamily}. Index its parameters by $0\leq i<r$, with $q_0=p$, and use the integral coefficients $(u_i,v_i)$ of \eqref{eq:integral-Kummer-matrices}. In the formulas of this subsection, a scalar named $v_i$ is a Kummer \emph{cochain entry}; the toric basis vector is written $\mathbf v_i$ to distinguish the two. Write the quotient vector as $w$, its Tate-dual as $a$, and the dual vectors to $\mathbf v_i$ as $b_i$.

Let $\lambda^0,\ldots,\lambda^{r-1}$ be the basis of $W^*$ dual to $q_0,\ldots,q_{r-1}$. Choose bounded overconvergent rank-one Herr representatives
\begin{equation}\label{eq:character-inputs}
 z^k=(c^k,d^k)\in(\mathcal E_K^\dagger)^{\Delta_K}\oplus(\mathcal E_K^\dagger)^{\Delta_K},
 \qquad (1-\gamma)c^k+(\varphi-1)d^k=0,
\end{equation}
for these characters. These are the ordinary rank-one coefficient inputs of the construction, with no semistability condition. Their existence is the overconvergent Herr comparison \cite{Her01,Liu08}; it is not inferred from a dimension count. To construct them, one may apply the coefficient functor to the two-dimensional representation $g\mapsto\left(\begin{smallmatrix}1&\lambda^k(g)\\0&1\end{smallmatrix}\right)$ and take the Frobenius and $\gamma$ differences of a lift of its quotient frame. The coefficient equivalence and overconvergence in \cite[Sections I.3--I.4 and Proposition III.3.2]{CC99} justify this rank-one procedure. Different choices change the representatives below by coboundaries. For example, over $\Qp$ with $q_1=1+p$, one can take $z^0=(-\epsilon,0)$ and
\[
 z^1=\left(0,-\frac{\epsilon\log\chi(\gamma)}{\log(1+p)}\right).
\]
The latter formula uses $\chi(\operatorname{rec}(u))=u^{-\epsilon}$ on principal units in the convention \eqref{eq:reciprocity-convention}.

\begin{theorem}\label{thm:fulldualrepresentatives}
For $B=(B_{ik})\in\slie_{\Qp}(W_{\Qp})$, define
\begin{equation}\label{eq:dual-obstruction-robba}
 c_i=\sum_kB_{ik}c^k,\qquad d_i=\sum_kB_{ik}d^k,\qquad
 S_B=\sum_i\bigl(u_i\varphi(d_i)-v_i\gamma(c_i)\bigr).
\end{equation}
Then $\Tr\operatorname{res}(S_B)=0$. If $(a_B,b_B)=\mathcal T_\gamma(S_B)$, the pair
\begin{equation}\label{eq:full-dual-cochain}
 Y_B=\left(\sum_i c_i b_i+a_Ba,\quad \sum_i d_i b_i+b_Ba\right)
\end{equation}
is an actual Herr cocycle in $D_{\mathrm{rig}}^\dagger(V^\dagger)$. Its class is the one corresponding to $B$ in \Cref{thm:integraltrace}. In particular, \eqref{eq:full-dual-cochain} represents all $r^2-1$ dual dimensions.
\end{theorem}
\begin{proof}
The trace of the residue of $S_B$, divided by the constant in \eqref{eq:normalized-residue}, is
\[
 \sum_i\lambda_i(q_i)=\sum_iB_{ii}=\tr B,
 \qquad \lambda_i=\sum_kB_{ik}\lambda^k.
\]
This uses the rank-one Kummer/character pairing \eqref{eq:reciprocity-convention}, the Herr cup comparison \cite{Her01,Liu08}, and the normalization \eqref{eq:normalized-residue}. Thus the trace-zero condition is exactly the degree-two obstruction, not an additional analytic assumption. \Cref{thm:explicitprimitive} provides the scalar entries.

Apply the differential to \eqref{eq:full-dual-cochain} using \eqref{eq:Robba-dual-matrices}. The $b_i$-coordinates vanish by \eqref{eq:character-inputs}. The $a$-coordinate is
\[
 (1-\chi(\gamma)\gamma)a_B+(\varphi-1)b_B
 -\sum_i\bigl(u_i\varphi(d_i)-v_i\gamma(c_i)\bigr)=0.
\]
This proves the cocycle assertion directly. Its quotient coordinates are the characters $\lambda_i$. In the long exact sequence for $V^\dagger$, the map to these quotient coordinates is injective and has image the trace-zero condition, by \Cref{thm:integraltrace}. Hence the class is exactly $B$. This argument identifies every class and does not restrict to the semistable characters. In particular, it constructs the $r(r-1)$ additional dimensions found in \Cref{cor:dimensionobstruction}.
\end{proof}

\subsection{Residue reciprocity and the integral lattices}

Let $A=(A_{ki})$ be an endomorphism of $W_{\Qp}$. Put
\begin{equation}\label{eq:primal-Herr-representative}
 U_i=\sum_k A_{ki}u_k,\qquad V_i=\sum_k A_{ki}v_k,\qquad
 X_A=\left(\sum_iU_i\mathbf v_i,\quad\sum_iV_i\mathbf v_i\right).
\end{equation}
This is the Kummer representative of $[A]\in\pgl_{\Qp}(W_{\Qp})$. The scalar matrix is a coboundary: it is the differential of $w$. Thus the displayed class depends only on $[A]$.

\begin{theorem}[Full residue reciprocity]\label{thm:fullresiduereciprocity}
With the representatives \eqref{eq:full-dual-cochain} and \eqref{eq:primal-Herr-representative}, the full cup coefficient is
\begin{equation}\label{eq:full-residue-formula}
 Z_{A,B}=\sum_i\bigl(U_i\varphi(d_i)-V_i\gamma(c_i)\bigr)\in\mathcal R_K,
 \qquad
 \operatorname{Inv}_{K,\gamma}([Z_{A,B}v_{\Tate}])=\tr(BA).
\end{equation}
The integral coefficient lattices are the module of \Cref{prop:integralRobba} and its Tate dual \eqref{eq:Robba-dual-matrices}. The cohomology lattices in the residue formula are
\[
 \{[A]:A\in\End_{\Zp}(W)\}=H^1(K,T),\qquad
 \{B\in\End_{\Zp}(W):\tr B=0\}=H^1(K,T^\dagger).
\]
On these lattices \eqref{eq:full-residue-formula} is $\Zp$-valued and perfect, including when $p\mid r$.
\end{theorem}
\begin{proof}
Use the actual module matrices \eqref{eq:Robba-matrices} and \eqref{eq:Robba-dual-matrices} in \eqref{eq:cupformula}. The toric vectors pair trivially with $a$, so the two primitive coordinates $a_B,b_B$ do not enter the cup coefficient. Pairing the $b_i$-coordinates gives exactly $Z_{A,B}$, with its displayed order and sign. In particular, this is a cup coefficient inside the unlocalized Robba ring. There is no separate residue operation on the two logarithmic summands of \eqref{eq:ZplusC}.

Apply \eqref{eq:normalized-residue} and distribute the sums. The rank-one pairing gives
\[
 \operatorname{Inv}_{K,\gamma}([Z_{A,B}v_{\Tate}])
 =\sum_{i,k}A_{ki}\lambda_i(q_k)
 =\sum_{i,k}A_{ki}B_{ik}=\tr(BA).
\]
The coefficient-lattice assertion was proved in \Cref{prop:integralRobba}. The two cohomology-lattice assertions are the integral identifications in \Cref{thm:integraltrace}. Their perfectness, and their distinction from a trace-zero self-pairing, follow from \Cref{prop:tracedefect}. These identifications show integrality of the invariant without asserting that coefficientwise integration in \Cref{thm:explicitprimitive} preserves every integral cochain. The latter operation need not preserve integral coefficients.
\end{proof}

For \Cref{ex:three}, the eight projective classes and their eight trace-zero dual classes therefore admit the representatives \eqref{eq:primal-Herr-representative} and \eqref{eq:full-dual-cochain}. In the dual bases of that example, the normalized residue matrix is the $8\times8$ identity. The self trace pairing on the trace-zero lattice still has discriminant divisible by $3$. The residue realization preserves the precise distinction already visible in the integral matrix calculation.

\subsection{Relative lifts for general \texorpdfstring{$1$}{1}-motives}\label{subsec:relativeHerr}

The scalar solution operator of \Cref{thm:explicitprimitive} applies beyond the toric basis family. It is the target of the lifting equation, rather than every graded coefficient of the representation, that must be a Tate module. We give the construction in a form that specifies its input.

Let $M=[Y\to G]$ be a $1$-motive over the unramified field $K$, with $Y=\mathbf Z^s$, and put $U=V_p(G)$. The following statement only requires the corresponding exact sequence of Galois representations. Write $D_{\mathrm{rig}}$ for the usual overconvergent Robba coefficient functor. It gives an exact sequence
\begin{equation}\label{eq:relative-coefficients}
0\longrightarrow\mathcal R_K(1)^s\longrightarrow
D_{\mathrm{rig}}(V_M^\dagger)\longrightarrow
D_{\mathrm{rig}}(U^\dagger)\longrightarrow0.
\end{equation}
The functor, its exactness, and the Herr comparison are established inputs \cite{CC99,Her01,Liu08}. Choose an $\mathcal R_K$-linear section $s_0$ of the last arrow. Averaging $\delta s_0\delta^{-1}$ over the torsion subgroup $\Delta$ makes the section $\Delta$-equivariant; denote this section by $s$. It need not commute with Frobenius or $\gamma$.

For a Herr cocycle $z=(c,d)$ in $D_{\mathrm{rig}}(U^\dagger)^\Delta$, lift its two entries by $s$ and define $S_j(z)\in\mathcal R_K$ by
\begin{equation}\label{eq:relative-error}
-d^1(s(c),s(d))
=\sum_{j=1}^s S_j(z)v_{\Tate,j},
\end{equation}
where $v_{\Tate,j}$ is the standard basis of the $j$th Tate summand. The expression lies in that submodule because $z$ is a cocycle. Each scalar coefficient has the required Tate $\Delta$-invariance. The module matrices of \eqref{eq:relative-coefficients} determine \eqref{eq:relative-error} by multiplication and subtraction.

\begin{theorem}[A relative construction of all dual Herr classes]\label{thm:relativeHerrlift}
Under the chosen coefficient data, the class of $z$ lifts to $H^1(K,V_M^\dagger)$ if and only if
\begin{equation}\label{eq:relative-residue-condition}
\Tr_{K/\Qp}\operatorname{res}(S_j(z))=0
\qquad(1\leq j\leq s).
\end{equation}
These conditions are precisely $\delta_M([z])=0$ from \eqref{eq:mixed-full-obstruction}. When they hold, compute
\[
(a_j,b_j)=\mathcal T_\gamma(S_j(z))
\]
with the explicit operator of \Cref{thm:explicitprimitive}. Then
\begin{equation}\label{eq:relative-Herr-lift}
Y_z=\left(s(c)+\sum_j a_jv_{\Tate,j},\
           s(d)+\sum_j b_jv_{\Tate,j}\right)
\end{equation}
is an actual Herr cocycle lifting $z$. Together with the images of the $s$ Tate-valued degree-one complexes, these cocycles represent all of $H^1(K,V_M^\dagger)$. The ambiguity in a lift is exactly the first term of \eqref{eq:mixed-dual-exact}.
\end{theorem}
\begin{proof}
The connecting class of $z$ in \eqref{eq:relative-coefficients} is represented by $d^1(s(c),s(d))$, which is minus \eqref{eq:relative-error}. The sign computation in \Cref{prop:fullmixedobstruction}, together with compatibility of the Herr and Galois connecting maps, gives
\[
\operatorname{Inv}_{K,\gamma}([S_j(z)v_{\Tate,j}])
=\delta_M([z])(y_j).
\]
By \eqref{eq:normalized-residue}, the invariant is a nonzero fixed multiple of $\Tr\operatorname{res}(S_j(z))$. Thus \eqref{eq:relative-residue-condition} is equivalent to vanishing of the connecting class and to the lifting criterion. For each vanishing residue, \Cref{thm:explicitprimitive} gives $d^1(a_jv_{\Tate,j},b_jv_{\Tate,j})=S_j(z)v_{\Tate,j}$. Adding these equalities to \eqref{eq:relative-error} proves that \eqref{eq:relative-Herr-lift} is closed.

Every full dual class projects to a class in $\ker\delta_M$. After subtracting a lift constructed in \eqref{eq:relative-Herr-lift}, its projection is a coboundary. Lift the corresponding degree-zero entry and subtract its coboundary as well. What remains is a cocycle in $\mathcal R_K(1)^s$. Quotienting its classes by the connecting image from $H^0(K,U^\dagger)$ is exactly \eqref{eq:mixed-dual-exact}. This proves the exhaustiveness and the stated ambiguity.
\end{proof}

\begin{corollary}[Compatibility with full local duality]\label{cor:relative-full-pairing}
Let $x\in H^1(K,U)$ and let $Y_z$ be a lift constructed in \eqref{eq:relative-Herr-lift}. If $j:U\to V_M$ is the inclusion, then
\begin{equation}\label{eq:relative-full-pairing}
\pair{j_*x}{[Y_z]}=\inv_K(x\cup[z]).
\end{equation}
In coefficient coordinates the scalar correction terms in \eqref{eq:relative-Herr-lift} contribute zero to this cup product. The right side uses the full pairing for the semi-abelian part, including the Poincar\'e pairing on its abelian graded quotient.
\end{corollary}
\begin{proof}
The submodule $U$ evaluates to zero against $Y_{\Qp}^*(1)$ in the dual extension. Pairing equivariance therefore removes the correction entries from both terms of \eqref{eq:cupformula}. The remaining cup is the evaluation cup for $U,U^\dagger$. Its normalized invariant is the arithmetic Tate pairing by the Herr comparison \cite{Her01,Liu08}. This proves \eqref{eq:relative-full-pairing}, equivalently the usual functoriality of local duality.
\end{proof}

In the toric basis family, \eqref{eq:relative-error} is exactly $S_B$ in \eqref{eq:dual-obstruction-robba}, and \Cref{thm:relativeHerrlift} specializes to \Cref{thm:fulldualrepresentatives}. In that family the first term of \eqref{eq:mixed-dual-exact} is zero, since its parameter classes form a basis of $W_{\Qp}$. Hence the quotient matrix specifies a unique dual cohomology class. For a nonsplit semi-abelian group, the additional input is its actual coefficient sequence \eqref{eq:relative-coefficients} and the chosen quotient cocycles. The relative theorem supplies the explicit Tate-valued correction once that input is given. The finite-dimensional maps $(R,\lambda)$ determine semistable extension cohomology by \Cref{thm:mixednormalform}; they do not supply the full abelian coefficient matrices. The integral lattice assertions of \Cref{prop:integralRobba,thm:fullresiduereciprocity} remain those proved in their stated toric setting, with the product extension of \Cref{cor:integralproduct}. Laurent integration in $\mathcal T_\gamma$ can introduce denominators, so \eqref{eq:relative-Herr-lift} is asserted over the rational Robba ring.

\section{Integral logarithms and divided powers}\label{sec:dividedpowers}

The integral trace defect of \Cref{prop:tracedefect} has a different origin from denominators in logarithmic integration. We record the latter exactly. For a $\Zp$-algebra $A_0$ without $p$-torsion,
\begin{equation}\label{eq:derivativecoker}
\coker\bigl(\partial_\ell:A_0[\ell]\to A_0[\ell]\bigr)
=\bigoplus_{m\geq1}A_0/mA_0\cdot\ell^{m-1}.
\end{equation}
This follows coefficient by coefficient from $\partial_\ell(\ell^m)=m\ell^{m-1}$. In particular no bounded power of $p$ annihilates the entire cokernel when $A_0=\Zp$.

Let $A_0\langle\ell\rangle$ be the divided-power polynomial algebra: it is free over $A_0$ on $\ell^{[m]}$, with
\[
\ell^{[i]}\ell^{[j]}=\binom{i+j}{i}\ell^{[i+j]},\qquad
\partial_\ell\ell^{[m]}=\ell^{[m-1]}\quad(m\geq1).
\]
The derivative is surjective and has kernel $A_0$.

\begin{proposition}\label{prop:pdhomotopy}
Let $D_0$ be finite free over $A_0$ with nilpotent endomorphism $N_0$. On $A_0\langle\ell\rangle\otimes D_0$, the horizontalization and homotopy of \Cref{prop:horizontal} are defined integrally by replacing $e^{\ell N_0}$ with
\[
\sum_{j\geq0}\ell^{[j]}N_0^j
\]
and integration by $I(\ell^{[m]})=\ell^{[m+1]}$. They satisfy the same two homotopy identities. A translation $\ell\mapsto\ell+c$ preserves this construction whenever $c$ has compatible divided powers in the coefficient ring.
\end{proposition}
\begin{proof}
Nilpotence makes each displayed exponential finite. The divided-power multiplication formula gives the product identity for these exponentials and the differentiation identity $\partial_\ell e^{\ell N_0}=N_0e^{\ell N_0}$. Moreover $\partial_\ell I=1$ and $I\partial_\ell=1-\ev_0$ on divided-power polynomials. The proof of \Cref{prop:horizontal} therefore applies word for word without division by an integer. For translations, the identity
$(\ell+c)^{[m]}=\sum_{i+j=m}\ell^{[i]}c^{[j]}$
gives the asserted compatibility precisely under the stated divided-power hypothesis.
\end{proof}

This proposition does not equip every $p$-adic shift with divided powers. In particular a rational change of logarithmic coordinate need not preserve a chosen integral period algebra. Both the trace defect $O/rO$ and the integration cokernel \eqref{eq:derivativecoker} must be tracked if one later proposes an integral regulator normalization.

\section{Division towers, Lie cohomology, and the arithmetic Iwasawa pairing}\label{sec:towers}

Return to the basis family and define
\[
K^{\mathrm{cyc}}=K(\mu_{p^\infty}),\qquad
\mathcal L=K^{\mathrm{cyc}}(q_i^{1/p^m}:0\leq i\leq n,\ m\geq1),
\]
\[
H=\operatorname{Gal}(\mathcal L/K^{\mathrm{cyc}}),\qquad
G=\operatorname{Gal}(\mathcal L/K).
\]
This is the division tower of the toric $1$-motive after adjoining the cyclotomic torsion.

\begin{proposition}\label{prop:towergroup}
The group $H$ is isomorphic to $\Zp^r$ and $G$ has dimension $r+1$. An open subgroup is a semidirect product $H\rtimes\Gamma_0$, with $\Gamma_0\simeq\Zp$ and conjugation $\gamma h\gamma^{-1}=\chi(\gamma)h$ in additive coordinates on $H$.
\end{proposition}
\begin{proof}
Kummer theory embeds $H$ as a closed submodule of $\Zp(1)^r$. If its rank were smaller than $r$, a nonzero $\Qp$-linear combination of the Kummer classes of the $q_i$ would restrict to zero over $K^{\mathrm{cyc}}$. Inflation--restriction places such a class in $H^1(\operatorname{Gal}(K^{\mathrm{cyc}}/K),\Qp(1))$. This group vanishes: on a procyclic open subgroup its degree-one cohomology is $\Qp(1)/(\gamma-1)\Qp(1)=0$ because $\chi(\gamma)-1\ne0$, and passage across the finite quotient does not change the vanishing. The Kummer basis is independent over $\Qp$, a contradiction. A closed rank-$r$ submodule of $\Zp^r$ has finite index and is abstractly $\Zp^r$.

Restrict to the inverse image of a pro-$p$ cyclotomic subgroup $\Gamma_0$. Lift a topological generator to $G$. Its procyclic closure is $\Zp$ and maps isomorphically onto $\Gamma_0$, so it gives a section. The Kummer cocycle identity gives the displayed conjugation action. This proves the dimension and semidirect product statements. No assertion about index one in a chosen root-coordinate lattice is needed.
\end{proof}

Let $E/\Qp$ be a finite coefficient field. The normalized Lie algebra of this open subgroup is
\begin{equation}\label{eq:liealgebra}
\mathfrak g=Eh\ltimes\mathfrak n,\qquad
\mathfrak n=\bigoplus_{i=0}^n Ex_i,\qquad[h,x_i]=x_i.
\end{equation}

\begin{proposition}\label{prop:characters}
Every one-dimensional locally analytic character of $G$ has differential zero on $\mathfrak n$ and has finite image on $H$. Scalar-character evaluations therefore do not detect every analytic distribution on $G$, where distributions are continuous linear functionals on locally analytic functions, with convolution as product.
\end{proposition}
\begin{proof}
A differential into an abelian Lie algebra kills $[\mathfrak g,\mathfrak g]=\mathfrak n$. At group level a scalar character kills $(\chi(\gamma)-1)H$, which is open in $H$. Its image on $H$ is therefore finite. The nonzero differential distribution $x_i$ evaluates to zero at every such character, giving the final assertion.
\end{proof}

For $k\in E$, write $E(k)$ for the Lie algebra character with $h$ acting by $k$ and $\mathfrak n$ acting trivially. For integral $k$ it is the differential of $\chi^k$.

\begin{proposition}\label{prop:liecohom}
For \eqref{eq:liealgebra},
\begin{equation}\label{eq:liecohom}
H^q(\mathfrak g,E(k))=
\begin{cases}
\bigwedge^q\mathfrak n^*\otimes E(k),&k=q,\\
h^*\wedge\bigwedge^{q-1}\mathfrak n^*\otimes E(k),&k=q-1,\\
0,&\text{otherwise},
\end{cases}
\end{equation}
with exterior powers outside $0,\ldots,r$ interpreted as zero. In particular
$H^{r+1}(\mathfrak g,E(1))=0$ since $r\geq2$, whereas
$H^{r+1}(\mathfrak g,E(r))=E$.
\end{proposition}
\begin{proof}
Decompose the Chevalley--Eilenberg complex as
\[
C^q=\bigl(\textstyle\bigwedge^q\mathfrak n^*\otimes E(k)\bigr)
\oplus\bigl(h^*\wedge\textstyle\bigwedge^{q-1}\mathfrak n^*\otimes E(k)\bigr).
\]
On the first summand in degree $q$, its differential is $(k-q)h^*\wedge-$: the coefficient action contributes $k$ and the $q$ brackets contribute $-q$. The differential of the second summand is zero. Each pair is acyclic unless its scalar $k-q$ vanishes, which proves the formula.
\end{proof}

The orientation coefficient of Lie algebra duality is $\det\mathfrak g=E(r)$ because $\tr(\operatorname{ad}h)=r$. More generally the cochain wedge pairing with values in $\det\mathfrak g$ yields
\[
H^q(\mathfrak g,M)\times
H^{r+1-q}(\mathfrak g,M^*\otimes\det\mathfrak g)\longrightarrow E
\]
as a perfect pairing for finite-dimensional $M$. The coefficient action cancels the bracket term in top degree, so the paired differentials are adjoint and cochain perfectness descends to cohomology. This is classical twisted Lie algebra duality \cite{Haz70}. Its orientation character must be distinguished from the arithmetic Tate twist $E(1)$. A comparison with period coefficients must account for their difference.

\subsection{The arithmetic inverse-limit pairing}

Take finite Galois subextensions $K_j/K$ cofinal in $\mathcal L/K$, set $G_j=\operatorname{Gal}(K_j/K)$, and define
\[
H^1_{\mathrm{Iw}}(\mathcal L/K,T)=\varprojlim_j H^1(K_j,T),\qquad
\Lambda(G)=O[[G]],
\]
using corestriction. For $x=(x_j)$ and $y=(y_j)$ in the two corresponding inverse limits put
\begin{equation}\label{eq:Iwpairlevel}
P_j(x_j,y_j)=\sum_{g\in G_j}\inv_{K_j}(x_j\cup gy_j)\,g.
\end{equation}

\begin{proposition}\label{prop:Iwpair}
The elements \eqref{eq:Iwpairlevel} are compatible and define $P_{\mathrm{Iw}}(x,y)\in\Lambda(G)$. If $\iota$ is the anti-involution $g\mapsto g^{-1}$, then
\[
P_{\mathrm{Iw}}(a x,y)=aP_{\mathrm{Iw}}(x,y),\qquad
P_{\mathrm{Iw}}(x,a y)=P_{\mathrm{Iw}}(x,y)\iota(a)
\quad(a\in\Lambda(G)).
\]
\end{proposition}
\begin{proof}
For $K_k\supset K_j$, write $J=\operatorname{Gal}(K_k/K_j)$ and choose a lift $\widetilde g$ of $g\in G_j$. The coefficient at $g$ of the projection of $P_k$ is
\[
\sum_{u\in J}\inv_{K_k}(x_k\cup\widetilde g u y_k)
=\inv_{K_j}(\operatorname{cor}x_k\cup g\operatorname{cor}y_k).
\]
The equality uses $\operatorname{res}\operatorname{cor}y_k=\sum_{u\in J}uy_k$ and restriction--corestriction adjointness of the local invariant. It equals the coefficient of $P_j$. The inverse limit therefore exists. For group elements $a$, reindexing $g$ by $a^{-1}g$ in the first argument gives left multiplication by $a$, and reindexing by $ga$ in the second gives right multiplication by $a^{-1}$. Linearity and continuity extend the identities to $\Lambda(G)$.
\end{proof}

The proposition constructs an arithmetic pairing; it does not assert perfectness on these degree-one inverse limits. In a noncommutative tower, derived duality and its degree shifts require a separate treatment.

\subsection{An integral period realization throughout the division tower}\label{subsec:towerperiods}

The division tower contains the cyclotomic extension, so its cohomology can be realized by a compatible system of cyclotomic coefficient modules. The transition maps are module traces. Cup products are formed in the finite-base coefficient rings, and their local invariants are then passed to the projective limit.

For a finite extension $F/K$ inside $\mathcal L$, set $F^{\mathrm{cyc}}=F(\mu_{p^\infty})$, $H_F=\operatorname{Gal}(\overline K/F^{\mathrm{cyc}})$, and
\[
 \mathcal A_F=\mathcal A^{H_F},\qquad
 D_F(T)=(\mathcal A\otimes_{\Zp}T)^{H_F}.
\]
Here $\mathcal A$ is Fontaine's integral unramified coefficient ring: it is the $p$-adic completion of the maximal unramified extension of the cyclotomic coefficient ring inside the Witt-vector coefficient ring. These definitions and the integral coefficient equivalence are recalled in \cite[Sections I.1--I.3]{CC99}. The operator $\psi_F$ on $D_F(T)$ is the normalized Frobenius trace. Thus $\psi_F\varphi=1$. Its normalization is Frobenius-theoretic; the transition between two finite fields will be an ordinary field trace.

Equip $\mathcal A_F$ and $D_F(T)$ with Fontaine's weak coefficient topology. The space $D_F(T)^{\psi=1}$ is compact in this topology by \cite[Proposition I.5.6(i)]{CC99}. Projective limits of these spaces below carry their compact projective-limit topology.

Choose a cofinal increasing sequence of finite Galois fields $F_m/K$ in $\mathcal L$, with $\mu_p\subset F_m$. All cyclotomic groups $\Gamma_{F_m}$ are then torsion free. If $F'\supset F$, the restriction identity
\[
 D_{F'}(T)=\mathcal A_{F'}\otimes_{\mathcal A_F}D_F(T)
\]
holds over the common cyclotomic coefficient fields, with the usual interpretation when $F'^{\mathrm{cyc}}=F^{\mathrm{cyc}}$. For the toric lattices in \Cref{prop:integralRobba}, this gives their coefficient modules at every level by scalar extension of the displayed integral matrices.

We specify the transition $\operatorname{Tr}^{\mathrm{Iw}}_{F'/F}$ on $\psi=1$ vectors. First replace $F$ by the finite intersection $F_0=F'\cap F^{\mathrm{cyc}}$. The step $F_0/F$ is a shift of cyclotomic base, and its transition is the identity on $\psi=1$ vectors. For the remaining step choose a sufficiently deep finite cyclotomic level $F_{0,N}$ such that $F'\cap F^{\mathrm{cyc}}\subset F_{0,N}$. At the compositum $F'F_{0,N}$, the two cyclotomic groups agree. The ordinary coefficient trace
\[
 \Tr_{\mathcal A_{F'}/\mathcal A_{F_0}}\otimes1
\]
then defines the transition. Equivalently, use that trace after this common cyclotomic base change and then undo the cyclotomic shift. The trace is integral, commutes with $\varphi$ and $\psi$, and is independent of increasing $N$. This prescription includes the case of a purely cyclotomic field step, for which no extra factor $[F':F]$ is inserted.

\begin{theorem}\label{thm:towerperiodrealization}
For every finite free $\Zp$-representation $T$ of $G_K$, there is a canonical integral isomorphism
\begin{equation}\label{eq:tower-period-isomorphism}
 \mathcal R_{\mathcal L,T}:H^1_{\mathrm{Iw}}(\mathcal L/K,T)
 \xrightarrow{\ \sim\ }
 \mathcal P_{\mathcal L}(T):=
 \varprojlim_m\bigl(D_{F_m}(T)^{\psi=1},\operatorname{Tr}^{\mathrm{Iw}}_{F_{m+1}/F_m}\bigr).
\end{equation}
It is independent of the cofinal sequence. For $V=T\otimes\Qp$, its rational target is $\mathcal P_{\mathcal L}(T)\otimes\Qp$. In particular, rationalization means a single denominator valid at every level.
\end{theorem}
\begin{proof}
At each finite base $F$, Fontaine's cyclotomic isomorphism, proved in \cite[Theorem II.1.3]{CC99}, is
\begin{equation}\label{eq:cyclotomic-psi-comparison}
 H^1_{\mathrm{Iw}}(F^{\mathrm{cyc}}/F,T)\simeq D_F(T)^{\psi=1}.
\end{equation}
Its compatibility with finite field corestriction can be checked before taking limits. After making the cyclotomic groups agree, use the finite free coefficient restriction identity above. Representatives of the relevant finite Galois cosets can then be chosen in $H_F$. Summing their actions on an $H_{F'}$-invariant period vector is the ordinary trace on its scalar factor. This trace commutes with Frobenius and the common group action, and hence gives a map of Herr complexes. To identify its cohomological map, compute a Galois cocycle using a period primitive as in \cite[Proposition I.4.1]{CC99}. The corestriction formula replaces that primitive by the sum of its coset translates. This is exactly the coefficient trace; the same primitive computation is written out for corestriction in \cite[proof of Lemma II.2.1]{CC99}. Thus coefficient trace induces field corestriction. For a cyclotomic shift, the cochain maps are those of \cite[Lemma II.2.1 and Remark II.2.2]{CC99}; on the $\psi=1$ parameter their transition is the identity. These two observations prove the trace prescription preceding the theorem.

Finite fields of the form $F_m(\mu_{p^N})$ are cofinal among the finite fields in $\mathcal L$. Taking the projective limit of \eqref{eq:cyclotomic-psi-comparison} over $m$ therefore gives \eqref{eq:tower-period-isomorphism}. Projective limits commute, so the order of the two field indices has no effect. The same cofinality proves independence of the sequence. Tensoring the resulting integral limit by $\Qp$ gives the stated rationalization. An unrestricted projective limit of $D_{F_m}(V)^{\psi=1}$ could contain systems with unbounded denominators; that larger space is not used in this theorem.
\end{proof}

The result is a period realization of the full division-tower cohomology, with integral traces and a specified $\psi$ operator. Its input \eqref{eq:cyclotomic-psi-comparison} is classical. The passage to the division tower is the trace-compatible projective limit just proved.

\subsection{A cochain formula for the period pairing}\label{subsec:periodtowerpairing}

We record finite specializations explicitly. For a finite $F$ containing $\mu_p$, let $r_F$ be defined by $\log\chi(\Gamma_F)=p^{r_F}\Zp$, and choose the unique generator $\gamma_F$ with $\log\chi(\gamma_F)=p^{r_F}$. This normalization makes the generator unit $\ell_F(\gamma_F)=\log\chi(\gamma_F)/p^{r_F}$ of \cite[Section I.4]{CC99} equal to one. It removes a unit factor from their canonical cyclotomic comparison.

For $u\in D_F(T)^{\psi=1}$ put
\begin{equation}\label{eq:psi-Herr-specialization}
 A_F(u)=(\gamma_F-1)^{-1}(\varphi-1)u,\qquad
 h_F(u)=(A_F(u),u).
\end{equation}
The inverse is on $D_F(T)^{\psi=0}$, since $\psi(\varphi-1)u=0$. Thus the formula is integral. After tensoring with $\Qp$, it is overconvergent by \cite[Proposition III.3.2(ii)--(iii)]{CC99}. It satisfies the Herr cocycle equation directly. The map from the Frobenius Herr complex to the $\psi$ complex is $(a,b)\mapsto(-\psi a,b)$ in degree one, as in \cite[Lemma I.5.2]{CC99}. It sends \eqref{eq:psi-Herr-specialization} to $(0,u)$, proving that its class is the specialization of \eqref{eq:cyclotomic-psi-comparison}.

If $\gamma'=\gamma_F^p$, the corestriction formula is
\begin{equation}\label{eq:cyclotomic-cochain-trace}
 (a,b)\longmapsto\bigl((1+\gamma_F+\cdots+\gamma_F^{p-1})a,b\bigr).
\end{equation}
Indeed, $(\gamma_F-1)(1+\cdots+\gamma_F^{p-1})=\gamma'-1$. The $\psi=0$ inverse therefore makes \eqref{eq:psi-Herr-specialization} compatible with this map. This is the explicit normalization check behind the cyclotomic shift used in \Cref{thm:towerperiodrealization}.

For an arbitrary finite $F$, denote the normalized Herr invariant on $H^2(D_F(\Zp(1)))$ by $\operatorname{Inv}_F$. It is the local invariant under the integral Herr comparison \cite{Her01}; equivalently, it is calibrated by the Kummer class of a uniformizer and the valuation character, as in \eqref{eq:reciprocity-convention}. Over our original unramified field it is \eqref{eq:normalized-residue}. At ramified levels one uses the residue and coefficient trace for $\mathcal A_F$; the unramified numerical constant in \eqref{eq:normalized-residue} is not reused without this normalization.

\begin{theorem}\label{thm:periodtowerreciprocity}
Let $u\in\mathcal P_{\mathcal L}(T)$ and $v\in\mathcal P_{\mathcal L}(T^\dagger)$. At a finite Galois level $F_j/K$ put $G_j=\operatorname{Gal}(F_j/K)$ and
\begin{align}\label{eq:period-Iw-pairing}
 \mathcal P_j(u,v)=\sum_{g\in G_j}\operatorname{Inv}_{F_j}\Bigl([
 &b(A_{F_j}(u_j),\varphi(gv_j))\notag\\
 &-b(u_j,\gamma_{F_j}A_{F_j}(gv_j))]\Bigr)\,g.
\end{align}
Here $u_j,v_j$ are the $\psi=1$ parameters at base $F_j$, and the actions on each argument are its coefficient-module actions. These elements define an integral pairing with values in $\Zp[[G]]$. If $x=\mathcal R_{\mathcal L,T}^{-1}(u)$ and $y=\mathcal R_{\mathcal L,T^\dagger}^{-1}(v)$, then
\begin{equation}\label{eq:period-Iw-adjointness}
 \mathcal P_{\mathcal L}(u,v)=P_{\mathrm{Iw}}(x,y).
\end{equation}
It has the same left-linearity and right $\iota$-linearity as \Cref{prop:Iwpair}.
\end{theorem}
\begin{proof}
The coefficient in square brackets is precisely the cup product \eqref{eq:cupformula} of $h_{F_j}(u_j)$ and $h_{F_j}(gv_j)$. By \eqref{eq:psi-Herr-specialization}, these represent $x_j$ and $gy_j$. The normalized Herr invariant is the local Tate invariant, so each coefficient of \eqref{eq:period-Iw-pairing} equals the corresponding coefficient of \eqref{eq:Iwpairlevel}. If a lift of $g$ changes by an element of $G_{F_j}$, its induced action on $H^1(F_j,T^\dagger)$ is unchanged; the displayed invariant is therefore independent of that lift. \Cref{thm:towerperiodrealization} identifies the transition traces with corestriction. The compatibility and sesquilinearity now follow from the coefficientwise calculation in \Cref{prop:Iwpair}. This proves \eqref{eq:period-Iw-adjointness}, including integrality.
\end{proof}

For the toric basis family, \Cref{prop:integralRobba} supplies explicit matrices in every coefficient module in \eqref{eq:tower-period-isomorphism}. At the original unramified level, \Cref{thm:fullresiduereciprocity} evaluates its finite-class pairing as matrix trace. \Cref{thm:periodtowerreciprocity} gives its compatibility throughout the noncommutative division tower. These statements furnish a cohomological period realization and its adjointness formula. An analytic regulator with values in a distribution algebra is a further structure, with a different target; it is discussed in \Cref{sec:remaining}.

\begin{corollary}[The division tower of a $1$-motive]\label{cor:motivicdivisiontower}
Let $M=[Y\to G]$ be a $1$-motive over $K$ and let $T_M=T_p(M)$. Define
\begin{equation}\label{eq:general-motivic-tower}
\mathcal L_M=K(T_M,\mu_{p^\infty}),
\qquad G_M=\operatorname{Gal}(\mathcal L_M/K).
\end{equation}
Here $K(T_M)$ is the fixed field of the kernel of the Tate-module action. This is a Galois $p$-adic Lie extension. If $Y$ is constant, it is the field generated by the coordinates of $G[p^\infty]$, all division points $[p^j]P=u(y)$, and the indicated cyclotomic roots. There is an integral period realization
\begin{equation}\label{eq:general-motivic-period-tower}
H^1_{\mathrm{Iw}}(\mathcal L_M/K,T_M)
\simeq\varprojlim_m D_{F_m}(T_M)^{\psi=1},
\end{equation}
with the transition maps, topology, and uniformly bounded rational denominators of \Cref{thm:towerperiodrealization}. Its pairing with the Tate-dual system is \eqref{eq:period-Iw-pairing} and equals the arithmetic Iwasawa pairing. These assertions do not require good reduction or semistability of $V_M$.
\end{corollary}
\begin{proof}
The Tate lattice is finite free, so the joint Tate and cyclotomic action has compact, hence closed, image in $\operatorname{GL}(T_M)\times\Zp^\times$. A closed subgroup of a $p$-adic linear group is a $p$-adic Lie group. Its fixed field is \eqref{eq:general-motivic-tower}. The description by torsion and division points follows from the finite exact realizations of $M$ given in \eqref{eq:finite-motive-realization}: the fiber above $y$ consists of the divisions of $u(y)$, with their toric and abelian torsion translates. Now apply \Cref{thm:towerperiodrealization,thm:periodtowerreciprocity} to $T_M$ and its Tate dual.
\end{proof}

The special toric calculation $H\simeq\Zp^r$ and its orientation character $\chi^r$ in \Cref{prop:characters,prop:liecohom} concern the multi-Kummer tower specified there. The field \eqref{eq:general-motivic-tower} contains an abelian torsion extension when $A\ne0$, whose Galois action can introduce further Lie directions. Its dimension, derived subgroup, and Lie orientation must be computed from $G_M$; they are not determined by the toric rank alone. The period-cohomological comparison \eqref{eq:general-motivic-period-tower} requires none of those special Lie-group identities.

\section{Consequences and further comparisons}\label{sec:remaining}

\subsection{From the semistable boundary to full duality}

The earlier semistable calculation and the full residue calculation describe different parts of the same local pairing. \Cref{thm:universalboundary,thm:boundary,thm:semistablereciprocity} identify the invariant monodromy boundary and its evaluation pairing with the unramified quotient. \Cref{thm:fulldualrepresentatives,thm:fullresiduereciprocity} place that calculation in the full coefficient complex. The latter construction includes the nonsemistable dual classes required by \Cref{cor:dimensionobstruction}; \Cref{thm:relativeHerrlift} extends its scalar correction to a general constant-lattice motive with specified coefficient data. In these coordinates, the trace-zero condition is equivalent to the vanishing of the cohomology class
$[S_Bv_{\Tate}]\in H^2(\mathcal R_K(1))$, or equivalently to
$\Tr_{K/\Qp}\operatorname{res}(S_B)=0$.
\Cref{thm:explicitprimitive} then constructs the scalar
cochain correction.

There are two distinct integral statements. \Cref{prop:integralRobba} identifies the integral coefficient lattices for all the toric parameter matrices under consideration. \Cref{thm:integraltrace,thm:fullresiduereciprocity} identify their full cohomology lattices and perfect residue pairing in the basis family. The algebraic integral complex of \Cref{sec:integralcomplex} is compared with integral Galois cohomology in that family by \Cref{thm:integralcomparison}; no such comparison is needed to construct the general coefficient lattice. These statements keep coefficient integrality, cohomology integrality, and the reduced semistable complex separate.

Monodromy consequently has two precise roles. The operator on the middle extension gives the invariant boundary $\partial_N$. In logarithmic coordinates it contributes the term $\mathcal C_N$ in \eqref{eq:ZplusC}. \Cref{thm:cupcorrection} proves the coordinate dependence of this term, while \Cref{thm:fullresiduereciprocity} evaluates the invariant full cup product. This provides a monodromy interpretation of the duality formula without requiring an additional coordinate-independent correction summand.

\subsection{The analytic target and known regulators}

The period target in \eqref{eq:tower-period-isomorphism} is a projective system of integral coefficient modules with normalized Frobenius traces. Its topology is the projective-limit topology, and its rationalization uses uniformly bounded denominators. The specialization \eqref{eq:psi-Herr-specialization} and pairing \eqref{eq:period-Iw-pairing} are specified on that target. Hence \Cref{thm:periodtowerreciprocity} proves a period-cohomological adjointness identity over the full division tower.

An analytic regulator into a distribution algebra would add a different target and interpolation statements. In the cyclotomic case those structures are supplied by Perrin-Riou \cite{PR01} and, for de Rham $(\varphi,\Gamma)$-modules, by Nakamura \cite{Nak14}. Locally analytic constructions over more general extensions can be compared with \eqref{eq:tower-period-isomorphism} through the period theories of Poyeton and Porat \cite{Poy25,Por24}. Such a comparison must identify its transition traces with the ones in \Cref{thm:towerperiodrealization} and recover \eqref{eq:period-Iw-adjointness}. This gives an exact compatibility statement for any proposed analytic target. Gao--Zhao's results \cite[Sections 9.2 and 10.6]{GZ25} explain why the choice of a cohomology complex matters in this passage.

The division tower has Lie directions invisible to scalar characters, by \Cref{prop:characters}. Its Lie orientation has character $\chi^r$, by \Cref{prop:liecohom}, whereas local Tate duality uses $\chi$. These calculations give independent tests of a distribution-valued or derived comparison. The projective-limit construction proves its pairing identity before making scalar specializations, so it retains every group-ring coefficient.

\subsection{The range of the \texorpdfstring{$1$}{1}-motivic extension}\label{subsec:mixedscope}

The generalization has three precise levels. \Cref{thm:universalboundary,prop:obstruction-duality} apply to every semistable representation over a finite extension of $\Qp$. \Cref{thm:mixedrealization,thm:mixednormalform,thm:boundary,thm:semistablereciprocity} use the split model with good abelian quotient and compute its semistable extensions from geometric point coordinates. \Cref{thm:relativeHerrlift} addresses full dual cohomology over our unramified field once the actual coefficient sequence is supplied. The last theorem solves the Tate-valued lifting equations even when the quotient contains abelian coefficients.

\begin{center}
\begin{longtable}{@{}>{\raggedright\arraybackslash}p{0.35\textwidth}>{\raggedright\arraybackslash}p{0.59\textwidth}@{}}
\toprule
Statement & Range and input\\
\midrule
Monodromy boundary and rank change & Any semistable $V$; use $(D/ND)^{p\bar\varphi=1}$, as in \Cref{thm:universalboundary}.\\[3pt]
\midrule
Reduced complex and exponential quotient & A $1$-motive with constant lattice, split torus, and model \eqref{eq:good-model}; use $R$ and $\lambda$ in $\operatorname{Lie}(G)$.\\[3pt]
\midrule
Ambient one-pole bound & Any semistable representation with filtration jumps in $\{-1,0\}$; \Cref{prop:generalpositive} uses the positive comparison.\\[3pt]\midrule
Logarithmic cup coordinates & Nilpotent monodromy and a compatible pairing; the first-order formula applies when $N^2=0$, hence to the mixed motives here.\\[3pt]\midrule
Full dual Herr lifts & Constant lattice $Y$ and a specified coefficient extension; solve $s$ scalar residue-zero equations by \Cref{thm:relativeHerrlift}.\\[3pt]\midrule
Integral projective/trace-zero duality & The toric basis family of \Cref{thm:integraltrace}, with the product extension and torsion hypotheses of \Cref{cor:integralproduct}.\\[3pt]\midrule
Division-tower period pairing & Any finite free Tate lattice and its Galois tower containing the cyclotomic extension; \Cref{cor:motivicdivisiontower}.\\
\bottomrule
\end{longtable}
\end{center}

For a nonsplit $G$, the Hodge quotient $\lambda$ is computed by \eqref{eq:geometric-lambda}, but the filtered crystalline module of $G$ contains the semi-abelian extension class. Its realization uses the Poincar\'e biextension. A closed formula for all coefficient matrices therefore needs abelian crystalline periods and the biextension data. The scalar functions of \Cref{sec:explicitrobba} compute the toric Kummer blocks; adjoining a scalar logarithmic variable does not by itself determine the additional matrices. Given those matrices, the relative solution in \Cref{thm:relativeHerrlift} is explicit and exhaustive. For a split product, \Cref{prop:splitproducts} gives the separate division-point cocycles and direct-sum comparison.

The integral statements require corresponding care. Rational Frobenius projectors can introduce denominators; the integral weight filtration need not split. Abelian $p$-power torsion can produce torsion in degree-one cohomology, as shown in the proof of \Cref{cor:integralproduct}. Thus the rational complex \eqref{eq:mixed-motivic-complex} is not asserted to identify the integral Galois lattice for every mixed motive. The log $p$-divisible construction supplies an integral geometric realization \cite{WZ24,BWZ26}; relating a particular integral reduced complex to that realization requires its lattice data. \Cref{thm:integralcomparison} and \Cref{cor:integralproduct} are the comparisons proved here.

When a torus or lattice becomes split over a finite unramified extension, \Cref{prop:uniformizationtransfer} gives rational descent, including the field-degree factor in the Tate invariant. Raynaud replacement handles an abelian quotient with semistable bad reduction when the stated uniformizing model is used. The added toric and lattice ranks, rather than the original algebraic weight ranks, enter \Cref{thm:mixednormalform}. Without such a replacement, the general monodromy target is the space in \eqref{eq:general-obstruction-spaces}; the original toric matrix alone need not describe it. Ramified potential semistability requires descent data and its own normalizations.

Formal-group and CM towers can have Lubin--Tate components. Their coefficient fields, Hodge embeddings, and analytic structures are further data for a comparison with the pairing \eqref{eq:period-Iw-pairing}. Existing formal-group and Lubin--Tate reciprocity laws \cite{TR11,SV26} provide that analytic foundation. The cohomological period pairing for the general $1$-motivic tower is already \Cref{cor:motivicdivisiontower}; a distribution-valued interpolation theorem is a separate comparison with those analytic theories.

\subsection{Equal characteristic}

Let $A$ be a coefficient ring of a global function field, let $F$ be a local field of characteristic $p$, and let $\phi$ be a Drinfeld $A$-module of rank $d$. For a prime $\mathfrak p$ away from the characteristic ideal $\ker(A\to F)$, the $\mathfrak p$-adic Tate module is free of rank $d$. At a nonzero characteristic ideal its \emph{\'etale} Tate-module rank depends on the height and can be smaller than $d$. Thus the coefficient prime and local characteristic are part of any comparison statement.

Existing Coleman and Wiles reciprocity constructions for Drinfeld modules, including \cite{BL09,AE23}, are the appropriate starting point. A precise further problem is to compare their logarithmic derivative or residue pairing with Frobenius-difference period coordinates at a fixed finite place, retaining trace and lattice normalization. Periods and Anderson generating functions at the infinite place do not automatically give such a local finite-place comparison. An equal-characteristic comparison would need its own coefficient ring and the appropriate finite-place Frobenius-difference realization before a residue identity can be matched with \Cref{thm:fullresiduereciprocity}.

\subsection{What the integral matrix calculation says modulo \texorpdfstring{$p$}{p}}

The scalar radical in \Cref{prop:tracedefect} explains the failure of the trace-zero self-pairing to remain perfect modulo $p$ when $p\mid r$. The projective/trace-zero pairing remains perfect after reduction, since its integral matrix in dual bases is the identity. This is a concrete modular representation-theoretic feature of the arithmetic cohomology lattices. A connection with modular local Langlands or a moduli problem would require an additional construction identifying these lattices in that setting.

\appendix
\section{Algebraic checks and normalization record}\label{app:checks}

The proofs use explicit algebra, analytic convergence estimates, and the comparison theorems cited at the relevant steps. The algebraic calculations can be reproduced directly from the
formulas in the paper. For the reduced differentials in
\eqref{eq:reducedcomplex} and \eqref{eq:mixed-reduced}, substitution gives
\[
\delta^1\delta^0(y)
=R(1-\sigma)y-(1-\sigma)Ry=0,
\]
because $R$ commutes with $\sigma$. The cohomology dimensions in
\eqref{eq:dimension-general} and \eqref{eq:mixed-dim-general}
follow from the trace normal forms and rank calculations proved in
\Cref{thm:threeblockcohomology,thm:normalform}.

The block pairing \eqref{eq:blockpair} is obtained by multiplying
the matrices in \eqref{eq:blockA} and \eqref{eq:blockB} and taking
their trace. Multiplication of matrix units gives the identity
pairing matrix for the two bases in \Cref{ex:three}.
For the self trace form in \Cref{prop:tracedefect}, the off-diagonal
matrix units form mutually dual pairs, while the diagonal Gram
matrix is $I_n+J_n$, where $J_n$ is the all-ones matrix. Since
\[
\det(I_n+J_n)=n+1=r,
\]
the full Gram determinant is $\pm r$, and its discriminant ideal
is $(r)$. Formal multiplication of the extension matrices gives
the scalar cocycle identity \eqref{eq:d-cocycle}.

For the mixed setting, the elimination of the acyclic summands
is described in \Cref{lem:threeblock}, and the resulting
cohomology is computed in \Cref{thm:threeblockcohomology}.
The monodromy target and rank-change criterion are established in
\Cref{thm:universalboundary,cor:rankjump}. Substitution of the
corrected cochain into the Herr differential verifies the lifting
identity in \Cref{thm:relativeHerrlift}.
The elliptic point count and the dimensions of the local
conditions are given in \Cref{ex:mixedelliptic} and
\eqref{eq:elliptic-local-table}. Analytic convergence is treated
in \Cref{sec:explicitrobba,thm:explicitprimitive}.

The algebraic identities follow from the displayed formulas
and the explicit arguments in the proofs referenced above. The conventions which enter these checks are collected here for reference:
\begin{equation*}
\begin{gathered}
\log_p p=0,\qquad N\ell_p=-1,\qquad N\varphi=p\varphi N,\\
d^0_{\mathrm{st}}=((1-\varphi),N,\mathrm{quot}),\qquad
d^1_{\mathrm{st}}(a,b,c)=Na-(1-p\varphi)b,\\
d^0_{F,G}=(F-1,G-1),\qquad d^1_{F,G}=(1-G,F-1),\\
\inv_K(\kappa(a)\cup\lambda)=\lambda(a),\qquad
\epsilon=\lambda_0(\operatorname{Frob}_K),\qquad\Tr(\vartheta)=1.
\end{gathered}
\end{equation*}
The unramified arithmetic quotient coordinate is $\epsilon^{-1}\Tr(u)$ in the reduced complex. The factor $\epsilon$ occurs in \eqref{eq:tau} and \eqref{eq:actualhdata}; it does not occur in the normalized arithmetic reciprocity formulas. The integral trace-one choice avoids a division by $n$. The only denominator in the rational trace-zero projection is $r$, whose exact lattice effect is \eqref{eq:trace-defect}. The integral coefficient construction uses the inverse theorem on the integral $\psi=0$ space in \Cref{prop:integralRobba}; the full residue invariant is integral on the cohomology lattices of \Cref{thm:fullresiduereciprocity}.

\section*{Computational and AI assistance}
Python was used for exploratory calculations and exact checks
of explicit algebraic computations and examples. Selected algebraic
results and their proofs were formalized and checked in Lean~4.
OpenAI's ChatGPT assisted with preparing the Python code and
Lean formalizations. The Lean~4 source, recorded verification output, and instructions
for reproducing the nine checks are available in the companion
repository~\cite{MohajerLean4}. The research conception, original mathematical
contributions, proof strategies, and mathematical arguments were
developed by the author.


\begin{thebibliography}{99}
\bibitem{AE23} M. Ala Eddine, \emph{Explicit reciprocity laws for formal Drinfeld modules}, J. Th\'eor. Nombres Bordeaux \textbf{35} (2023), no. 3, 675--695. \href{https://doi.org/10.5802/jtnb.1260}{doi:10.5802/jtnb.1260}.

\bibitem{BL09} F. Bars and I. Longhi, \emph{Coleman's power series and Wiles' reciprocity for rank 1 Drinfeld modules}, J. Number Theory \textbf{129} (2009), no. 4, 789--805. \href{https://doi.org/10.1016/j.jnt.2008.12.005}{doi:10.1016/j.jnt.2008.12.005}.

\bibitem{Ben11} D. Benois, \emph{A generalization of Greenberg's $\mathcal L$-invariant}, Amer. J. Math. \textbf{133} (2011), no. 6, 1573--1632. \href{https://doi.org/10.1353/ajm.2011.0043}{doi:10.1353/ajm.2011.0043}; \href{https://arxiv.org/abs/0906.2857}{arXiv:0906.2857}.

\bibitem{Ber02} L. Berger, \emph{Repr\'esentations $p$-adiques et \'equations diff\'erentielles}, Invent. Math. \textbf{148} (2002), no. 2, 219--284. \href{https://perso.ens-lyon.fr/laurent.berger/articles/article02.pdf}{Author's version}. See the normalization correction in \cite[Appendix B]{Ber08}.

\bibitem{Ber08} L. Berger, \emph{\'Equations diff\'erentielles $p$-adiques et $(\varphi,N)$-modules filtr\'es}, Ast\'erisque \textbf{319} (2008), 13--38. \href{https://perso.ens-lyon.fr/laurent.berger/articles/article10.pdf}{Author's version}.

\bibitem{BerIntro02} L. Berger, \emph{An introduction to the theory of $p$-adic representations}, 2002 preprint. \href{https://arxiv.org/abs/math/0210184}{arXiv:math/0210184}. Section references are to this version.

\bibitem{BWZ26} A. Bertapelle, S. Wang, and H. Zhao, \emph{Log $p$-divisible groups and semi-stable representations}, preprint, version 4, 20 May 2026. \href{https://arxiv.org/abs/2302.11030v4}{arXiv:2302.11030v4}.

\bibitem{BLZ16} A. Besser, D. Loeffler, and S. L. Zerbes, \emph{Finite polynomial cohomology for general varieties}, Ann. Math. Qu\'ebec \textbf{40} (2016), no. 1, 203--220. \href{https://doi.org/10.1007/s40316-015-0041-7}{doi:10.1007/s40316-015-0041-7}.

\bibitem{BK90} S. Bloch and K. Kato, \emph{$L$-functions and Tamagawa numbers of motives}, in \emph{The Grothendieck Festschrift}, Vol. I, Progr. Math. 86, Birkh\"auser, Boston, 1990, 333--400. \href{https://virtualmath1.stanford.edu/~conrad/BSDseminar/refs/BKTamagawa.pdf}{Public scan}.

\bibitem{BC09} O. Brinon and B. Conrad, \emph{CMI Summer School Notes on $p$-adic Hodge Theory}, preliminary version, 2009. \url{https://math.stanford.edu/~conrad/papers/notes.pdf}.

\bibitem{CC99} F. Cherbonnier and P. Colmez, \emph{Th\'eorie d'Iwasawa des repr\'esentations $p$-adiques d'un corps local}, J. Amer. Math. Soc. \textbf{12} (1999), no. 1, 241--268. \href{https://doi.org/10.1090/S0894-0347-99-00281-7}{doi:10.1090/S0894-0347-99-00281-7}; \href{https://webusers.imj-prg.fr/~pierre.colmez/CCjams.pdf}{Author's version}.

\bibitem{CGY25} A. Chitrao, E. Ghate, and S. Yasuda, \emph{Semistable representations as limits of crystalline representations}, Algebra Number Theory \textbf{19} (2025), no. 6, 1049--1097. \href{https://doi.org/10.2140/ant.2025.19.1049}{doi:10.2140/ant.2025.19.1049}; \href{https://arxiv.org/abs/2109.13676}{arXiv:2109.13676}.

\bibitem{CI99} R. Coleman and A. Iovita, \emph{The Frobenius and monodromy operators for curves and abelian varieties}, Duke Math. J. \textbf{97} (1999), no. 1, 171--215. \href{https://doi.org/10.1215/S0012-7094-99-09708-9}{doi:10.1215/S0012-7094-99-09708-9}; \href{https://arxiv.org/abs/math/9701229}{arXiv:math/9701229}.

\bibitem{Col04} P. Colmez, \emph{Fontaine's rings and $p$-adic $L$-functions}, Tsinghua University course notes, autumn 2004. \href{https://www.mathi.uni-heidelberg.de/~ariedel/plfun/colmez2.pdf}{Online course notes}.

\bibitem{Col08} P. Colmez, \emph{Repr\'esentations triangulines de dimension $2$}, Ast\'erisque \textbf{319} (2008), 213--258. \href{https://webusers.imj-prg.fr/~pierre.colmez/Trianguline-ast.pdf}{Author's version}.

\bibitem{FO} J.-M. Fontaine and Y. Ouyang, \emph{Theory of $p$-adic Galois Representations}, book manuscript. \url{https://www.math.u-psud.fr/~fontaine/galoisrep.pdf}.

\bibitem{GZ25} H. Gao and L. Zhao, \emph{Cohomology of $(\varphi,\tau)$-modules}, Adv. Math. \textbf{476} (2025), article 110363. \href{https://doi.org/10.1016/j.aim.2025.110363}{doi:10.1016/j.aim.2025.110363}; \href{https://arxiv.org/abs/2402.02492}{arXiv:2402.02492}.

\bibitem{HS05} D. Harari and T. Szamuely, \emph{Arithmetic duality theorems for $1$-motives}, J. Reine Angew. Math. \textbf{578} (2005), 93--128. \href{https://doi.org/10.1515/crll.2005.2005.578.93}{doi:10.1515/crll.2005.2005.578.93}; \href{https://arxiv.org/abs/math/0304480}{arXiv:math/0304480}. See the authors' \href{https://www.imo.universite-paris-saclay.fr/~david.harari/errata/erratcrelle.pdf}{updated corrigenda}.

\bibitem{Haz70} M. Hazewinkel, \emph{A duality theorem for the cohomology of Lie algebras}, Math. USSR-Sb. \textbf{12} (1970), no. 4, 638--644. \href{https://doi.org/10.1070/SM1970v012n04ABEH000942}{doi:10.1070/SM1970v012n04ABEH000942}.

\bibitem{Her01} L. Herr, \emph{Une approche nouvelle de la dualit\'e locale de Tate}, Math. Ann. \textbf{320} (2001), no. 2, 307--337. \href{https://doi.org/10.1007/PL00004476}{doi:10.1007/PL00004476}.

\bibitem{Liu08} R. Liu, \emph{Cohomology and duality for $(\varphi,\Gamma)$-modules over the Robba ring}, Int. Math. Res. Not. IMRN (2008), article rnm150, 32 pp. \href{https://doi.org/10.1093/imrn/rnm150}{doi:10.1093/imrn/rnm150}.

\bibitem{Milne06} J. S. Milne, \emph{Arithmetic Duality Theorems}, second ed., BookSurge, 2006. \url{https://www.jmilne.org/math/Books/ADTnot.pdf}.

\bibitem{MilneAV} J. S. Milne, \emph{Abelian Varieties}, course notes, version 2.00, 2008. \url{https://www.jmilne.org/math/CourseNotes/AV.pdf}.

\bibitem{MohajerLean4}
M.~Mohajer,
\emph{Selected Lean 4 checks for ``Semistable 1-motives:
monodromy, explicit periods, and local duality''}
[Computer software], 2026.
\url{https://github.com/mohammadreza-mohajer/semistable-1-motives--lean4-checks}.

\bibitem{Nak14} K. Nakamura, \emph{Iwasawa theory of de Rham $(\varphi,\Gamma)$-modules over the Robba ring}, J. Inst. Math. Jussieu \textbf{13} (2014), no. 1, 65--118. \href{https://doi.org/10.1017/S1474748013000078}{doi:10.1017/S1474748013000078}.

\bibitem{PR94} B. Perrin-Riou, \emph{Th\'eorie d'Iwasawa des repr\'esentations $p$-adiques sur un corps local}, with an appendix by J.-M. Fontaine, Invent. Math. \textbf{115} (1994), 81--161.

\bibitem{PR01} B. Perrin-Riou, \emph{Th\'eorie d'Iwasawa des repr\'esentations $p$-adiques semi-stables}, M\'em. Soc. Math. France (N.S.) \textbf{84} (2001), 115 pp. \href{https://doi.org/10.24033/msmf.397}{doi:10.24033/msmf.397}.

\bibitem{Por24} G. Porat, \emph{Locally analytic vector bundles on the Fargues--Fontaine curve}, Algebra Number Theory \textbf{18} (2024), no. 5, 899--946. \href{https://doi.org/10.2140/ant.2024.18.899}{doi:10.2140/ant.2024.18.899}.

\bibitem{Poy25} L. Poyeton, \emph{Locally analytic vectors and rings of periods}, Doc. Math. \textbf{30} (2025), no. 4, 839--886. \href{https://doi.org/10.4171/DM/995}{doi:10.4171/DM/995}.

\bibitem{Ray94} M. Raynaud, \emph{1-motifs et monodromie g\'eom\'etrique}, Ast\'erisque \textbf{223} (1994), 295--319.

\bibitem{Rie13} A. Riedel, \emph{On Perrin-Riou's exponential map and reciprocity laws for $(\varphi,\Gamma)$-modules}, doctoral thesis, Universit\"at Heidelberg, 2013. URN: urn:nbn:de:bsz:16-heidok-149359. \url{https://www.mathi.uni-heidelberg.de/~ariedel/}.

\bibitem{SV26} P. Schneider and O. Venjakob, \emph{Reciprocity Laws for $(\varphi_L,\Gamma_L)$-Modules over Lubin--Tate Extensions}, EMS Press, cover year 2026, published online 9 October 2025, 194 pp. \href{https://doi.org/10.4171/MEMS/24}{doi:10.4171/MEMS/24}.

\bibitem{TR11} F. Tavares Ribeiro, \emph{An explicit formula for the Hilbert symbol of a formal group}, Ann. Inst. Fourier \textbf{61} (2011), no. 1, 261--318. \href{https://doi.org/10.5802/aif.2602}{doi:10.5802/aif.2602}.

\bibitem{WZ24} M. W\"urthen and H. Zhao, \emph{Log $p$-divisible groups associated with log 1-motives}, Canad. J. Math. \textbf{76} (2024), no. 3, 946--983. \href{https://doi.org/10.4153/S0008414X23000287}{doi:10.4153/S0008414X23000287}.

\end{thebibliography}
\end{document}